\documentclass[a4paper, 12pt]{article}
\usepackage{latexsym}
\usepackage{amsmath}
\usepackage{graphics}
\usepackage{amsthm}
\usepackage{amssymb}
\newfont{\Fr}{eufm10}
\newfont{\Sc}{eusm10}
\newfont{\Bb}{msbm10}
\newfont{\Am}{msam10}
\newfont{\am}{msam7}
\numberwithin{equation}{subsection}
\newtheorem{theorem}{Theorem}[section]
\newtheorem{proposition}[theorem]{Proposition}
\newtheorem{lemma}[theorem]{Lemma}
\newtheorem{corollary}[theorem]{Corollary}
\newtheorem{claim}{Claim}
\newtheorem{conjecture}{Conjecture}
\newtheorem{ftheorem}{Theorem}{\bf}{\it}
\newtheorem{fproposition}[ftheorem]{Proposition}{\bf}{\it}
\theoremstyle{definition}
\newtheorem{definition}[theorem]{Definition}
\newtheorem{fdefinition}[ftheorem]{Definition}{\bf}{\rm}

\theoremstyle{remark}
\newtheorem{example}[theorem]{Example}
\newtheorem{remark}[theorem]{Remark}
\newtheorem{fremark}[ftheorem]{Remark}
\newtheorem{definition and corollary}[theorem]{Definition and Corollary}
\newtheorem{fexamples}[ftheorem]{Examples}{\it}{\rm}
\newtheorem{fexample}[ftheorem]{Example}{\it}{\rm}

\newcommand{\C}{{\bf C}}

\newcommand{\Hom}{\mbox{\rm Hom}}

\newcommand{\Rep}{\mbox{\rm Rep}}

\newcommand{\Ind}{\mbox{\rm Ind}}

\newcommand{\id}{\mbox{\rm id}}

\newcommand{\Spec}{\mbox{\rm Spec}}

\newcommand{\Pic}{\mbox{\rm Pic}}

\newcommand{\MID}{\! \! \mid}

\newcommand{\h}{\mathfrac{\h}}

\newcommand{\Q}{\mathbb{Q}}
\newcommand{\R}{\mathbb{R}}
\newcommand{\Z}{\mathbb{Z}}

\begin{document}
\title{Equivariant vector bundles on group completions}
\author{KATO, Syu\footnote{The author is supported by JSPS Research Fellowship for Young Scientists} \footnote{Graduate School of Mathematical Sciences, University of Tokyo, 3-8-1 Komaba Meguro 153-8914 Tokyo Japan}\footnote{E-mail: syuchan@ms.u-tokyo.ac.jp} }
\maketitle

\abstract{
In this paper, we describe the category of bi-equivariant vector bundles on a bi-equivariant smooth (partial) compactification of a connected reductive algebraic group with normal crossing boundary divisors. Our result is a generalization of the description of the category of equivariant vector bundles on toric varieties established by A.A. Klyachko [Math. USSR. Izvestiya {\bf 35} No.2 (1990)]. As an application, we prove splitting of equivariant vector bundles of low rank on the wonderful compactification of an adjoint simple group in the sense of C. De Concini and C. Procesi [Lecture Notes in Math. {\bf 996} (1983)]. Moreover, we present an answer to a problem raised by B. Kostant in the case of complex groups.
}
\tableofcontents
\section*{Introduction}
Let $G$ be a connected reductive algebraic group over an algebraically closed field $k$ of characteristic zero. Then we have a $G \times G$-action on $G$ defined as $(g _1, g _2). x = g _1. x. g _2 ^{-1}$. Let $X$ be a $G \times G$-equivariant smooth partial compactification of $G$ with normally crossing boundary divisors. (Here, a nonsingular toric variety (cf. Oda \cite{Oda}) and the wonderful compactification of an adjoint semisimple group in the sense of De Concini-Procesi \cite{DP} satisfy the condition.) By the general theory (Theorem \ref{fundamental}), such a partial compactification is described by a fan $\Sigma$. 

A ($G \times G$-) equivariant vector bundle is a vector bundle ${\mathcal E}$ on $X$ with an action of $G \times G$ which is linear on each fiber and makes the
following diagram commutative:
$$
\begin{matrix}
V \left( {\mathcal E} \right) & \stackrel{g ^*}{\rightarrow} & V
\left( {\mathcal E} \right) & \\
\downarrow & & \downarrow & \forall g = ( g _1, g _2 ) \in G \times G.\\
X & \stackrel{g ^*}{\rightarrow} & X &
\end{matrix}
$$
Here, $V \left( {\mathcal E} \right)$ denotes the total space of
${\mathcal E}$.

Let $EV \left( X \right)$ be the category of $G \times G$-equivariant vector bundles and morphisms of $G \times G$-equivariant coherent ${\mathcal O} _X$-modules. (Our actual setting is to handle a finite cover $\tilde{G}$ of $G$ instead of $G$. Here we restrict the situation for simplicity.) Our goal (cf. Theorem \ref{tinymain}) is to describe this category explicitly via linear algebra data.

Let $e$ be the point of $X$ corresponding to the identity element of
$G$. Choose a maximal torus $T$ of $G$. An equivariant vector bundle ${\mathcal E}$ is completely determined by the
following data:

\begin{itemize}
\item ${\mathcal E} \otimes k ( e )$: the fiber at $e$;
\item A family of subspaces
$$F ^{\tau} \left( n, {\mathcal E} \right) := \{ v \in {\mathcal E}
\otimes k ( e ); \exists \lim _{t \rightarrow \infty} t ^{n}
( 1 \times \tau ( t ) ) v \in V \left( {\mathcal E} \right) \} \subset {\mathcal E}
\otimes k ( e )$$
for every one-parameter subgroup $\tau : \mathbb{G} _m \rightarrow T$ and every $n \in \Z$. In other words, the asymptotic behavior of elements of ${\mathcal E} \otimes k ( e )$
with respect to all one-parameter subgroups.
\end{itemize}

Thus, the problem is to find a sufficient condition for the existence
of an equivariant vector bundle which corresponds to a given asymptotic datum.

Klyachko \cite{Kl} described the category of equivariant
vector bundles on toric varieties by using the above asymptotic data $\left( {\mathcal E} \otimes k ( e ), \{ F ^{\tau} \left( n, {\mathcal E} \right) \} _{\tau} \right)$. His description is based on the fact that there are no ``local" obstruction for the existence of an equivariant vector bundle with a given asymptotic behavior. (For another description in this case, see Kaneyama \cite{Kane1, Kane2}.)

In this paper, we extend such a description to the compactification of an arbitrary reductive algebraic group by introducing a new constraint which we call the transversality condition.

Let $\Sigma ( 1 )$ be the set of one-skeletons of the fan $\Sigma$ which is identified with a set of one-parameter subgroups of $T$. We put $P ^{\tau} := \{g \in G; \exists \lim _{t \rightarrow 0} \tau ( t ) g \tau ( t ) ^{- 1} \in G \}$. Put $\mathfrak{g} := \mathrm{Lie} G$. For each root $\alpha$, let us denote the $\alpha$-root space of $\mathfrak g$ by $\mathfrak g _{\alpha}$.

\begin{fdefinition}
We define the category ${\mathfrak C} \left( X \right) ^+$ by the following:
\begin{itemize}
\item ${\bf (Objects)}$ Pairs $\left( V, \{ F ^ {\tau} \left( \bullet \right) \} _{\tau \in \Sigma ( 1 )} \right)$ consisting of a $G$-module $V$ and exhausting $\Z$-indexed decreasing $P ^{\tau}$-filtrations $F ^ {\tau} \left( \bullet \right)$ of $V$ such that:
\begin{itemize}
\item For each $\sigma \in \Sigma$, there exists a basis $B ^{\sigma}$ of $V$ which spans every $F ^{\tau} ( n )$ when $\tau \in \sigma ( 1 )$.
\end{itemize}
\item ${\bf (Morphisms)}$ Let $\left( V _1, \{ F ^ {\tau} _1 \left( \bullet \right) \} _{\tau} \right), \left( V _2, \{ F ^ {\tau} _ 2\left( \bullet \right) \} _{\tau} \right) \in {\bf Ob} {\mathfrak C} \left( X \right) ^+$. Then, we define
\begin{align*}
& \Hom _{{\mathfrak C} \left( X \right) ^+} \left( \left( V _1, \{ F ^ {\tau} _1 \left( \bullet \right) \} _{\tau \in \Sigma ( 1 )} \right), \left( V _2, \{ F ^ {\tau} _ 2\left( \bullet \right) \} _{\tau \in \Sigma ( 1 )} \right) \right) \\
& := \{ f \in \Hom _G \left( V _1, V _2 \right) ; f \left( F ^{\tau} _1 ( n ) \right) \subset F ^{\tau} _2 ( n ) \text{ for every }n \in \Z \text{ and every } \tau \in \Sigma ( 1 ) \}.
\end{align*}
\end{itemize}
We also define a full-subcategory ${\mathfrak C} \left( X \right)$ of ${\mathfrak C} \left( X \right) ^+$ by the following:
\begin{itemize}
\item $\left( V, \{ F ^ {\tau} \left( \bullet \right) \} _{\tau \in \Sigma ( 1 )} \right) \in {\bf Ob} {\mathfrak C} \left( X \right) ^+$ is in ${\bf Ob} {\mathfrak C} \left( X \right)$ if, and only if, $F ^{\tau} \left( \bullet \right)$ satisfies the following ($\tau$-) transversality condition for each $\tau \in \Sigma ( 1 )$.
\end{itemize}
$\text{\bf (Transversality Condition)}$ For every $n \in \Z$ and each root $\alpha$ of $\mathfrak{g}$ such that $\left< \tau, \alpha \right> < 0$, we have $\mathfrak g _{\alpha} F ^{\tau} \left( n \right) \subset F ^{\tau} \left( n + \left< \tau, \alpha \right> \right)$.
\end{fdefinition}

\begin{fexamples}\label{iexample}
\begin{enumerate}
\item If $X = G$, then $\mathfrak C \left( X \right) = \mathfrak C \left( X \right) ^+ = \mathrm{Rep} G$.
\item {\bf (Case $G = \mathop{PGL} _2$ and $X = \mathbb P ^3$)} Fix a triangular decomposition $\mathfrak{sl} _2 = \mathfrak n^- \oplus \mathfrak b$. Then, an object of $\mathfrak C \left( X \right) ^+$ is a pair consisting of a $\mathop{PGL} _2$-module $V$ and its $\mathfrak b$-stable decreasing filtration $F ( \bullet )$. Transversality condition says $\mathfrak n ^- F ( n ) \subset F ( n - 1 )$ for every $n$.
\end{enumerate}
\end{fexamples}

\begin{fproposition}[cf. Proposition \ref{naive functor} + \ref{last}]\label{nf}
The map
$$\Xi :{\bf Ob} EV \left( X \right)\ni {\mathcal E} \mapsto \left( {\mathcal E} \otimes k ( e ), \{ F ^{\tau} \left( n, \mathcal E \right) \} _{\tau \in \Sigma ( 1 )} \right) \in {\bf Ob} {\mathfrak C} \left( X \right) ^+$$
is injective. Moreover, it gives a faithful covariant functor $EV \left( X \right) \rightarrow {\mathfrak C} \left( X \right) ^+$ which we denote by the same letter $\Xi$.
\end{fproposition}

\begin{fremark}
Let $D _{\tau}$ be the $G \times G$-equivariant divisor corresponding to $\tau \in \Sigma ( 1 )$. Then, tensoring $\mathcal O _X ( n D _{\tau} )$ at the LHS corresponds to shifting the degree of $F ^{\tau}$ by $- n$ at the RHS.
\end{fremark}

The idea in the proof of Proposition \ref{nf} is to measure the asymptotic behavior of an equivariant vector bundle by bounding it with a certain family of ``standard" vector bundles instead of using the asymptotic behavior of elements on the identity fiber.

Next, our key result is as follows:

\begin{fproposition}[cf. Proposition \ref{induced action}]\label{ia}
The image of the map
$$\Xi : {\bf Ob} EV \left( X \right) \rightarrow {\bf Ob} {\mathfrak C} \left( X \right) ^+$$
is contained in ${\bf Ob} {\mathfrak C} \left( X \right)$.
\end{fproposition}

In the above, there are two choices of families of ``standard" vector bundles. It turns out that each of them yields different constraints about filtrations. Proposition \ref{ia} follows from the comparison of these constraints by means of the theory of Tannakian categories (cf. \cite{DMOS}).

Now we state our main result.

\begin{ftheorem}[cf. Theorem \ref{Main Theorem}]\label{tinymain}
We have a category equivalence
$$\Xi : EV \left( X \right) \stackrel{\cong}{\rightarrow} {\mathfrak C} \left( X \right).$$
\end{ftheorem}

For the proof of Theorem \ref{tinymain}, we construct the inverse functor $\Phi$ of $\Xi$. The construction of $\Phi$ depends on the fact that the boundary behavior of equivariant vector bundles is completely controlled by their restriction to certain toric varieties. Thus, we devote ourselves to check that the reconstruction process prescribed by $\Phi$ is compatible with the action of the ``unipotent part".

If $G$ is commutative, the left and the right multiplications are essentially the same. In this case, the transversality condition is a void condition. However, we cannot identify the left and the right multiplications if $G$ is non-commutative. To establish Theorem \ref{tinymain} in this setting, we need to control the $G \times G$-action via a single group $G$, which is the stabilizer of the $G \times G$-action on $G$ at the identity element of $G$. This is why the transversality condition appears. Moreover, it imposes a strong restriction on the existence of equivariant vector bundles. A typical example of this kind of phenomena is the following.

\begin{fexample} In the same setting as in Examples \ref{iexample} 2), simple objects are classified as follows: Let $V _n$ be an irreducible $\mathop{PGL} _2$-module of dimension $2n + 1$. Then, every irreducible equivariant vector bundle of rank $2n + 1$ has $V _n$ as its identity fiber. From this fact, one can prove that irreducible equivariant vector bundles are classified by Dynkin quivers of type $A _{2 n + 1}$ up to line bundle twist. In particular, there are only finitely many (namely $2 ^{2n}$) possibilities.
\end{fexample}

For more details, see Example \ref{example}. In general, a simple equivariant vector bundle may have nontrivial equivariant deformations even in the case of the wonderful compactification.

Kostant has raised the question of the
existence of a canonical extension of an equivariant vector bundle on
a symmetric space to its wonderful compactification in order to deduce representation theoretic
data from asymptotic expansions of matrix coefficients (see \S \ref{canonical extension}). Our description gives the following answer to
his problem in the case of complex groups.

\begin{ftheorem}[= Theorem \ref{Answer to Kostant}]
Let $G$ be a semisimple adjoint group, and let $X$ be its wonderful compactification. For every $G$-module $V$, there exists a unique $G \times G$-equivariant vector bundle ${\mathcal E} _V$ on $X$ which satisfies the following
properties:
\begin{enumerate}
\item ${\mathcal E} _V \otimes k ( e ) \cong V$ as a
$G$-module;
\item For every $v \in {\mathcal E} _V \otimes k ( e )$ and every one-parameter subgroup $\tau : \mathbb{G} _m \rightarrow G$, there exist a
limit value $\lim _{t
\rightarrow 0} ( \tau ( t ) \times \tau ( t ) ^{- 1} ) v$ in $V \left( {\mathcal E} _V \right)$;
\item Every $G \times G$-equivariant vector bundle ${\mathcal E}$ with the above two properties can be $G \times G$-equivariantly embedded into ${\mathcal E} _V$.
\end{enumerate}
\end{ftheorem}

Though the above Kostant problem is known by some experts, there is no existing literature. Thus, we also present the whole picture of his problem in \S \ref{canonical extension}. (The author learned about this problem from Prof. Brion, Prof. Kostant, and Prof. Uzawa. He wants to express gratitude to them.)

As a bonus of our description, we have the following result.

\begin{ftheorem}[= Corollary \ref{split via rank-statement}]\label{srpre}
Let $G$ be an adjoint simple group and let $X$ be its wonderful compactification. Then, every $G \times G$-equivariant vector bundle of rank less than or equal to $r = \mathrm{rk} G$ splits into a direct sum of
line bundles.
\end{ftheorem}
This kind of material is treated in \S \ref{compare}.
For a classical (simple) group, the wonderful compactification is
obtained by successive blowing-ups of a partial flag variety of an overgroup
of $G$ (cf. Brion \cite{Br2}). On projective spaces or Grassmannians, we have
a splitting criterion of vector bundles in terms of the vanishing of
intermediate cohomologies (see Horrocks \cite{H}, Ottaviani \cite{O},
or Arrondo and Gra\~na \cite{AG}). Thus, if we have an analogous result of
\cite{Kl} 1.2.1 (a vector bundle is equivariant if its infinitesimal
deformation is zero), we may find a splitting criterion of vector bundles using our result. There are plenty of vector bundles with nontrivial infinitesimal deformations, but a
direct sum of line bundles on the wonderful compactification has no
nontrivial infinitesimal deformation (cf. Tchoudjem \cite{T}, \cite{T2}, \cite{T3} and K \cite{Ka}). However, there exists a line bundle with non-vanishing intermediate cohomologies in this case. Hence, a naive reformulation
of the splitting criterion of vector bundles in terms of cohomology vanishing has a counter-example in this case.

The organization of this paper is as follows. In \S
\ref{statement}, we fix our notation and introduce objects which we concern. In
particular, we assume the notation and terminology introduced in
\S \ref{statement} in the whole
paper unless stated otherwise. In \S \ref{setting up}, we present some fundamental results which is needed to formulate our main theorem. More precisely, in \S \ref{prelim}, we devote ourselves to the preparation of the working ground. In \S \ref{reconst}, we construct a functor $\Xi$ in its full generality (Proposition \ref{naive functor}). In \S \ref{toric slice}, we develop a technique to discuss the boundary behavior of equivariant vector bundles and prove that the ``image" of $\Xi$ is contained in ${\mathfrak C} \left( X \right)$ (Proposition \ref{induced action}). In \S \ref{MT}, we state and prove Theorem \ref{Main Theorem}, which is our main result. Finally, we deal with some consequences of Theorem \ref{Main Theorem} in \S \ref{sectioncor}.

This paper is the main body of the author's Doctoral Dissertation at University of Tokyo.
\section{Notation and Terminology}\label{statement}
\subsection{Notation on algebraic groups}\label{NoteAlg}

The general reference for the material in this subsection is Springer's book \cite{Sp}.

Let $G$ be a connected reductive group of rank $r$ over an
algebraically closed field $k$ of characteristic zero. There exists a connected
finite cover $\tilde{G}$ of $G$ which is isomorphic to the direct
product of a torus $\tilde{T} _0$
and a simply connected semisimple algebraic group $\tilde{G} _s$. We put
$\tilde{Z} \left( G \right) := \ker [ \tilde{G} \rightarrow G ]$. For a group $H$, we denote its center by $Z ( H )$. We have $\tilde{Z} \left( G \right) \subset Z ( \tilde{G} )$. We put $G _{ad} := G / Z ( G )$, $T _{ad} := T / Z ( G )$, and $\tilde{T} _s := \tilde{G} _s \cap \tilde{T}$.

Let $B$ and $B ^-$ be (mutually opposite) Borel subgroups of $G$, with a unique common maximal torus $T$;
let $N _G \left( T \right)$ be its normalizer, and put $W = N _G \left(
T \right) / T$ the corresponding Weyl group. Let $U$ and $U ^-$ be the unipotent radicals of $B$ and $B ^-$, respectively. We
denote by $\tilde{B}$, $\tilde{T}$... the preimages of $B$, $T$... in
$\tilde{G}$.

For a torus $S$, we denote the weight lattice $\Hom \left( S,
\mathbb{G} _m \right)$ of $S$ by $X ^* \left( S \right)$ and the coweight lattice $\Hom
\left( \mathbb{G} _m, S \right)$ of $S$ by $X _* \left(
S \right)$. We put $X ^* ( S ) _{\mathbb R} := X ^* ( S ) \otimes _{\mathbb Z} \mathbb R$. We regard $X ^* \left( T \right)$ as a subset of $X ^* ( \tilde{T} )$. We have a natural $\Z$-bilinear pairing
$$\left< , \right> : X _* \left( T \right) \times X ^* ( \tilde{T} ) \to \Q.$$

Let $\triangle \subset X ^* \left( T \right)$ be the root system of $( G, T
)$ and let $\triangle ^+$ be its subset of positive roots defined by $B$. Let $\mathfrak{g}$, $\mathfrak{b}$, $\mathfrak{t}$... be the Lie algebras of $G$,
$B$, $T$... For each $\alpha \in \triangle$, we fix a root vector $e _{\alpha} \in \mathfrak{g}$. Let $\Pi = \{ \alpha _1, \ldots, \alpha _{\ell}
\} \subset \triangle ^+$ be the set of simple roots. We denote by $\Pi ^{\vee}= \{ \alpha _1 ^{\vee}, \ldots, \alpha _{\ell} ^{\vee} \} \subset X _* ( \tilde{T} _s )$ the set of simple coroots, which we consider as a subset of $X _* ( \tilde{T} )$ via natural inclusion $X_* ( \tilde{T} _s ) \subset X _* ( \tilde{T} )$. The set of fundamental coweights $\{ \omega _1 ^{\vee}, \omega _2
^{\vee}, \ldots, \omega _{\ell} ^{\vee} \}$ is defined as the set of $\Z$-linear
forms on $X ^* \left( T \right)$ such that $\omega _i ^{\vee} \left( \alpha _j
\right) = \delta _{i, j}$ ($1 \le i, j \le \ell$). Every $\omega ^{\vee} _i$ defines an element of $X _* \left( T _{ad} \right)$,
which we also denote by $\omega _i ^{\vee}$.

Let $\tilde{Z} ( G ) ^{\vee}$ be the character group of $\tilde{Z} ( G )$ and let $h$ be its order. Since $\tilde{Z} ( G )$ is contained in $Z ( \tilde{G} )$, we have natural surjection $X ^* ( \tilde{T} ) \rightarrow \tilde{Z} ( G ) ^{\vee}$. For each $\lambda \in X ^* ( \tilde{T} )$, we denote its image in $\tilde{Z} ( G ) ^{\vee}$ by $\bar{\lambda}$.
$\lambda \in X ^* ( \tilde{T} )$ is called a dominant weight if, and only if, $\left< \alpha ^{\vee} _i, \lambda \right> \ge 0$ for every $1 \le i \le \ell$. For a dominant weight $\lambda$, we denote by $V_{\lambda}$ the irreducible rational representation of $\tilde{G}$ with highest weight $\lambda$. Let $v _{w _0 \lambda}$ be a lowest weight vector of $V _{\lambda}$, where $w _0$ is the longest element of $W$.

For every $\tau \in X _* \left( T \right)$, we define the following three Lie subalgebras of $\mathfrak{g}$:
$$\mathfrak{l} ^{\tau} := \mathfrak{t} \oplus \bigoplus _{\alpha \in \triangle ; \left<
\tau, \alpha \right> = 0} k e _{\alpha}, \mathfrak{u} ^{\tau} _+ := \bigoplus _{\alpha \in \triangle ; \left<
\tau, \alpha \right> > 0} k e _{\alpha}, \text{ and } \mathfrak{p}
^{\tau} := \mathfrak{l} ^{\tau} \oplus \mathfrak{u} ^{\tau} _+.$$
We also put $\mathfrak{u} ^{\tau} _- := \mathfrak{u} ^{- \tau} _+$. We denote by $L ^{\tau}, U ^{\tau} _+...$ the corresponding subgroups of
$G$. $P ^{\tau}$ is a parabolic subgroup of $G$. $P ^{\tau} = L
^{\tau} U ^{\tau} _+$ is the Levi decomposition of $P ^{\tau}$ such that
$T \subset L ^{\tau}$.

For a group $H$, we denote the diagonal embedding $H \hookrightarrow H \times H$ by $\triangle ^d$. Also, we denote by $V ^{H}$ the space of $H$-fixed vectors of $V$ for a $H$-module $V$.

\subsection{Notation on partial compactifications}\label{NoteVar}
Here we assume that readers are familiar with the standard material in the theory of toric varieties (cf. Oda \cite{Oda} Chapter 1). The contents in this subsection are found in the papers by De
Concini-Procesi \cite{DP2}, Uzawa \cite{U}, and Knop's survey \cite{Kn}.

For every fan $\Sigma$ in $X _* ( T ) _{\mathbb R}$, we denote by $\Sigma ( n )$ the set of
$n$-dimensional cones of $\Sigma$. By abuse of notation, we denote the
fan consisting of a cone $\sigma$ of $\Sigma$ and its faces by the
same letter $\sigma$. Moreover, we denote the integral generator of a
one-dimensional cone $\tau$ by the same letter. (Thus, $\tau$ defines a one-parameter subgroup of $T$ when $\Sigma$ is a fan of $X _* \left( T \right) _{\mathbb R}$.)

We define a fan $\Sigma _0$ of $X _* \left( T _{ad} \right) _{\mathbb R}$ by
$$\Sigma _0 := \left\{ \Sigma _{i \in S} \R _{\ge 0} \omega _{i} ^{\vee} ; S \subset \{1, 2, \ldots, \ell \} \right\}.$$

$\Sigma _0$ is the fan consisting of the dominant Weyl co-chamber and its faces.

Let $X _0$ be the wonderful compactification of $\left( G _{ad}
\times G _{ad} \right) / \triangle ^d \left( G _{ad} \right)$
in the sense of De Concini-Procesi \cite{DP}. We consider $X _0$
as a $G \times G$-equivariant compactification of $G _{ad}$ via
the quotient map $G \times G \rightarrow G _{ad} \times G _{ad}$.

\begin{definition}[Regular embeddings]
A regular embedding $X$ of $G$ is a $G \times G$-equivariant smooth partial compactification of $X$ with the following conditions:
\begin{itemize}
\item[1.] $X \backslash G$ is a union of normal crossing divisors $D _1, \ldots, D _{p}$;
\item[2.] Each $D _i$ is smooth and is the closure of a single $G \times G$-orbit;
\item[3.] Every $G$-orbit closure in $X$ is a certain intersection of $D _1, \ldots, D _p$;
\item[4.] For each $x \in X$, the total space of the normal bundle of $\overline{( G \times G ) x}$ in $X$ contains a dense $G \times G$-orbit.
\end{itemize}
\end{definition}

\begin{theorem}[Uzawa \cite{U} 3.5. See also \cite{DP2}]\label{fundamental}
Let $\Sigma$ be a fan of $X _* \left( T \right) _{\mathbb R}$ such that the following two conditions hold:
\begin{enumerate}
\item The natural quotient map $q : X _* \left( T \right)
\rightarrow X _* \left( T _{ad} \right)$ yields a
morphism $\Sigma \rightarrow \Sigma _0$ of fans;
\item The toric variety $T \left( \Sigma \right)$ corresponding to
$\Sigma$ is nonsingular.
\end{enumerate}
Then, there exists a unique regular embedding $X \left( \Sigma \right)$ of $G$ such that:
\begin{enumerate}
\item[a.] We have an embedding $T ( \Sigma ) \subset X ( \Sigma )$;
\item[b.] Each $G \times G$-orbit in $X ( \Sigma )$ intersects with a unique $T$-orbit in $T ( \Sigma )$;
\item[c.] There exists a dominant $G \times G$-equivariant morphism $\pi : X \left( \Sigma \right) \rightarrow X _0$.
\end{enumerate}
Its converse is also true. We denote by $\overline{\mathbf O} _{\sigma}$ the closure of a $G \times G$-orbit $\mathbf O _{\sigma}$ corresponding to $\sigma \in \Sigma$. If $\tau \in \Sigma ( 1 )$, we also denote $\overline{\mathbf O} _{\tau}$ by $D _{\tau}$. (We regard it as a prime divisor.)
\end{theorem}

Theorem \ref{fundamental} is more or less a consequence of the results in the references we list above and known by experts. However, since the author does not know an appropriate reference, we provide a proof. In the proof, we need the following modification of a theorem of Strickland, Theorem \ref{structure theorem} which is obtained by pulling back the original version ($X \left( \Sigma \right) = X _0$ case) via $\pi$.

\begin{theorem}[Local structure theorem \cite{Str} cf. \cite{Br} 1.1 or \cite{BLV}]\label{structure theorem}
Under the same settings as in Theorem \ref{fundamental}, the map
$$U ^- \times T \left( \Sigma \right) \times U \longrightarrow X \left( \Sigma \right), \quad (g, x, h) \mapsto (g, h).x,$$
is an open embedding.
\end{theorem}

\begin{proof}[Proof of Theorem \ref{fundamental}]
We write the images of $B$ and $B ^-$ in $G _{ad}$ by $B _{ad}$ and $B
_{ad} ^-$, respectively. Then $B _{ad} B ^- _{ad} \subset G _{ad}$ is an open dense subset. Hence, $G _{ad} \cong G _{ad} \times G _{ad} / \triangle ^d ( G _{ad} )$ and its compactification is a spherical variety. For a reductive group $K$ and its (mutually opposite) Borel subgroups $B _K$ and $B _K ^-$ with a common torus $T _K$, we define $\Lambda _K$ as follows:
$$\Lambda _K := \left( k \left( ( K \times K ) / \triangle ^d ( K ) \right) ^{( B _K \times B ^- _K)} - \{ 0 \} \right) / k ^*.$$
Here superscript $(B _K \times B ^- _K)$ means the eigenpart with respect to the $B _K \times B ^- _K$-action. Then, we have $\Hom \left( \Lambda _K, \Z \right) \cong X _* \left( T _K \right)$. As is described in Brion \cite{Br} (2.2 Remarques i)), the wonderful compactification corresponds to the colored fan $( \Sigma _0, \emptyset )$ (= usual fan in $X _* \left( T _{ad} \right)$). Then, by Knop \cite{Kn} 3.3 and 4.1, we have a $G \times G$-equivariant (partial) compactification $X ^{\prime} \left( \Sigma \right)$ of $(G \times G) / \triangle ^d ( G )$ with a $G \times G$-equivariant dominant map
$$X ^{\prime} \left( \Sigma \right) \rightarrow X _0$$
for each fan $\Sigma$ satisfying 1). Here each cone of $\Sigma$ corresponds to a unique $G \times G$-orbit of $X ^{\prime} \left( \Sigma \right)$ by \cite{Kn} 3.2 and 3.3. Since $\Lambda _G \cong \Lambda _T$, the closure of $T \subset X ^{\prime} \left( \Sigma \right)$ contains $T \left( \Sigma \right)$ by \cite{Kn} 4.1. Since the smoothness is an open condition, Theorem \ref{structure
theorem} and the equality $(G \times G). T \left( \Sigma \right) = X ^{\prime} \left( \Sigma \right)$ asserts that $X ^{\prime} \left( \Sigma \right)$ is smooth if, and only if, $T \left( \Sigma \right)$ is smooth. Therefore, setting $X \left( \Sigma \right) := X ^{\prime} \left( \Sigma \right)$ completes the proof of Theorem \ref{fundamental}.
\end{proof}

From now on, we always assume the assumptions of Theorem \ref{fundamental}. For simplicity, we may write $X$ instead of $X \left( \Sigma \right)$.

Let $\iota _{\sigma}$ be the inclusion $\overline{\mathbf O} _{\sigma} \hookrightarrow X \left( \Sigma \right)$ corresponding to $\sigma \in \Sigma$. We define $\mathfrak{u} ^{\sigma} _+ := \sum _{\tau \in \sigma ( 1 )} \mathfrak{u} ^{\tau} _+ \subset \mathfrak{g}$. Similarly, we write $\mathfrak{l} ^{\sigma}$ for $\cap _{\tau \in \sigma ( 1 )} \mathfrak{l} ^{\tau}$. Let $U ^{\sigma} _+, U ^{\sigma} _-$, and $L ^{\sigma}$ be the algebraic subgroups of $G$ corresponding to $\mathfrak{u} ^{\sigma} _+$, $\mathfrak{u} ^{- \sigma} _+$, and $\mathfrak{l} ^{\sigma}$, respectively. Denote $L ^{\sigma} U ^{\sigma} _+$ by $P ^{\sigma}$ (this is a parabolic subgroup of $G$).

Since $X \left( \Sigma \right)$ dominates $X _0$, the same arguments as in \cite{DP} \S 5.2 assert the existence of a $G \times G$-equivariant fibration $\pi _{\sigma} : \overline{\mathbf O} _{\sigma} \rightarrow G / P ^{\sigma} \times G / P ^{- \sigma}$. We denote the point of $X \left( \Sigma \right)$ corresponding to the
identity element of $G$ by $e$. We have $T \left( \Sigma \right) \subset \overline{T .e} \subset X$. Hence, $x _{\sigma} := \lim _{t \rightarrow \infty} (1 \times \Pi _{\tau \in \sigma ( 1 )} \tau ( t )) e$ exists in $T \left( \Sigma \right)$ for each $\sigma \in \Sigma$. Let $\mathbb{G} _{m} ^{\tau}$ be the image of $\mathbb{G} _m$ via $( 1 \times \tau ^{- 1}) : \mathbb{G} _m \rightarrow T \times T$ for each $\tau \in \Sigma ( 1 )$. We denote by $G ^{\sigma}$ the stabilizer of the $G \times G$-action at $x _{\sigma}$. Then we have
$$G ^{\sigma} = \triangle ^d ( L ^{\sigma} ) ( \prod _{\tau \in \sigma ( 1 )} \mathbb{G} ^{\tau} _{m}) ( U ^{\sigma} _+ \times U ^{\sigma} _- ) = \triangle ^d ( L ^{\sigma} ) ( \prod _{\tau \in \sigma ( 1 )} \mathrm{Im} \tau \times 1 ) ( U ^{\sigma} _+ \times U ^{\sigma} _- ).$$
This is another consequence of the fact that $X \left( \Sigma \right)$ dominates $X _0$.

For simplicity, we denote $\otimes _{{\mathcal O} _X}$ by $\otimes _X$, or even by $\otimes$ when there is no risk of confusion.

\subsection{The category ${\mathfrak C} \left( \Sigma \right)$}\label{NoteCat}
In order to formulate our main theorem, we introduce some notation and a category ${\mathcal C} \left( \Sigma \right) _c$ which contains the category $\mathcal C \left( X \right) = \mathcal C \left( X ( \Sigma ) \right)$ in the introduction as a fullsubcategory. This enhancement is necessary in order to handle all line bundles on $X$ in the main theorem (cf. Theorem \ref{Steinberg Theorem}).

We denote the universal enveloping algebra of a Lie algebra $\mathfrak{a}$ by $U \left( \mathfrak{a} \right)$. For each $\tau \in \Sigma ( 1 )$ and every $n \in \Z$, we define
$$U \left( \mathfrak{g} \right) _n ^{\tau} := \{ X \in U \left( \mathfrak{g} \right) ; \tau ( t ) X \tau ( t ) ^{- 1} = t ^n X \text{ for every }t \in \mathbb{G} _m ( k ) \cong k ^{\times} \}.$$
We have $U \left( \mathfrak{g} \right) = \oplus _{n \in \Z} U \left( \mathfrak{g} \right) _n ^{\tau}$. For each $\tau$-stable
subalgebra $\mathfrak{f}$ of $\mathfrak{g}$, we put $U \left(
\mathfrak{f} \right) ^{\tau} _n := U \left(
\mathfrak{f} \right) \cap U \left(
\mathfrak{g} \right) ^{\tau} _n$. Since the $\tau$-action on
$\mathfrak{g}$ is semi-simple, we have $U \left(
\mathfrak{f} \right) = \oplus _{n \in \Z} U \left(
\mathfrak{f} \right) ^{\tau} _n$. Every Lie algebra defined in \S \ref{NoteAlg} satisfies this property.

\begin{remark}
We have $U \left( \mbox{\Fr l} ^{\tau} \right) = U \left( \mbox{\Fr p} ^{\tau} \right) ^{\tau} _0$ and $U \left( \mbox{\Fr u} ^{\tau} _+ \right)  ^{\tau} _{- n} = U \left( \mbox{\Fr u} ^{\tau} _- \right) ^{\tau} _{n} = 0$ for all $n > 0$.
\end{remark}

Now we introduce our linear algebra data, the category ${\mathfrak C} \left( \Sigma \right)$. The most remarkable property of ${\mathfrak C} \left( \Sigma \right)$ is the ($\tau$-) transversality condition.

\begin{definition}\label{adm}
Let $\tau \in \Sigma ( 1 )$. Let $V$ be a $\tilde{G} \times \tilde{Z} \left( G \right)$-module. A $\tau$-standard filtration $F \left( \bullet \right)$ is a decreasing
filtration of $V$ indexed by $\Z$ such that:
\begin{enumerate}
\item For each $n \in \Z$, $F \left( n \right)$ is a $\tilde{P} ^{\tau} \times \tilde{Z} \left( G \right)$-module via restriction $\tilde{P} ^{\tau} \subset \tilde{G}$;
\item $F \left( - n \right) = V$ and $F \left( n \right) = \{ 0 \}$ for $n >> 0$.
\end{enumerate}
A $\tau$-standard filtration $F \left( \bullet \right)$ is called an $\tau$-transversal filtration if, and only if, the following condition holds:
\begin{itemize}
\item $\text{{\bf (Transversality condition)}}$ For every $n, m \in \Z$, we have$$U ( \mathfrak{u} ^{\tau}
_- ) ^{\tau} _m F \left( n \right) \subset F \left( n + m \right).$$
\end{itemize}
\end{definition}

\begin{definition}\label{distributive lattice}
Let $V$ be a finite dimensional vector space. A family of linear
subspaces $\{U _{\lambda}\} _{\lambda \in \Lambda}$ forms a
distributive lattice if, and only if, there exists a basis $B$ of $V$ such that $B \cap U _{\lambda}$ is a basis of $U _{\lambda}$ for every $\lambda \in \Lambda$.
\end{definition}

For a category ${\mathcal Z}$, we denote the class of objects of ${\mathcal Z}$ by
${\bf Ob} {\mathcal Z}$. Also, for every ${\mathcal X},
{\mathcal Y} \in {\bf Ob} {\mathcal Z}$, we denote the set of morphisms of ${\mathcal X}$ to ${\mathcal Y}$ by $\Hom _{{\mathcal
Z}}\left( {\mathcal X}, {\mathcal Y} \right)$.

\begin{definition}[Category ${\mathfrak C} \left( \Sigma \right) _c$]\label{cat}
Let ${\mathfrak C} \left( \Sigma \right) _{c}$ and ${\mathfrak C} \left( \Sigma \right) _{c} ^l$ be categories defined as follows:\\
{\bf (Objects)} We define ${\bf Ob} {\mathfrak C} \left( \Sigma \right) _{c} ^l$ as pairs $\left(V,
\{F ^{\tau} \left( \bullet \right) \}_{\tau \in \Sigma ( 1 )} \right)$
such that the following three conditions hold:
\begin{enumerate}
\item $V$ is a $\tilde{G} \times \tilde{Z} \left( G \right)$-module;
\item For each $\tau \in \Sigma ( 1 )$, $F ^{\tau} \left( \bullet
\right)$ is an $\tau$-standard filtration of $V$;
\item For each $\sigma \in \Sigma$, a family of subspaces $\{F^{\tau}
\left( n \right) \} _{\tau \in \sigma ( 1 ), n \in \Z}$ of $V$ forms a distributive lattice.
\end{enumerate}
Let ${\mathfrak C} \left( \Sigma \right) _{c}$ be the full-subcategory of ${\mathfrak C} \left( \Sigma \right) _{c} ^l$ obtained by replacing condition 2. by the following stronger condition:
\begin{itemize}
\item[4.] For each $\tau \in \Sigma ( 1 )$, $F ^{\tau} \left( \bullet
\right)$ is an $\tau$-transversal filtration of $V$.
\end{itemize}
{\bf (Morphisms)} For $\left(V _1, \{F ^{\tau} _1 \left( \bullet \right) \}_{\tau \in \Sigma ( 1 )} \right), \left(V _2, \{F ^{\tau} _2 \left( \bullet \right) \}_{\tau \in \Sigma ( 1 )} \right) \in {\bf Ob} {\mathfrak C} \left( \Sigma \right) _{c} ^l$, we define
\begin{align*}
\Hom &  _{{\mathfrak C} \left( \Sigma \right) _{c} ^l} \left( \left(V _1, \{F ^{\tau} _1 \left( \bullet
\right) \}_{\tau \in \Sigma ( 1 )} \right), \left(V _2, \{F ^{\tau} _2 \left( \bullet
\right) \}_{\tau \in \Sigma ( 1 )} \right) \right)\\
:= \{ f & \in \Hom _{\tilde{G}\times \tilde{Z} ( G )} \left( V _1, V _2 \right); f \left( F ^{\tau} _1 \left( n
\right) \right) \subset F ^{\tau} _2 \left( n \right) \text{ for every } n \in \Z \text{ and every } \tau \in \Sigma ( 1 )\}
\end{align*}
For simplicity, we may write $\left(V,
\{F ^{\tau} \left( \bullet \right) \} \right)$ instead of $\left(V,
\{F ^{\tau} \left( \bullet \right) \}_{\tau \in \Sigma ( 1 )} \right)$.
\end{definition}

\begin{remark}
Replacing $\tilde{G} \times \tilde{Z} ( G )$-modules in Definition \ref{cat} by $G \times 1$-modules, we obtain the two fullsubcategories $\mathcal C ( X ) ^+ \subset \mathcal C ( \Sigma ) ^l _c$ and $\mathcal C ( X ) \subset \mathcal C ( \Sigma ) _c$ defined in the introduction.
\end{remark}

\begin{definition}[Category $\mathfrak C ( \Sigma )$]
We define subcategories $\mathfrak C ( \Sigma ) \subset \mathfrak C ( \Sigma ) _c$ and $\mathfrak C ( \Sigma ) ^l \subset \mathfrak C ( \Sigma ) _c ^l$ as follows:\\
{\bf (Objects)} ${\bf Ob} \mathfrak C ( \Sigma ) := {\bf Ob} \mathfrak C ( \Sigma ) _c$ and ${\bf Ob} \mathfrak C ( \Sigma ) ^l := {\bf Ob} \mathfrak C ( \Sigma ) _c ^l$.\\
{\bf (Morphisms)} For $\left(V _1, \{F ^{\tau} _1 \left( \bullet \right) \} \right), \left(V _2, \{F ^{\tau} _2 \left( \bullet \right) \} \right) \in {\bf Ob} {\mathfrak C} \left( \Sigma \right) ^l$, we define
\begin{align*}
\Hom & _{{\mathfrak C} \left( \Sigma \right) ^l} \left( \left(V _1, \{F ^{\tau} _1 \left( \bullet
\right) \}_{\tau \in \Sigma ( 1 )} \right), \left( V _2, \{F ^{\tau} _2 \left( \bullet
\right) \}_{\tau \in \Sigma ( 1 )} \right) \right)\\
:= &  \{ f \in \Hom _{\tilde{G}\times \tilde{Z} ( G )} \left( V _1, V _2 \right); f \text{ satisfies the conditions (L) and (R), where}\\
(L) & \quad f \left( F ^{\tau} _1 \left( n
\right) \right) = f \left( V _1\right) \cap F ^{\tau} _2 \left( n
\right) \text{ for every } n \in \Z \text{ and every } \tau \in \Sigma ( 1 );\\
(R) & \quad \left\{
f \left( V _1 \right), \{F^{\tau} _2
\left( n \right) \} _{\tau \in \sigma ( 1 ), n \in \Z}\right\} \text{
 forms a distributive lattice for every } \sigma \in \Sigma.\}
\end{align*}
We regard ${\mathfrak C} \left( \Sigma \right)$ as a fullsubcategory of ${\mathfrak C} \left( \Sigma \right) ^l$.
\end{definition}

\subsection{The category $EV \left( \Sigma \right)$ and the map $\Xi$}\label{CatEV}

For a $G$-variety $Y$, we denote the category of $G$-equivariant coherent sheaves on $Y$ by $\mathfrak{Coh} ^G Y$.

\begin{definition}[Category of equivariant bundles]
We define two subcategories $EV \left( \Sigma \right)$ and $EV \left( \Sigma \right) _{c}$ of $\mathfrak{Coh} ^{\tilde{G} \times \tilde{G}} X ( \Sigma )$ as follows:
\begin{itemize}
\item ${\bf (Objects)}$ We define ${\bf Ob} EV \left( \Sigma \right)$ and ${\bf Ob} EV \left( \Sigma \right) _c$ as
$$\{ \text{$\tilde{G} \times \tilde{G}$-equivariant vector bundles on $X \left( \Sigma \right)$} \}.$$
\item ${\bf (Morphisms)}$ We define the morphisms of $EV \left( \Sigma \right) _{c}$ as the morphisms of $\tilde{G} \times \tilde{G}$-equivariant coherent ${\mathcal O} _{X \left( \Sigma \right)}$-modules. We define the morphisms of $EV \left( \Sigma \right)$ as morphisms in $EV \left( \Sigma \right) _c$ such that both their kernel and cokernel exist in $EV \left( \Sigma \right)$. We call a morphism of $EV \left( \Sigma \right)$ a ($\tilde{G} \times \tilde{G}$-equivariant) vector bundle morphism.
\end{itemize}
\end{definition}

\begin{remark}
The category $EV \left( \Sigma \right) _c$ introduced above contains the category $EV \left( X \right) = EV \left( X ( \Sigma ) \right)$ in the introduction as a fullsubcategory.
\end{remark}

We denote the total space of a $\tilde{G} \times \tilde{G}$-equivariant vector bundle ${\mathcal E}$ by $V \left( {\mathcal E} \right)$.

\begin{definition}\label{naive Xi}
For each $\tilde{G} \times G$-equivariant vector bundle ${\mathcal E}$, we define a pair $\Xi \left( {\mathcal E} \right)$ as follows.
$$\Xi \left( {\mathcal E} \right) := \left( B ( \mathcal E ), \{ F ^{\tau} \left( n \right) \} _{\tau \in \Sigma ( 1 ), n \in \Z} \right)$$

Here $B ( \mathcal E )$ and $F ^{\tau} \left( n \right)$ are vector spaces such that:
\begin{itemize}
\item $B ( \mathcal E ) := {\mathcal E} \otimes _X k ( e )$ is the identity fiber of ${\mathcal E}$;
\item $F ^{\tau} \left( n \right) := \{ v \in B ( \mathcal E ) ; \exists \lim _{t
\rightarrow \infty} t ^n \left( 1 \times \tau ( t ) \right) v \in V
\left( {\mathcal E} \right) \}$ for every $\tau \in \Sigma ( 1 )$ and
every $n \in \Z$.
\end{itemize}
In particular, $F ^{\tau} ( \bullet )$ is a decreasing filtration of $V$ for each $\tau \in \Sigma ( 1 )$.
\end{definition}
A refinement of $\Xi$ gives an equivalence of $EV ( \Sigma ) _c$ to $\mathfrak C ( \Sigma ) _c$ in \S \ref{MT}.
\section{Foundational results}\label{setting up} 
\subsection{Preliminaries}\label{prelim}
\subsubsection{Isotypical decompositions}\label{isodecomp}
Let $X = X \left( \Sigma \right)$ be the $G \times G$-equivariant partial compactification of $G$ associated with $\Sigma$. Then, the natural projection $\tilde{G} \times \tilde{G} \rightarrow G \times G$ induces a $\tilde{G} \times \tilde{G}$-action on $X \left( \Sigma \right)$. Thus, by the definition of $\tilde{Z} \left( G \right)$, the $1 \times \tilde{Z} \left( G \right)$-action on $X \left( \Sigma \right)$ is trivial. Let ${\mathcal E}$ be a $1 \times \tilde{Z} \left( G \right)$-equivariant vector bundle on $X \left( \Sigma \right)$ (e.g. a $\tilde{G} \times \tilde{G}$-equivariant vector bundle). Then, $1 \times \tilde{Z} \left( G \right)$ operates on ${\mathcal E}
\otimes _X k ( x )$ for every $x \in X \left( \Sigma\right)$. For each $\chi \in \tilde{Z} \left( G \right) ^{\vee}$, we define a sheaf of abelian groups ${\mathcal E} _{\chi}$ on $X \left( \Sigma \right)$ as follows:
$$
{\mathcal E} _{\chi} \left( {\mathcal U} \right) := \Hom _{1 \times \tilde{Z} \left( G \right)} \left( k \boxtimes \chi, {\mathcal E} \left( {\mathcal U} \right) \right) \text{ for every Zariski open subset } {\mathcal U} \subset X \left( \Sigma \right).
$$
Then, $\mathcal E _{\chi}$ is an ${\mathcal O} _X$-submodule of ${\mathcal E}$. The following lemmas easily follow from standard representation theory.

\begin{lemma}\label{isotypical decomposition of morphism}
Let ${\mathcal E}$ and ${\mathcal F}$ be $1 \times \tilde{Z} \left( G \right)$-equivariant vector bundles. Choose $\chi, \xi \in \tilde{Z} \left( G \right) ^{\vee}$. If $\chi \neq \xi$, then we have $\Hom _{( {\mathcal O} _X, 1 \times \tilde{Z} \left( G \right) )} \left( {\mathcal E} _{\chi}, {\mathcal F} _{\xi} \right) = 0$. $\Box$
\end{lemma}

\begin{lemma}\label{isotypical decomposition of vector bundle}
Every $\tilde{G} \times \tilde{G}$-equivariant vector bundle ${\mathcal E}$ admits the following isotypical decomposition:
$${\mathcal E} \cong \oplus _{\chi \in \tilde{Z} \left( G \right) ^{\vee}} {\mathcal E} _{\chi}.$$
Moreover, each ${\mathcal E} _{\chi}$ is a $\tilde{G} \times \tilde{G}$-equivariant vector bundle. $\Box$
\end{lemma}

\begin{definition}[Isotypical components of the category $EV \left( \Sigma \right) _c$]\label{isoeqcat}
For each $\chi \in \tilde{Z} \left( G \right) ^{\vee}$, we define two categories $EV \left( \Sigma, \chi \right) _{c}$ and $EV \left( \Sigma, \chi \right)$ as follows:
\begin{itemize}
\item $({\bf Objects})$ We define ${\bf Ob} EV \left( \Sigma, \chi \right) _{c}$ and ${\bf Ob} EV \left( \Sigma, \chi \right)$ by
$$\{ {\mathcal E} \in {\bf Ob} EV \left( \Sigma \right) _{c} ; {\mathcal E} _{\chi} \cong {\mathcal E} \}.$$
\item $({\bf Morphisms})$ We regard $EV \left( \Sigma, \chi \right) _{c}$ and $EV \left( \Sigma, \chi \right)$ as full-subcategories of $EV \left( \Sigma \right) _{c}$ and $EV \left( \Sigma \right)$, respectively.
\end{itemize}
\end{definition}

Notice that an object of $EV \left( \Sigma, 1 \right)$ can be identified with a $\tilde{G} \times G$-equivariant vector bundle on $X$.

\begin{corollary}[Isotypical decomposition of $EV \left( \Sigma \right) _{c}$]\label{isotypical decomposition of category of vector bundles}
We have a direct sum decomposition
$$EV \left( \Sigma \right) _{c} \cong \oplus_{\chi \in \tilde{Z} \left( G \right) ^{\vee}} EV \left( \Sigma , \chi \right) _{c}$$
as categories.
\end{corollary}

\begin{proof}
This is a direct consequence of Lemma \ref{isotypical decomposition of morphism} and Lemma \ref{isotypical decomposition of vector bundle}.
\end{proof}

\begin{definition}[Isotypical components of the category ${\mathfrak C} \left( \Sigma \right) _{c}$]\label{isotypical components of C}
For each $\chi \in \tilde{Z} \left( G \right) ^{\vee}$, we define three categories ${\mathfrak C} \left( \Sigma , \chi \right)$, ${\mathfrak C} \left( \Sigma , \chi \right) _{c}$, and ${\mathfrak C} \left( \Sigma , \chi \right) _{c} ^l$ as follows:
\begin{itemize}
\item $({\bf Objects})$ We define ${\bf Ob} {\mathfrak C} \left( \Sigma, \chi \right) _{c}$ and ${\bf Ob} {\mathfrak C} \left( \Sigma, \chi \right)$ by
$$\{ \left( V, \{F ^{\tau} \left( \bullet \right) \} \right) \in {\bf Ob} {\mathfrak C} \left( \Sigma \right) _{c} ; \Hom _{1 \times \tilde{Z} \left( G \right)} \left( k \boxtimes {\chi}, V \right) \cong V \}.$$
We also define ${\bf Ob} {\mathfrak C} \left( \Sigma, \chi \right) _{c} ^l$ by
$$\{ \left( V, \{F ^{\tau} \left( \bullet \right) \} \right) \in {\bf Ob} {\mathfrak C} \left( \Sigma \right) _{c} ^l ; \Hom _{1 \times \tilde{Z} \left( G \right)} \left( k \boxtimes {\chi}, V \right) \cong V \}.$$
\item $({\bf Morphisms})$ We consider ${\mathfrak C} \left( \Sigma , \chi \right)$, ${\mathfrak C} \left( \Sigma , \chi \right) _{c}$, and ${\mathfrak C} \left( \Sigma , \chi \right) _{c} ^l$ as full-subcategories of ${\mathfrak C} \left( \Sigma \right)$, ${\mathfrak C} \left( \Sigma \right) _{c}$, and ${\mathfrak C} \left( \Sigma \right) _{c} ^l$, respectively.
\end{itemize}
\end{definition}

We have an analogous decomposition to that of Corollary \ref{isotypical decomposition of category of vector bundles} for ${\mathfrak C} \left( \Sigma \right) _{c}$.

\begin{corollary}[Isotypical decomposition of ${\mathfrak C} \left( \Sigma \right) _{c}$]\label{isotypical decomposition for C}
We have a direct sum decomposition
$${\mathfrak C} \left( \Sigma \right) _{c} \cong \oplus_{\chi \in \tilde{Z} \left( G \right) ^{\vee}} {\mathfrak C} \left( \Sigma , \chi \right) _{c}$$
as categories.
\end{corollary}

\begin{proof}
The assertion follows directly from the complete reducibility of a $\tilde{Z} \left( G \right)$-module.
\end{proof}

\subsubsection{Some lemmas about equivariant structures}\label{lemma on equivariant str}
For equivariant structures on line bundles, we have the following celebrated result of Steinberg.

\begin{theorem}[Steinberg \cite{St} and GIT
\cite{M} 1.4-1.6]\label{Steinberg Theorem}Assume that a connected, simply
connected semisimple linear algebraic group $K$ acts on a normal projective algebraic
variety $Y$. Then, every line bundle on $Y$ admits a unique $K$-linearization. $\Box$
\end{theorem}

Until Lemma \ref{uniq pair}, we assume that $X$ is complete. In vector bundle case, a naive extension of Theorem \ref{Steinberg Theorem} is false. For an arbitrary irreducible $\tilde{G}$-module $V$ of dimension $\ge 2$, $V \times X$ has at least three $\tilde{G} \times \tilde{G}$-equivariant vector bundle structures. One is given by
$$( \tilde{G} \times \tilde{G} ) \times \left( V \times X \right) \ni ( g _1, g _2, v, x ) \mapsto (v, ( g _1, g _2 ). x) \in V \times X,$$
that is $\tilde{G} \times \tilde{G}$-equivariantly isomorphic to a direct sum ${\mathcal O} _X ^{\oplus \dim V}$ of trivial line bundles. The others are given by
\begin{eqnarray}\label{left basic module}
( \tilde{G} \times \tilde{G} ) \times \left( V \times X \right) \ni ( g _1, g _2, v, x ) \mapsto (g _1 . v, ( g _1, g _2 ). x) \in V \times X,\\\label{right basic module}
( \tilde{G} \times \tilde{G} ) \times \left( V \times X \right) \ni ( g _1, g _2, v, x ) \mapsto (g _2 . v, ( g _1, g _2 ). x) \in V \times X.
\end{eqnarray}
We denote the $\tilde{G} \times \tilde{G}$-equivariant vector bundles (\ref{left basic module}) and (\ref{right basic module}) by $V \otimes {\mathcal O} _X$ and ${\mathcal O} _X \otimes V$, respectively. Moreover, we also write $V \otimes {\mathcal L}$ and ${\mathcal L} \otimes V$ for their twist by a $\tilde{G} \times \tilde{G}$-equivariant line bundle ${\mathcal L}$.

Hence, to obtain a proper analogue of Steinberg's theorem for vector bundles, we must impose some auxiliary condition.

\begin{lemma}\label{uniq pair}
Let ${\mathcal E} \hookrightarrow {\mathcal F}$ be an inclusion of vector bundles on $X$ as coherent ${\mathcal O} _X$-modules. Assume that both ${\mathcal E}$ and ${\mathcal F}$ have $\tilde{G} \times \tilde{G}$-equivariant structures such that ${\mathcal E} \MID _G \hookrightarrow {\mathcal F} \MID _G$ as $\tilde{G} \times \tilde{G}$-equivariant ${\mathcal O} _G$-modules. Then, the $\tilde{G} \times \tilde{G}$-equivariant structure of ${\mathcal E}$ is given by the restriction of that of ${\mathcal F}$. $\Box$
\end{lemma}

\begin{corollary}\label{uniqueness}
Let ${\mathcal E}$ be a vector bundle on $X$. If we fix a compatible $( \tilde{G} \times \tilde{G}, k [ G ] )$-module structure on ${\mathcal E} \left( G \right) (= \Gamma ( G, {\mathcal E} ) )$, then $\tilde{G} \times \tilde{G}$-equivariant vector bundle structure on ${\mathcal E}$ is unique if it exists. $\Box$
\end{corollary}

Let $\{ 0 \}$ be the fan consisting of a unique cone $\{ 0 \} \subset X _* \left( T \right) _{\mathbb R}$.

\begin{lemma}\label{Basic Situation}
We have an equivalence of categories $\Ind : {\mathfrak C} \left( \{ 0 \} \right) \stackrel{\cong}{\rightarrow} EV \left( \{ 0 \} \right)$. Its inverse functor is given by the restriction to the fiber at $e$.
\end{lemma}

\begin{proof}
We have $X \left( \{ 0 \} \right) = G$. In particular, $X \left( \{ 0 \} \right)$ is a homogeneous space under $\tilde{G} \times \tilde{G}$-action. Its isotropy group at $e$ is isomorphic to $\tilde{G} \times \tilde{Z} \left( G \right) \ni \left( g, h \right) \mapsto \left(g, gh \right) \in \tilde{G} \times \tilde{G}$. By definition, the set of one-dimensional cones of $\{ 0 \}$ is an empty set. Hence, we have ${\mathfrak C} \left( \{ 0 \} \right) \cong \Rep \tilde{G} \times \tilde{Z} \left( G \right)$. Thus, the desired equivalence is standard (cf. Chriss and Ginzburg \cite{CG} 5.2.16).
\end{proof}

\begin{lemma}\label{switch}
Assume that $X = G$ $($i.e. $\Sigma = \{ 0 \})$. Let $\Ind$ be as in Lemma $\ref{Basic Situation}$. Then, we have
$$\Ind \left( V_{\lambda} \boxtimes k \right) \cong V_{\lambda} \otimes {\mathcal O} _G \cong \Ind \left( k \boxtimes \bar\lambda ^{- 1} \right) \otimes V_{\lambda}$$
for every $\lambda \in X ^* ( \tilde{T} )$. Here $\bar\lambda$ is the image of $\lambda$ in $\tilde{Z} ( G ) ^{\vee}$ as in \S \ref{NoteVar}.
\end{lemma}

\begin{proof}
As a $\tilde{G} \times \tilde{Z} \left( G \right)$-module, we have $V_{\lambda} \otimes {\mathcal O}_X \otimes _X k \left( e \right) \cong V_{\lambda} \boxtimes k$. Hence, we have the first isomorphism. Next, we prove $\Ind \left( V_{\lambda} \boxtimes k \right) \cong \Ind \left( k \boxtimes \bar\lambda ^{- 1} \right) \otimes V_{\lambda}$ to complete the proof. $k \boxtimes \bar\lambda ^{- 1}$ is trivial as a $\tilde{G} \times 1 ( \subset \tilde{G} \times \tilde{Z} \left( G \right) )$-module. Thus, we have a $\tilde{G} \times 1$-module isomorphism $\left( \Ind \left( k \boxtimes \bar\lambda ^{- 1} \right) \otimes V _{\lambda} \right) \otimes _X k \left( e \right) \cong V _{\lambda}$. Here the action of $1 \times \tilde{Z} \left( G \right) \subset \tilde{G} \times \tilde{Z} \left( G \right)$ is the right action. This action is trivial since $\tilde{Z} \left( G \right)$ acts by $- \bar\lambda + \bar\lambda = 0$. As a result, we have $\Ind \left( V_{\lambda} \boxtimes k \right) \cong \Ind \left( k \boxtimes \bar\lambda ^{- 1} \right) \otimes V_{\lambda}$.
\end{proof}

\subsubsection{Equivariant divisors}
In this section, we introduce and explain the notion of equivariant divisors. First, we present their definition. Recall that $\overline{\mathbf O} _{\sigma}$ is the closure of the $G \times G$-orbit of $X$ corresponding to $\sigma \in \Sigma$.

\begin{definition}[Equivariant Divisors]
An equivariant divisor is a $\Z$-linear formal sum of $D _{\tau}$ for every $\tau \in \Sigma ( 1 )$. For two equivariant divisors $D _1 = \sum n _1 ^{\tau} D _{\tau}$ and $D _2 = \sum n _2 ^{\tau} D _{\tau}$, we say $D _1 \ge D _2$ if, and only if, we have $n _1 ^{\tau} \ge n _2 ^{\tau}$ for every $\tau \in \Sigma ( 1 )$. Moreover, $D _1$ is said to be sufficiently large if, and only if, $n _1 ^{\tau} >> 0$ for every $\tau \in \Sigma ( 1 )$.
\end{definition}

Before we exploit some properties of equivariant divisors, we need a result.

\begin{lemma}\label{ambient}
There exists a fan $\Sigma _+$ such that $X (\Sigma _+)$ is a complete regular embedding of $G$ which (equivariantly) contains $X ( \Sigma )$.
\end{lemma}

\begin{proof}
By the theory of spherical embeddings, we have a $G \times G$-equivariant complete embedding $Y$ of $G$ which contains $X ( \Sigma )$ (cf. \cite{Kn}). Then, we successively blow-up along the singular locus (which is $G \times G$-stable) to obtain a regular embedding. Now Theorem \ref{fundamental} gives the result.
\end{proof}

We fix one $\Sigma _+$ of Lemma \ref{ambient} hereafter.

\begin{theorem}[See \cite{Br} 2.2, \cite{Bi} 2.4, and \cite{Kl} 2.1-2.2]\label{embed ed}
We have the following:
\begin{enumerate}
\item For each $\sigma \in \Sigma ( r )$, the restriction map
$$\iota _{\sigma} ^* : \bigoplus _{\tau \in \sigma ( 1 )} \Z D _{\tau} \hookrightarrow \Pic ^{1 \times \tilde{G}} \overline{\mathbf O} _{\sigma}$$
is an injection.
\item We have the following short exact sequence.
$$0 \rightarrow \bigoplus _{\tau \in \Sigma ( 1 )} \Z D _{\tau} \rightarrow \Pic ^{1 \times \tilde{G}} X \left( \Sigma \right) \stackrel{\kappa}{\rightarrow} \tilde{Z} \left( G \right) ^{\vee} \rightarrow 0$$
\item For each equivariant divisor $D = \sum _{\tau \in \Sigma ( 1 )} n _{\tau} D _{\tau}$, $\mathbb{G} _m ^{\tau}$ acts on ${\mathcal O} _X \left( D \right) \otimes _X k ( x _{\tau} )$ by weight $- n _{\tau}$. $(\mathbb{G} _m ^{\tau}$ is the image of $\mathbb{G} _m$ under $1 \times \tau ^{- 1}$. See \S \ref{NoteVar}.$)$
\end{enumerate}
Here $\Pic ^{1 \times \tilde{G}} \bullet$ is the $1 \times \tilde{G}$-equivariant Picard group. Moreover, the image of $\kappa$ determines the right $1 \times \tilde{Z} ( G ) (\subset \tilde{G} \times \tilde{Z} ( G ) )$-module structure of the identity fiber.
\end{theorem}

\begin{remark}
The origin of Theorem \ref{embed ed} is somewhat complicated. If $G$ is an adjoint semisimple group and $X$ is complete, the above formulation is essentially due to Bifet \cite{Bi}. For another extreme of our scope, namely toric varieties, the description of the equivariant Picard group is a corollary of Klyachko's theorem \cite{Kl}. In the meantime, Brion \cite{Br} established a general description of Picard groups using $B \times B ^-$-orbits. Hence, Theorem \ref{embed ed} is essentially a corollary of his result (but not a direct consequence). Anyway, since the author could not find a proper reference to this form of the theorem, we provide a proof.
\end{remark}

For the proof of Theorem \ref{embed ed}, we need some preparation.

The following is a modification of Strickland's Theorem \cite{Str} 2.4 to our setting. It is easily deduced from \cite{Str} 2.4 by twisting a character of $\tilde{T} _0 \times \tilde{T} _0$.

\begin{theorem}[Strickland's theorem cf. \cite{Str} 2.4]\label{Strickland}
Let $p _0$ be the unique $B \times B ^-$-fixed point in $X _0$. For each $\lambda \in X ^* ( \tilde{T} )$, there exists a $\tilde{G} \times \tilde{G}$-equivariant line bundle $\bar{\mathcal L} _{\lambda}$ on $X _0$ such that $\bar{\mathcal L} _{\lambda} \otimes _{{\mathcal O} _{X _0}} k ( p _0 ) \cong \lambda ^{-1} \boxtimes \lambda$ as $\tilde{T} \times \tilde{T}$-modules. Moreover, $\bar{\mathcal L} _{\lambda}$ is unique up to isomorphism.
\end{theorem}

Recall the morphism $\pi : X \left( \Sigma \right) \rightarrow X _0$ introduced in Theorem \ref{fundamental}.

\begin{corollary}\label{boundary injection}
For each $\sigma \in \Sigma ( r )$, the map
$$\jmath _{\sigma} : \Pic ^{1 \times
\tilde{G}} X \left( \Sigma \right) \ni {\mathcal L} \mapsto
{\mathcal L} \otimes _X k ( x _{\sigma} ) \in  X ^* ( \tilde{T} ) \cong \Rep ( 1 \times \tilde{T} )$$
is a surjection.
\end{corollary}

\begin{proof}
$\mathbf O _{\sigma}$ is a closed $G \times G$-orbit of $X \left( \Sigma \right)$ with dimension $\dim G - r$. Thus, $\pi \left( \mathbf O _{\sigma} \right)$ is the closed $G \times G$-orbit of $X _0$. As a result, $\pi ( x _{\sigma} )$ is the unique $B \times B ^-$-fixed point of $X _0$. Hence, for each $\lambda \in X ^{*} ( \tilde{T} )$, we have $\jmath _{\sigma} \left(  [ \pi ^* \bar{\mathcal L} _{\lambda} ] \right) = \lambda$.
\end{proof}

Next, we review the Key Proposition of \cite{Kl}.

\begin{theorem}[Klyachko \cite{Kl} 2.1.1]\label{Klyachko's Proposition}
For each $\sigma \in \Sigma$, we have the following:
\begin{enumerate}
\item Every $1 \times T$-equivariant vector bundle ${\mathcal E}$ on $T \left( \sigma \right)$ is uniquely written as $E \times T \left( \sigma \right)$ by a $T$-module $E$. $($Here $T$ acts on both $E$ and $T \left( \sigma \right)$.$)$
\item Let $E, F$ be $T$-modules. Then, two $1 \times T$-equivariant vector bundles $E \times T \left( \sigma \right)$ and $F \times T \left( \sigma \right)$ are isomorphic if, and only if, $E \cong F$ as $T$-modules.
\item Let $E$ be a $T$-module. Let $E = \oplus _{\lambda \in X ^* \left( T \right)} E ^{\lambda}$ be its $T$-isotypical decomposition. Choose $v \in E ^{\lambda} \subset E \otimes k ( e )$ and $\tau \in \Sigma ( 1 )$. Then,
$$\lim _{t \rightarrow \infty} t ^{n} ( 1 \times \tau ( t ) ) v$$
exists in $E \times T \left( \sigma \right)$ if, and only if, $\left< \tau, \lambda \right> \le - n$.
\end{enumerate}
\end{theorem}

\begin{proof}[Proof of Theorem \ref{embed ed}]
First, we prove 1). We restrict our attention to $T \left( \Sigma \right)$. Let $\{ \tau _1, \tau _2, \ldots \tau _r \}$ be the $\Z$-basis of $X _* \left( T \right)$ which spans $\sigma$. Let $\{ \mu _1, \mu _2, \ldots \mu _r \}$ be its dual basis (i.e. $\left< \tau _i, \mu _j \right> = \delta _{i, j}$). By Theorem \ref{structure theorem}, $T \left( \sigma \right) \cap D _{\tau}$ is a $1 \times T$-stable divisor if $\tau \in \sigma ( 1 )$. Moreover, by the definition of the coordinate ring of toric varieties \cite{Oda} 1.2, we have ${\mathcal O} _X \left( D _{\tau _i} \right) \otimes _X k ( x _{\sigma} ) = \mu _i$ as $\tilde{T}$-module for each $1 \le i \le r$. As a result, the image of $D _{\tau _i}$ under the composition map $\oplus _{i = 1} ^r \Z D _{\tau _i} \rightarrow \Pic ^{1 \times \tilde{G}} X \left( \Sigma \right) \rightarrow X ^* ( \tilde{T} )$ is $\mu _i ( \in X ^* ( T ) )$ for each $i$. In particular, 1) follows.
Combined with Theorem \ref{Klyachko's Proposition} 1), we have 3) if $\Sigma = \Sigma _+$. Thus, 3) follows from Lemma \ref{ambient} by restriction.

Now we prove 2). If we have $\sum _{\tau \in \Sigma ( 1 )} a _{\tau} D _{\tau} \sim 0$ in $\mathrm{Pic} ^{1 \times \tilde{G}} X \left( \Sigma _+ \right)$, then there exist $f \in k [ G ] ^{\times}$ such that $\sum _{\tau \in \Sigma _+ ( 1 )} a _{\tau} D _{\tau} = \mathrm{div} ( f )$. By a theorem of Rosenlicht, we have $k [ G ] ^{\times} = k ^{\times} X ^* \left( G \right)$. Hence, 1) yields $a _{\tau} = 0$ for every $\tau \in \Sigma ( 1 )$. Therefore, we have a short exact sequence
$$0 \rightarrow \bigoplus _{\tau \in \Sigma _+ ( 1 )} \Z D _{\tau} \rightarrow \Pic ^{1 \times \tilde{G}} X \left( \Sigma _+ \right) \rightarrow \Pic ^{1 \times \tilde{G}} G \rightarrow 0.$$
We have $\Pic ^{1 \times \tilde{G}} G \cong \tilde{Z} ( G ) ^{\vee}$. By restricting the above short exact sequence to $X ( \Sigma )$, we obtain 2) and the assertion about $\kappa$.
\end{proof}

By Theorem \ref{embed ed}, the pullback of the line bundle $\bar{\mathcal L} _{\lambda}$ defined in Theorem \ref{Strickland} corresponds to an element $D ^{\lambda}$ of $\Pic ^{1 \times \tilde{G}} X \left( \Sigma \right)$. We denote the $\tilde{G} \times \tilde{G}$-equivariant line bundle $\pi ^* \bar{\mathcal L} _{\lambda}$ by ${\mathcal L} _{\lambda}$ (when we want to stress characters) or ${\mathcal O} _X \left( D ^{\lambda} \right)$ (when we want to stress divisors).

\begin{lemma}\label{divisor}
Let $\lambda \in X ^* ( \tilde{T} )$. We have $D ^{\lambda} = \sum _{\tau \in \Sigma ( 1 )} \left< \tau, \lambda \right> D _{\tau}$ in $\Pic ^{1 \times \tilde{G}} X \left( \Sigma \right) \otimes _{\Z} \Q$. In particular, the RHS is an element of $\Pic ^{1 \times \tilde{G}} X \left( \Sigma \right)$.
\end{lemma}

\begin{proof}
First, we consider the case $\Sigma = \Sigma _+$. We use the same notation
 as in the proof of Theorem \ref{embed ed} 1). Let $\sigma \in \Sigma ( r )$. From the last step of the proof of Theorem \ref{embed ed} 1), we have $\Pic
^{1 \times \tilde{G}} X \left( \sigma \right) \supset \oplus _{i = 1} ^r \Z \mu _i \cong X ^* ( T ) \subset X ^* ( \tilde{T}
)$. By Corollary \ref{boundary injection}, $\Pic
^{1 \times \tilde{G}} X \left( \sigma \right) \rightarrow X ^* ( \tilde{T}
)$ is a surjection. By Theorem \ref{embed ed} 2) (applied to the case $\Sigma = \Sigma _+ = \sigma$), we conclude $\Pic
^{1 \times \tilde{G}} X \left( \sigma \right) \cong X ^* ( \tilde{T}
)$ since we have $| \tilde{Z} ( G ) ^{\vee} | = | X ^* ( \tilde{T} ) / X ^* ( T ) |$. Thus, we have $D ^{\lambda} \MID _{X ( \sigma )} = \sum _{\tau \in \sigma
( 1 )} \left< \tau, \lambda \right> D _{\tau} \in \Pic ^{1
\times \tilde{G}} X \left( \sigma \right)$. We have the following exact sequence of free $\Z$-modules:
$$0 \rightarrow \bigoplus _{\tau \in \Sigma _+ ( 1 ) \backslash \sigma ( 1 )} \Z D _{\tau} \rightarrow \Pic ^{1 \times \tilde{G}} X \left( \Sigma _+ \right) \rightarrow \Pic ^{1 \times \tilde{G}} X \left( \sigma \right)$$
Hence, $\oplus _{\sigma \in
\Sigma _+ ( r )}
\jmath _{\sigma} : \Pic ^{1 \times \tilde{G}} X \left( \Sigma _ + \right)
\rightarrow \oplus _{\sigma \in \Sigma + ( r )} X ^* ( \tilde{T} )$
is injective by Theorem \ref{embed ed} 1)
and the above arguments. As a result, we conclude $D ^{\lambda} = \sum _{\tau \in
\Sigma _+
( 1 )} \left< \tau, \lambda \right> D _{\tau} \in \Pic ^{1
\times \tilde{G}} \left( \Sigma _+ \right)$. Therefore, we have the result if $\Sigma = \Sigma _+$. Since $\pi : X \left( \Sigma \right) \rightarrow X _0$ factors through $X \left( \Sigma _+\right)$, the general case follows.
\end{proof}

For each $\tilde{G} \times \tilde{G}$-equivariant line bundle ${\mathcal L}$, $\otimes _X {\mathcal
L}$ yields the following category (auto-) equivalence.
\begin{align*}
\otimes _X {\mathcal L} & : {\bf Ob} EV \left( \Sigma \right) _{c} \ni
{\mathcal E} \mapsto {\mathcal E} \otimes _X {\mathcal L} \in {\bf Ob} EV
\left( \Sigma \right) _{c}\\
\otimes _X {\mathcal L} & : \Hom _{EV \left( \Sigma \right) _{c}} \left(
{\mathcal E}, {\mathcal F} \right) \stackrel{\cong}{\rightarrow} \Hom _{EV \left( \Sigma \right) _{c}} \left(
{\mathcal E} \otimes _X {\mathcal L}, {\mathcal F} \otimes _X {\mathcal L}
\right) \text{ for every } {\mathcal E}, {\mathcal F} \in {\bf Ob} EV \left( \Sigma \right) _{c}.
\end{align*}
Thus, we have the following.

\begin{corollary}\label{reduction to trivial character}
We have a category equivalence $EV \left( \Sigma, \chi \right) _{c} \cong EV \left( \Sigma, \xi \right) _{c}$ for every $\chi, \xi \in \tilde{Z} \left( G \right) ^{\vee}$.
\end{corollary}

\begin{proof}
By Theorem \ref{embed ed} 2), we have a $\tilde{G} \times \tilde{G}$-equivariant line bundle ${\mathcal L}$ such that $\kappa \left( {\mathcal L} \right) = \xi - \chi$. Hence, the result follows from the fact that $\otimes _X {\mathcal L}$ yields an auto-equivalence of $EV \left( \Sigma \right) _{c}$ and the definition of $EV \left( \Sigma, \chi \right) _{c}$.
\end{proof}

\subsection{Redefinition of the functor $\Xi$}\label{reconst}
In this subsection, we rewrite the definition of the functor $\Xi$ in order to handle subtle behaviors at the boundaries. Main ingredients in this subsection are the introduction of another functor $\Xi ^{\prime}$ (Definition \ref{main functor for trivial factor} and Definition \ref{general redef}), the proof of $\Xi = \Xi ^{\prime}$ (Proposition \ref{identification}), and a description of $\Xi$ (Proposition \ref{naive functor}).

\subsubsection{Bounding an equivariant bundle by standard ones}

\begin{lemma}\label{incl}
Let ${\mathcal E}$ be a $\tilde{G} \times \tilde{G}$-equivariant vector bundle. Assume that we have $V \boxtimes k \subset \Gamma \left( X, {\mathcal E} \right)$ $($as $\tilde{G} \times \tilde{G}$-modules$)$ for some $\tilde{G}$-module $V$. Then we have a unique $\tilde{G} \times \tilde{G}$-equivariant inclusion $V \otimes {\mathcal O} _X \hookrightarrow {\mathcal E}$ such that its image is generated by $V \boxtimes k$.
\end{lemma}

\begin{proof}
There exists a $\tilde{G}$-module $Z$ such that we have $Z \boxtimes k \cong {\mathcal
E} _1 \otimes _X k ( e )$ as $\tilde{G} \times \tilde{Z} ( G
)$-modules (cf. Lemma \ref{Basic Situation}). Then, we have $Z \otimes {\mathcal O} _G \cong {\mathcal
E} _1 \MID _{G}$ by Lemma \ref{switch}. By looking at $1 \times \tilde{Z} ( G )$-action, we have $\Gamma \left(G, {\mathcal E}
\right) ^{1 \times \tilde{G}} = \Gamma \left(G, {\mathcal E} _1
\right) ^{1 \times \tilde{G}}$. (Here superscript $1 \times \tilde{G}$ means the fixed part as in \S \ref{NoteAlg}.) Hence, we have $\Gamma \left(G,
{\mathcal E} \right) ^{1 \times \tilde{G}} \cong Z \boxtimes k$ as $\tilde{G} \times \tilde{G}$-modules. As a
consequence, we have an induced inclusion $V \boxtimes k \subset Z
\boxtimes k$ as $\tilde{G} \times \tilde{G}$-submodules of
$\Gamma \left( X, {\mathcal E} \right)$. Thus, the induced ${\mathcal O} _X$-module morphism
$\phi : V \otimes {\mathcal O}_X \rightarrow {\mathcal E}$ is a $\tilde{G} \times \tilde{G}$-equivariant inclusion when restricted to $G$ from Lemma \ref{switch} and Lemma \ref{Basic Situation}. Hence, Lemma \ref{uniq pair} yields the result.
\end{proof}

For a $\tilde{G} \times \tilde{G}$-equivariant vector bundle $\mathcal E$ and an equivariant divisor $D$, we often write $\mathcal E ( D )$ instead of $\mathcal E \otimes _X \mathcal O _X ( D )$.

\begin{corollary}\label{uniq-incl}
Let ${\mathcal E}$ be a $\tilde{G} \times \tilde{G}$-equivariant vector bundle. Let $W$ be a $\tilde{G}$-module. Assume that we have $V \boxtimes k = \Gamma \left( G, {\mathcal E} \right) ^{1 \times \tilde{G}}$ for some $\tilde{G}$-module $V$. Then, for a sufficiently large equivariant divisor $D$, every $\tilde{G} \times \tilde{G}$-equivariant morphism $W \otimes {\mathcal O} _X ( - D ) \rightarrow {\mathcal E}$ factors through $V \otimes {\mathcal O} _X ( - D )$. In particular, we have a natural isomorphism
$$\Hom _{EV \left( \Sigma, 1 \right) _{c}} \left( W \otimes {\mathcal O} _X \left( - D \right), {\mathcal E} \right) \cong \Hom _{\tilde{G}} \left( W, V \right).$$
\end{corollary}

\begin{proof}
Let us choose a (sufficiently large) equivariant divisor $D _0$ such that $V \boxtimes k \subset \Gamma \left( X, {\mathcal E} \left( D _0 \right) \right)$ as a $\tilde{G} \times \tilde{G}$-submodule. We have $V \boxtimes k
\subset \Gamma \left( X, {\mathcal E} \left( D \right) \right)$ for all $D \ge D_0$. Thus, we have $V \boxtimes k = \Gamma \left( X, {\mathcal E} \left( D \right) \right) ^{1 \times \tilde{G}}$. By Lemma \ref{incl}, we have a $\tilde{G} \times \tilde{G}$-equivariant inclusion $V \otimes
{\mathcal O} _X \subset {\mathcal E} \left( D \right)$. A $\tilde{G} \times \tilde{G}$-equivariant morphism $W \otimes {\mathcal O}_X \rightarrow {\mathcal E} \left( D \right)$ defines a $\tilde{G} \times \tilde{G}$-module morphism between their sections $\Gamma \left(X, W \otimes {\mathcal O}_X \right) = W \boxtimes k \rightarrow V \boxtimes k = \Gamma \left(X, {\mathcal E} \left( D \right) \right) ^{1 \times \tilde{G}}$. Since $W \otimes {\mathcal O}_X$ is generated by its global sections, every $\tilde{G} \times \tilde{G}$-equivariant morphism $W \otimes {\mathcal O}_X \rightarrow {\mathcal E}\left( D \right)$ factors through $V \otimes {\mathcal O}_X$. Therefore, twisting by ${\mathcal O} _X \left( - D \right)$ completes the proof.
\end{proof}

From now on, we sometimes deal with $\tilde{G} \times G$-equivariant vector bundles ($=$ objects in $EV ( \Sigma, 1 )$). This strange restriction comes from two rather technical reasons. On one hand, the representation categories of $\tilde{G}$ and $G$ are essentially different. Hence, we cannot forget the $\tilde{G}$-action completely. On the other hand, we use equivariant divisors which correspond only to $G \times G$-equivariant line bundles.

For each $\tilde{G} \times G$-equivariant vector bundle $\mathcal E$ and an equivariant divisor $D$, we denote the sheaf $B ( \mathcal E ) \otimes _X \mathcal O _X ( D )$ by $\mathcal B ( \mathcal E ) ^D$.

\begin{corollary}\label{sandwich}
For every ${\mathcal E} \in {\bf Ob} EV \left( \Sigma, 1 \right) _{c}$, there exists a sufficiently large equivariant divisor $D _0$ such that:
\begin{itemize}
\item $\mathcal B ( {\mathcal E} ) ^{- D} \subset {\mathcal E} \subset \mathcal B ( {\mathcal E} ) ^{D}$ are inclusions of $\tilde{G} \times \tilde{G}$-equivariant coherent sheaves for every $D \ge D _0$;
\item The composition map $\mathcal B ( {\mathcal E} ) ^{- D} \subset \mathcal B ( {\mathcal E} ) ^{D}$ is the tensor product of the identity map of $B ( {\mathcal E} )$ and a unique $\tilde{G} \times \tilde{G}$-equivariant inclusion ${\mathcal O} _X ( - D ) \subset {\mathcal O} _X ( D )$.
\end{itemize}
Moreover, such a series of inclusions is unique up to automorphism of $B ( {\mathcal E} )$ as a $\tilde{G}$-module.
\end{corollary}

\begin{proof}
We have $B ( {\mathcal E} )
\boxtimes k \subset B ( {\mathcal E} ) \otimes k [ G ] \cong \Gamma \left( G, {\mathcal E} \right)$ as $\tilde{G} \times \tilde{G}$-modules by the algebraic Peter-Weyl theorem and Lemma \ref{switch}. Hence, there exists a sufficiently large equivariant
divisor $D _0$ such that we have $( B ( {\mathcal E} )
\boxtimes k ) \subset \Gamma \left( X, {\mathcal E} ( D ) \right)$ for every $D \ge D _0$. Thus, we have a $\tilde{G} \times \tilde{G}$-equivariant embedding $\mathcal B ( {\mathcal E} ) ^{- D}
\subset {\mathcal E}$ by Lemma
\ref{incl}. Moreover, such an embedding is uniquely determined up to an automorphism of $B ( {\mathcal E} )$ as a $\tilde{G}$-module. Next, we
apply the same argument to ${\mathcal E} ^{\vee}$. Then, enlarge $D
_0$ if necessary, we obtain $B ( {\mathcal E} ) ^* \otimes {\mathcal O} _X ( - D )
\subset {\mathcal E} ^{\vee}$ for every $D \ge D _0$. Hence, we have a $\tilde{G} \times \tilde{G}$-equivariant inclusion ${\mathcal E}
\subset \mathcal B ( {\mathcal E} ) ^{D}$ since the both vector bundles have the same rank. By Corollary \ref{uniq-incl}, we can use an automorphism of $B ( {\mathcal E} )$ (as $\tilde{G}$-modules) to obtain the desired series of inclusions.
\end{proof}

\begin{lemma}[Main observation of this subsection]\label{intersection of equivariant vector bundles}
Consider two objects ${\mathcal E}$ and ${\mathcal F}$ of $EV \left( \Sigma, 1 \right)
_{c}$ such that $B ( {\mathcal E} ) \cong B ( {\mathcal F} )$ as $\tilde{G}$-modules. Then, we can take the intersection ${\mathcal E} \cap _{D _0} {\mathcal F}$
of ${\mathcal E}$ and ${\mathcal F}$ as $\tilde{G} \times \tilde{G}$-equivariant coherent subsheaves of $\mathcal B ( {\mathcal E} ) ^{D _0}$ for a sufficiently large equivariant divisor $D _0$. Moreover, we have the following commutative diagram of $\tilde{G} \times \tilde{G}$-equivariant coherent sheaves:
\begin{eqnarray}
\begin{matrix}
\mathcal B ( {\mathcal E} ) ^{D _0}& \hookrightarrow & \mathcal B ( {\mathcal E} ) ^{D} &\\
\bigcup & & \bigcup & \text{ for every }D \ge D _0.\\
{\mathcal E} \cap _{D _0} {\mathcal F} & \cong & {\mathcal E} \cap _{D} {\mathcal F} &
\end{matrix}
\end{eqnarray}
\end{lemma}

\begin{proof}
We have an embedding ${\mathcal E} \hookrightarrow \mathcal B ( {\mathcal E} ) ^{D _0}$ as $\tilde{G} \times \tilde{G}$-equivariant coherent sheaves by Corollary
\ref{sandwich}. By the isomorphism $B ( {\mathcal E} ) \cong B ( {\mathcal F} )$, we also have an embedding ${\mathcal F}
\hookrightarrow \mathcal B ( {\mathcal E} ) ^{D _0}$ as $\tilde{G} \times \tilde{G}$-equivariant coherent sheaves.
Hence, we can take their intersection. Here, ${\mathcal E} \cap _{D _0}
{\mathcal F}$ is a $\tilde{G} \times \tilde{G}$-equivariant coherent subsheaf of $\mathcal B ( {\mathcal E} ) ^{D _0}$ because the $\tilde{G} \times \tilde{G}$-equivariant structure of $\mathcal B ( {\mathcal E} ) ^{D _0}$
preserves both ${\mathcal E}$ and ${\mathcal F}$ by Lemma \ref{uniq pair}. By Corollary \ref{uniq-incl}, every $\tilde{G} \times \tilde{G}$-equivariant embedding $B ( {\mathcal E} ) ^* \otimes {\mathcal O} _X ( - D ) \hookrightarrow {\mathcal E} ^{\vee}$(resp. ${\mathcal F} ^{\vee}$) factors through $B ( {\mathcal E} ) ^* \otimes {\mathcal O} _X ( - D _0 ) \hookrightarrow {\mathcal E} ^{\vee}$(resp. ${\mathcal F} ^{\vee}$). By taking their dual, we obtain the second assertion.
\end{proof}

\begin{remark}
From now on, we freely use the notion of the intersection of two
 $\tilde{G} \times \tilde{G}$-equivariant vector bundles in the sense of Lemma \ref{intersection of
 equivariant vector bundles}. Notice that such an intersection is
 defined if, and only if, we specify an isomorphism between
 identity fibers (as $\tilde{G} \times \tilde{Z} ( G )$-modules) and a sufficiently large equivariant divisor.
\end{remark}

\subsubsection{Redefinition of $\Xi$}
Now we redefine the map $\Xi$ introduced in \S \ref{CatEV}. We first introduce another map $\Xi ^{\prime}$ and prove $\Xi = \Xi ^{\prime}$ at Proposition \ref{identification}. At the same time, we extend the definition of the map $\Xi$ to the whole of ${\bf Ob} EV \left( \Sigma \right) _{c}$.

\begin{definition}[$\Xi ^{\prime}$ for $EV \left( \Sigma, 1 \right) _{c}$]\label{main functor for trivial factor}
For each ${\mathcal E} \in {\bf Ob} EV
\left( \Sigma, 1 \right) _{c}$, we define a pair
$$\Xi ^{\prime} \left( {\mathcal E} \right) := \left( B \left(
{\mathcal E} \right), \{ {}^{\prime} F ^{\tau} \left( \bullet,
{\mathcal E} \right) \} _{\tau \in \Sigma ( 1 )} \right)$$
by the following rules:
\begin{itemize}
\item $B \left( {\mathcal E} \right) := {\mathcal E} \otimes _X k ( e )$ as a $\tilde{G} \times \tilde{Z} ( G )$-module. $($By assumption, the $1 \times \tilde{Z} ( G )$-action is trivial.$)$
\item For each $\tau \in \Sigma ( 1 )$, let $\phi _{\tau} : B \left( {\mathcal E} \right) \otimes
{\mathcal O} _X \rightarrow B \left( {\mathcal E} \right) \otimes k (
x _{\tau} )$ be the residual map at $x _{\tau}$. Then we define
$${} ^{\prime} F ^{\tau} \left( n, {\mathcal E} \right) := \phi _{\tau} \left( {\mathcal E} ( - n D _{\tau} ) \cap _D B \left( {\mathcal E} \right) \otimes {\mathcal O} _X \right)$$
for every $n \in \Z$ and a sufficiently large equivariant divisor $D$.
\end{itemize}
We call this map $\Xi ^{\prime}$. Here the existence (and the stability) of intersections with respect to $D$ is guaranteed by Lemma
\ref{intersection of equivariant vector bundles}. By fixing an isomorphism $k (e) \cong k ( x _{\tau} )$, we have an inclusion ${} ^{\prime} F ^{\tau} \left( n, {\mathcal E} \right) \subset B \left( {\mathcal E} \right)$ for every $\tau \in \Sigma ( 1 )$ and every $n \in \Z$.
\end{definition}

\begin{remark}
The isomorphism $k (e) \cong k ( x _{\tau} )$ in Definition \ref{main
functor for trivial factor} looks like quite artificial. However,
since we deal with only vector spaces, we cannot detect the diagonal
$\mathbb{G} _{m}$-action for each vector space. If $X$ is complete, we
can use the space of the global sections of $B \left( {\mathcal E} \right) \otimes {\mathcal O} _X$ to canonicalize  $k (e) \cong k ( x _{\tau} )$.
\end{remark}

For a rational number $x$, $\lfloor x \rfloor$ denotes the maximal integer which does not exceed $x$.

\begin{definition}[$\Xi ^{\prime}$ for $EV \left( \Sigma \right) _{c}$]\label{general redef}
Let $\Lambda = \{ \lambda _1, \lambda _2, \ldots, \lambda _{h} \}$ be a set of complete representative of $X ^* ( \tilde{T} )
\rightarrow \tilde{Z} ( G ) ^{\vee}$. For each $1 \le i \le h$, we have a $\tilde{G} \times \tilde{G}$-equivariant line bundle
${\mathcal L} _{\lambda _i}$ defined before Lemma
\ref{divisor}. By Lemma \ref{isotypical decomposition of vector bundle}, every $\tilde{G} \times \tilde{G}$-equivariant vector bundle ${\mathcal E}$ is uniquely written as ${\mathcal E} \cong \oplus _{i = 1} ^{h} {\mathcal E} _i \otimes _X {\mathcal L} _{\lambda _i}$ for some ${\mathcal E} _1, {\mathcal E} _2, \ldots, {\mathcal E} _{h} \in {\bf Ob} EV \left( \Sigma, 1 \right) _c$. Using the above decomposition, we define a pair $\Xi ^{\prime} _{\Lambda} \left( {\mathcal E} \right)$ as the direct sum
$$\Xi ^{\prime} _{\Lambda} \left( {\mathcal E} \right) = \left( \bigoplus _{i = 1} ^{h} B \left( {\mathcal E} _i \otimes _X {\mathcal L} _{\lambda _i} \right), \left\{ \bigoplus _{i = 1} ^{h} {}^{\prime} F ^{\tau} \left( \bullet, {\mathcal E} _i \otimes _X {\mathcal L} _{\lambda _i} \right) \right\} _{\tau \in \Sigma ( 1 )} \right)$$
of pairs $\left( B \left( {\mathcal E} _i \otimes _X {\mathcal L} _{\lambda _i} \right), \{ {}^{\prime} F ^{\tau} \left( \bullet, {\mathcal E} _i \otimes _X {\mathcal L} _{\lambda _i} \right) \} _{\tau \in \Sigma ( 1 )} \right)$ determined by the following rules:
\begin{itemize}
\item $B \left( {\mathcal E} _i \otimes _X {\mathcal L} _{\lambda _i} \right) := B \left( {\mathcal E} _i \right) \otimes ( k \boxtimes \bar{\lambda} _i )$ as a $\tilde{G} \times \tilde{Z} ( G )$-module;
\item For each $\tau \in \Sigma ( 1 )$, let $\phi _{\tau} ^i : B \left( {\mathcal E} _i \right) \otimes
{\mathcal O} _X \rightarrow B \left( {\mathcal E} _i \right) \otimes k (
x _{\tau} )$ be the residual map at $x _{\tau}$. Then, we define
$${} ^{\prime} F ^{\tau} \left( n + \lfloor \tau ( \lambda _i ) \rfloor, {\mathcal E} _i \otimes {\mathcal L}_{\lambda _i} \right) := \phi _{\tau} ^i \left( {\mathcal E} _i ( - n D _{\tau} ) \cap _D B \left( {\mathcal E} _i \right) \otimes {\mathcal O} _X \right) \otimes ( k \boxtimes \bar{\lambda} _i )$$ for every $n \in \Z$ and a sufficiently large equivariant divisor $D$.
\end{itemize}
We have an inclusion
$$
{} ^{\prime} F ^{\tau} \left( n, {\mathcal E} \right) = \bigoplus _{i = 1} ^{h} {} ^{\prime} F ^{\tau} \left( n, {\mathcal E} _{i} \otimes _X {\mathcal L} _{\lambda _i} \right) \subset \bigoplus _{i = 1} ^{h} B  \left( {\mathcal E} _{i} \otimes _X {\mathcal L} _{\lambda _i} \right) = B \left( {\mathcal E} \right)
$$
for every $\tau \in \Sigma ( 1 )$ and every $n \in \Z$.
\end{definition}

\begin{lemma}\label{lambda indep}
$\Xi ^{\prime} _{\Lambda}$ is independent of the choice of $\Lambda$.
\end{lemma}

\begin{proof}
Let $\Lambda _1 = \{ \lambda ^{\prime} _1, \lambda ^{\prime} _2, \ldots \}$ and $\Lambda _2 = \{ \lambda _{1}, \lambda _2, \ldots \}$ be two choices of $\Lambda$. By changing the numbers of the elements of $\Lambda _1$, we can assume $\lambda _i - \lambda ^{\prime} _i \in X ^* ( T )$ for every $1 \le i \le h$. We put $\Lambda ^i := \{ \lambda _1, \lambda _2, \ldots, \lambda _i , \lambda ^{\prime} _{i + 1}, \ldots\}$. Then, it suffices to check the assertion $\Xi ^{\prime} _{\Lambda ^i} = \Xi ^{\prime} _{\Lambda ^{i + 1}}$ for an arbitrary $0 \le i \le h$. For each $\lambda \in X ^* \left( T \right)$, we have $\left< \tau, \lambda \right> \in \Z$ for every $\tau \in \Sigma ( 1 )$. Hence, $D
^{\lambda ^{\prime} _i} - D ^{\lambda _i} = \sum _{\tau \in \Sigma ( 1 )} \left< \tau, \lambda ^{\prime} _i - \lambda _i \right> D _{\tau}$ is an equivariant divisor by Lemma \ref{divisor}. As a consequence, we have
\begin{align*}
& {} ^{\prime} F ^{\tau} \left( n + \lfloor \left< \tau, \lambda _i \right> \rfloor,
{\mathcal E} ( D ^{\lambda _i} )
\right)\\
& = {}^{\prime} F ^{\tau} \left( n + \lfloor \left< \tau, \lambda
 _i \right> \rfloor
+ \left< \tau, \lambda ^{\prime} _i - \lambda _i \right>,
{\mathcal E} ( D ^{\lambda _i} + \left< \tau,
\lambda ^{\prime} _i  - \lambda _i \right> D _{\tau} )
\right)\\
& = {}^{\prime} F ^{\tau} \left( n + \lfloor \left< \tau, \lambda _i
 \right> \rfloor
+ \left< \tau, \lambda ^{\prime} _i - \lambda _i \right>,
{\mathcal E} ( D ^{\lambda _i} + D ^{\lambda
^{\prime} _i - \lambda _i} )
\right)\\
& = {}^{\prime} F ^{\tau} \left( n + \lfloor \left< \tau, \lambda
^{\prime} _i \right> \rfloor,
{\mathcal E} ( D ^{\lambda ^{\prime} _i} )
\right)
\end{align*}
for every $\tau \in \Sigma ( 1 )$ by the construction of Definition \ref{general redef}.
Therefore, Definition \ref{general redef} does not depend on the choice of $\lambda _i$.
\end{proof}

Now we come to one of the main result of this subsection.

\begin{proposition}\label{identification}
For each ${\mathcal E} \in {\bf Ob} EV \left( \Sigma, 1 \right) _{c}$, we have $\Xi ^{\prime} \left( {\mathcal E} \right) \cong \Xi \left( {\mathcal E} \right)$ as vector spaces equipped with families of filtrations.
\end{proposition}

\begin{proof}
Let $\tau \in \Sigma ( 1 )$. The definition of $F ^{\tau} \left( \bullet, {\mathcal E} \right)$ (\S
\ref{CatEV}) makes sense if we restrict ourselves to $T \left( \tau
\right)$. Thus, we restrict our attention to $T \left( \tau
\right)$. ${\mathcal E}
\MID_{T \left( \tau \right)}$ is isomorphic to $E \times T \left(
\tau \right)$ for some $T$-module $E$ by Theorem \ref{Klyachko's Proposition} 1). Let $E =
\oplus _{\lambda \in X ^* ( T )} E ^{\lambda}$ be the isotypical decomposition of $E$.

Let $E _0$ be a trivial $\mathbb{G} ^{\tau} _m$-representation such
that $\dim E = \dim E _0$. By Theorem \ref{Klyachko's Proposition} 3), we have $F ^{\tau}
( n, {\mathcal E}) = \bigoplus _{\left< \tau, \lambda \right> + n \le 0} E
^{\lambda}$. We write $T \left( \tau \right) = \Spec R$, where $R$ is the $( 1 \times T )$-algebra $\bigoplus _{\mu \in X ^* ( T ) ; \left< \tau, \mu \right> \le 0} \mu$. Then, Definition \ref{main functor for trivial
factor} turns into the following:  Let $\phi ^{\prime }_{\tau} :
R \rightarrow k ( x _{\tau} )$ and $\phi ^{\prime } _{e} : R
\rightarrow k ( e )$ be the residual maps at $x _{\tau}$ and $e$, respectively. We
define $\phi _{\tau} := \id \otimes \phi ^{\prime} _{\tau} : E _0
\otimes R \rightarrow E _0 \otimes k ( x _{\tau} )$ and $\phi _{e} := \id \otimes \phi ^{\prime} _{e} : E _0 \otimes R
\rightarrow E _0 \otimes k ( x _{e} )$. We have an isomorphism $E
\otimes k ( x _{\tau} ) \cong E
\otimes k ( e )$ induced from the constant section $E _0 \otimes k
\subset E _0 \otimes \bigoplus _{\mu \in X ^* ( T ) ; \left< \tau, \mu
\right> \le 0} \mu = E _0 \otimes R$. For a sufficiently large integer $N$, we have
$E \otimes R \subset E _0 \otimes \bigoplus _{\mu \in X ^* ( T ) ; \left< \tau, \mu
\right> \le N} \mu$ as compatible $\left( 1 \times T, R \right)$-modules. (For any one-dimensional $T$-module $\mu _0$ such that $\left< \tau, \mu _0 \right> = N$, the RHS is isomorphic to $E _0 \otimes \mu _0 \otimes R$.) Thus, we have
\begin{eqnarray}
E \otimes R \hookrightarrow \left( E _0 \otimes \bigoplus _{\mu \in X ^* ( T ) ; \left< \tau, \mu
\right> \le N} \mu \right) \hookleftarrow E _0 \otimes R\label{inclusions}
\end{eqnarray}
as compatible $( 1 \times T, R )$-modules. Since the image of $\mu$ under $\phi ^{\prime} _e$ is $1$ for every $\mu \in X ^* ( T )$, (\ref{inclusions}) induces isomorphisms
$$E \otimes R \otimes _R k ( e ) \cong \left( E _0 \otimes \bigoplus _{\mu \in X ^* ( T ) ; \left< \tau, \mu
\right> \le N} \mu \right) \otimes _R k ( e ) \cong E _0 \otimes R
\otimes _R k ( e ) \cong  E _0 \otimes k ( e ).$$
Hence, we have $F ^{\tau}
( n, {\mathcal E}) \subset E \otimes R \otimes _R k ( e ) \cong E _0 \otimes R
\otimes _R k ( e )$ for every $n \in \Z$.
In this setting, we have
\begin{eqnarray*}
{} ^{\prime} F ^{\tau} ( n, {\mathcal E}) &=& \phi _{\tau} \left(
 E _0
\otimes R \cap \bigoplus _{\mu , \rho \in X ^* ( T ) ; \left< \tau,
\mu \right> \le - n} E ^{\rho} \otimes \mu \right) \\
&=& \bigoplus _{\mu , \rho \in X ^* ( T ) ; \left< \tau,
\mu \right> \le - n} \left( E _0 \otimes k \cap E ^{\rho} \otimes \mu \right) \subset E _0 \otimes k ( x _{\tau} ) \cong E _0 \otimes k ( e ).
\end{eqnarray*}
By comparing term by term, we conclude ${} ^{\prime} F ^{\tau} ( n, {\mathcal E}) = F ^{\tau} ( n, {\mathcal E})$ as subspaces of the identity fiber $E _0 \otimes k ( e )$ for each $\tau$.
\end{proof}

Thanks to Proposition \ref{identification}, we do not need to distinguish $\Xi$ and $\Xi ^{\prime}$. Thus, we consider $\Xi ^{\prime} _{\Lambda}$ in Definition \ref{general redef} to be the primary definition of $\Xi$ hereafter because it prolongs the original $\Xi$ in \S \ref{CatEV} and is independent of the choice of $\Lambda$ by Lemma \ref{lambda indep}.

\subsubsection{Making $\Xi$ into a functor}
In the definition of $\Xi$ in \S \ref{CatEV}, we only need the closure of $( 1 \times T ) e$ (inside $X ( \Sigma )$) to define $F ^{\tau}$. As a consequence, our definition of $\Xi$ coincides with the composition of the restriction to $T \left( \Sigma \right)$ and Klyachko's functor \cite{Kl} 0.1. Therefore, we obtain the following from Klyachko's description \cite{Kl} 2.3 by using Definition \ref{general redef} and Proposition \ref{identification}.

\begin{corollary}\label{saturatedness}
Let ${\mathcal E} \in {\bf Ob} EV \left( \Sigma \right) _{c}$. For each $\tau \in \Sigma ( 1 )$, $F ^{\tau} \left( \bullet, {\mathcal E} \right)$ is a decreasing filtration of $B ( {\mathcal E} )$. Moreover, we have $F ^{\tau} \left( - n, {\mathcal E} \right) = B ( {\mathcal E} )$ and $F ^{\tau} \left( n, {\mathcal E} \right) = 0$ for $n >> 0$.
\end{corollary}

\begin{corollary}\label{distributiveness}
Let ${\mathcal E} \in {\bf Ob} EV \left( \Sigma \right) _{c}$. For each $\sigma \in \Sigma$, $\{ F ^{\tau} \left( n, {\mathcal E} \right) \} _{\tau \in \sigma ( 1 ), n \in \Z}$ forms a distributive lattice.
\end{corollary}

\begin{lemma}\label{module}
Let ${\mathcal E} \in {\bf Ob} EV \left( \Sigma \right) _{c}$. For each $\tau \in \Sigma( 1 )$ and $n \in \Z$, $F ^{\tau} \left( n, {\mathcal E} \right)$ is a ${\tilde P} ^{\tau} \times \tilde{Z} ( G )$-submodule of $B ( {\mathcal E} )$.
\end{lemma}

\begin{proof}
By Definition \ref{general redef}, it suffices to prove the assertion for ${\bf Ob} EV \left( \Sigma, 1 \right) _{c}$. By Lemma \ref{intersection of equivariant vector bundles}, ${\mathcal E} \left( - n H_{\tau} \right) \cap _D B ( \mathcal E ) \otimes {\mathcal O}_X$ is a $\tilde{G} \times G$-equivariant ${\mathcal O}_X$-module for a sufficiently large equivariant divisor $H$. Thus, its fiber at $x_{\tau}$ admits an action of the stabilizer at $x _{\tau}$. Hence, $F^{\tau} \left( n, {\mathcal E} \right)$ is a $\tilde{G} ^{\tau}$-submodule of $B \left( {\mathcal E}\right)$. By the definition of $\triangle ^d$, we have $\tilde{P} ^{\tau} = \tilde{L} ^{\tau} U ^{\tau} _+ \hookrightarrow \triangle ^d \tilde{L} ^{\tau} ( U ^{\tau} _+ \times 1 ) \subset \tilde{G} ^{\tau}$. By Lemma \ref{switch}, the $\tilde{P} ^{\tau} = \tilde{L} ^{\tau} U ^{\tau} _+$-action on $B \left( {\mathcal E} \right) ( \cong B ( \mathcal E ) \otimes {\mathcal O} _X \otimes _X k ( x _{\tau} ) )$ coincides with the action of stabilizer $\triangle ^d \tilde{P} ^{\tau}$ on $B \left( {\mathcal E} \right) ( \cong B ( \mathcal E ) \otimes {\mathcal O}_X \otimes _X k ( e ))$ (notice that we have fixed an isomorphism $k ( e ) \cong k ( x _{\tau} )$ in Definition \ref{main functor for trivial factor}).
\end{proof}

\begin{lemma}\label{hom of naive functor}
Let ${\mathcal E} _1, {\mathcal E} _2 \in {\bf Ob} EV \left( \Sigma \right) _{c}$ and let $f: {\mathcal E} _1 \rightarrow {\mathcal E} _2$ be a morphism of $EV \left( \Sigma \right) _{c}$. Then, we have a $\tilde{P} ^{\tau} \times \tilde{Z} ( G )$-module morphism $f ^{\prime} : F ^{\tau} \left( n, {\mathcal E} _1 \right) \rightarrow F ^{\tau} \left( n, {\mathcal E} _2 \right)$ such that the morphism of their ambient spaces $f ^{\prime} : B ( {\mathcal E} _1 ) \rightarrow B ( {\mathcal E} _2 )$ is obtained by the specialization of $f$ at $e$ for every $\tau \in \Sigma ( 1 )$ and every $n \in \Z$.
\end{lemma}

\begin{proof}
By Definition \ref{general redef} and Corollary \ref{isotypical decomposition of category of vector bundles}, it suffices to prove the assertion for ${\bf Ob} EV \left( \Sigma, 1 \right) _{c}$. 
We have a natural morphism $f ^{\prime} : B ( {\mathcal E} _1 ) \rightarrow B ( {\mathcal E} _2 )$ induced by $f$. We have $B ( {\mathcal E} _i ) = \Gamma \left( G, \left( {\mathcal E} _i \right) _1 \right)^{1 \times \tilde{G}}$ for $i = 1, 2$. Let $D$ be a sufficiently large divisor. By Lemma \ref{incl} and Corollary \ref{uniq-incl}, we have a $\tilde{G} \times \tilde{G}$-equivariant inclusion $\Gamma \left( G, {\mathcal E} _1 \right)^{1 \times \tilde{G}} \otimes {\mathcal O}_X \left( - D \right) \subset {\mathcal E} _1$ such that $f$ induces a morphism $\Gamma \left( G, {\mathcal E} _1 \right)^{1 \times \tilde{G}} \otimes {\mathcal O}_X \left( -D \right) \stackrel{f ^{\prime} \otimes \mathrm{id}}{\rightarrow} \Gamma \left( G, {\mathcal E} _2 \right)^{1 \times \tilde{G}} \otimes {\mathcal O} _X \left( - D \right) \subset {\mathcal E} _2$. Enlarging $D$ if necessary, we have the following commutative diagram of $\tilde{G} \times \tilde{G}$-equivariant coherent sheaves by Corollary \ref{sandwich}.
\begin{eqnarray*}
\begin{matrix}
{\mathcal E} _1 \left( - n D _{\tau} \right) & \stackrel{f \otimes \mathrm{id}}{\rightarrow} & {\mathcal E} _2 \left( - n D _{\tau} \right)\\
\bigcap & & \bigcap\\
\mathcal B ( {\mathcal E} _1 ) ^D & \stackrel{f ^{\prime} \otimes \mathrm{id}}{\rightarrow} & \mathcal B ( {\mathcal E} _1 ) ^D\\
\bigcup & & \bigcup\\
B ( {\mathcal E} _1 ) \otimes {\mathcal O} _X & \stackrel{f ^{\prime} \otimes \mathrm{id}}{ \rightarrow} & B ( {\mathcal E} _1 ) \otimes {\mathcal O} _X
\end{matrix}
\end{eqnarray*}
Here we can take the intersection of ${\mathcal E} _1 \left( - n D _{\tau} \right) \stackrel{f \otimes \mathrm{id}}{\rightarrow} {\mathcal E} _2 \left( - n D _{\tau} \right)$ and $B ( {\mathcal E} _1 ) \otimes {\mathcal O} _X \stackrel{f ^{\prime} \otimes \mathrm{id}}{\rightarrow} B ( {\mathcal E} _2 ) \otimes {\mathcal O} _X$ in $\mathcal B ( {\mathcal E} _1 ) ^D \stackrel{f ^{\prime} \otimes \mathrm{id}}{\rightarrow} \mathcal B ( {\mathcal E} _2 ) ^D$ termwise to define a morphism
$${\mathcal E} _1 \left( - n D _{\tau} \right) \cap _D B ( {\mathcal E} _1 ) \otimes {\mathcal O} _X \rightarrow {\mathcal E} _2 \left( - n D _{\tau} \right) \cap _D B ( {\mathcal E} _2 ) \otimes {\mathcal O} _X.$$
Applying $\phi _{\tau}$ (in Definition \ref{main functor for trivial factor}) to both sides, we obtain a $\tilde{G} ^{\tau}$-module morphism $F ^{\tau} \left( n, {\mathcal E} _1 \right) \rightarrow F ^{\tau} \left( n, {\mathcal E} _2 \right)$ because all sheaves and morphisms in the above diagram are $\tilde{G} \times \tilde{G}$-equivariant. Hence, we obtain the result from the same argument as in the proof of Lemma \ref{module}.
\end{proof}

\begin{proposition}\label{naive functor}
$\Xi$ defines a faithful covariant functor $\Xi : EV \left( \Sigma \right) _{c} \rightarrow {\mathfrak C} \left( \Sigma \right) ^l _{c}$.
\end{proposition}

\begin{proof}
By Corollary \ref{saturatedness}, Corollary \ref{distributiveness} and Lemma \ref{module}, the range of the map $\Xi$ is contained in ${\bf Ob} {\mathfrak C} \left( \Sigma \right) ^l _{c}$. By Lemma
\ref{hom of naive functor}, we have a morphism (of sets)
$$\Xi : \Hom _{EV \left( \Sigma \right) _{c}} \left( {\mathcal E}, {\mathcal F} \right) \rightarrow \Hom _{{\mathfrak C} \left( \Sigma \right) ^l _{c}} \left( \Xi \left( {\mathcal E} \right), \Xi \left( {\mathcal F} \right)\right)$$
for every ${\mathcal E}, {\mathcal F} \in {\bf Ob} EV \left( \Sigma \right) _{c}$. By the similar argument as in the proof of Lemma \ref{hom of naive functor}, we can deduce that $\Xi \left( f \right) \circ \Xi \left( g \right) = \Xi \left( f \circ g \right)$ for two morphisms $f, g$ in $EV \left( \Sigma \right) _{c}$. Thus, $\Xi$ becomes a covariant functor from $EV
\left( \Sigma \right) _{c}$ to ${\mathfrak C} \left( \Sigma \right) ^l _{c}$.
For each nonzero morphism of $EV \left( \Sigma
\right) _{c}$, its restriction to $G$ is nontrivial since a locally
free ${\mathcal O} _X$-module cannot contain torsion submodule.
Then, it induces a nontrivial morphism between the identity
fibers by Lemma \ref{Basic Situation}. Hence, $\Xi$ must be faithful.
\end{proof}
\section{Actions and toric slice}\label{toric slice}
\subsection{Toric decomposition}\label{infinitesimal}
Here we develop a method to study of $\tilde{G} \times G$-equivariant vector bundles via the study of equivariant vector bundles of a certain one-dimensional toric variety.

\subsubsection{Reduction to $\mathbb A ^1$}\label{2red}
First, we fix an arbitrary one-dimensional cone $\tau$ of $\Sigma$.

We define $\mathbb{A} _{\tau} := \overline{\mathbb{G} ^{\tau} _m e} \subset T \left( \Sigma \right) \backslash \{ x _{- \tau} \}$. We have $\mathbb A _{\tau} = \mathbb G _m ^{\tau} \cup \{ x _{\tau} \}$ (as sets) and $\mathbb A _{\tau} \cong \mathbb A ^1$ (as $\mathbb G _m ^{\tau}$-toric varieties). By Theorem \ref{structure theorem}, we have an open embedding $B ^- . \mathbb A _{\tau} . B \hookrightarrow X ( \Sigma )$. Hence, we have the following composition of faithful functors
$$\mathrm{Res} ^{\tau} : \mathfrak{Coh} ^{\tilde{G} \times G} X ( \Sigma ) \rightarrow \mathfrak{Coh} ^{\tilde{B} ^- \times B} ( B ^- . \mathbb A _{\tau} . B ) \rightarrow \mathfrak{Coh} ^{\mathbb G _m ^{\tau}} \mathbb A _{\tau},$$
where the second functor is the composition of the restriction (cf. \cite{CG} 5.2.16) and restriction of the group. Since all of the previous constructions use faithful functors, we can restrict the bounding sheaves of Corollary \ref{sandwich} to $\mathbb A _{\tau}$ by $\mathrm{Res} ^{\tau}$. In particular, the definition of $F ^{\tau} ( \bullet, {\mathcal E} )$ (see Definition \ref{naive Xi}) factors through $\mathrm{Res} ^{\tau}$. Let $t$ be a local uniformizer of $x _{\tau}$ on $\mathbb{A} _{\tau}$ which is an eigenfunction with respect to $\mathbb{G} ^{\tau} _m$. 

\begin{lemma}\label{small sandwich}
For each ${\mathcal E} \in {\bf Ob} EV \left( \Sigma, 1 \right) _{c}$ and a sufficiently large integer $N$, we have the following inclusions of $( \triangle ^d ( \tilde{L} ^{\tau}) \mathbb{G} ^{\tau} _m, k[t] )$-modules:

$$B( {\mathcal E} ) \otimes t ^N k[t] \subset \Gamma \left( \mathbb{A} _{\tau}, {\mathcal E} \MID _{\mathbb{A} _{\tau}} \right) \subset B( {\mathcal E} ) \otimes t ^{- N} k[t] = \Gamma ( \mathbb{A} _{\tau}, {\mathcal B} \left( {\mathcal E} \right)^{ N D _{\tau}}  \MID _{\mathbb{A} _{\tau}} ).$$
Here $\triangle ^d ( \tilde{L} ^{\tau} )$ acts on $B( {\mathcal E} )$ as a subgroup of $\tilde{G}$ via $\tilde{L} ^{\tau} \subset \tilde{G}$ and acts on $t$ trivially. $\mathbb{G} ^{\tau} _m$ acts on $B( {\mathcal E} )$ trivially and acts on $t$ by weight one.
\end{lemma}

\begin{proof}
By Corollary \ref{sandwich}, all we have to do is to check the compatibility with respect to the group
actions. Since $\triangle ^d ( \tilde{L} ^{\tau} )$ fixes ${\mathbb{A} _{\tau}}$ pointwise, it acts on $t$ trivially
(cf. \S \ref{NoteVar}). By the definition of $B( {\mathcal E} ) \otimes {\mathcal
O} _X$ (\S \ref{lemma on equivariant str}), $\triangle ^d (
\tilde{L} ^{\tau} )$ acts on $B( {\mathcal E} )$ as $\tilde{L} ^{\tau} \subset \tilde{G}$ and $\mathbb{G} ^{\tau} _m$ acts on $B( {\mathcal E} )$ trivially. Since $x _{\tau}$ is a limit of $t .e$ with $t \rightarrow 0$ in $\mathbb{G} ^{\tau} _m$, $\mathbb{G} ^{\tau} _m$ acts on $t$ by degree one.
\end{proof}

\begin{corollary}\label{cutting down}
Let ${\mathcal E} \in {\bf Ob} EV \left( \Sigma, 1 \right) _{c}$. Let $N$ be an integer. Then, we have
$$\Gamma ( \mathbb{A} _{\tau}, ( {\mathcal B} \left( {\mathcal E} \right) ^{N D _{\tau}} \cap _D {\mathcal E} ) \MID
_{\mathbb{A} _{\tau}} ) = \bigoplus _{n = - \infty}
^{N} F ^{\tau} \left( n, {\mathcal E} \right) \otimes k t ^{- n}.$$
as $\triangle ^d ( \tilde{L} ^{\tau} ) \mathbb{G} _m^{\tau}$-submodules of $B \left(
{\mathcal E} \right) \otimes t ^{- N} k[t]$ for a sufficiently large equivariant divisor $D$. In
particular, the RHS is the $\mathbb{G}
^{\tau} _m$-isotypical decomposition.
\end{corollary}

\begin{proof}
By Lemma \ref{small sandwich}, $\mathbb{G} _m
^{\tau}$ acts on $B( {\mathcal E} ) \otimes k t ^n \subset B ( {\mathcal
E} ) \otimes t ^{- N} k [t] = \Gamma ( \mathbb{A} _{\tau}, {\mathcal B} \left(
{\mathcal E} \right) ^{N D _{\tau}} \MID
_{\mathbb{A} _{\tau}} )$ by degree $n$. Thus, we have the $\mathbb{G}
^{\tau} _m$-isotypical decomposition
$$\Gamma \left( \mathbb{A}
_{\tau}, {\mathcal B} ( {\mathcal E}
\right) ^{N D _{\tau}} ) \cong \bigoplus _{n = - \infty}
^{N} B \left( {\mathcal E} \right) \otimes k t ^{- n}.$$
By Corollary \ref{sandwich}, we have
$$\Gamma \left( \mathbb{A} _{\tau}, 
{\mathcal E} \MID _{\mathbb{A} _{\tau}} \right) \subset \Gamma ( \mathbb{A} _{\tau}, {\mathcal B} \left( {\mathcal E} \right) ^{M D _{\tau}} \MID _{\mathbb{A} _{\tau}} )$$
as a compatible $( \mathbb{G} _m^{\tau}, k[t] )$-submodule for $M >>
0$. Hence, we have the $\mathbb{G} ^{\tau} _m$-isotypical decomposition
$$\Gamma \left( \mathbb{A} _{\tau}, 
{\mathcal E} \MID _{\mathbb{A} _{\tau}} \right) = \bigoplus _{n = - \infty}
^{\infty} {} ^{\prime} F ^{\tau} \left( n, {\mathcal E} \right) \otimes k t ^{- n}.$$
Since the latter is $k [ t ]$-free, ${} ^{\prime} F ^{\tau} \left( \bullet, {\mathcal E} \right)$ is a decreasing filtration of $B ( {\mathcal E} )$. If $\tau \neq \xi \in \Sigma ( 1 )$, then twisting by $D _{\xi}$ has no effect on coherent shaves on $\mathbb{A} _{\tau}$. Moreover, twisting $D _{\tau}$ is equivalent to multiplying $t ^{- 1}$. Thus, we have
\begin{align*}
F ^{\tau} ( n, {\mathcal E} )  & = \phi _{\tau} \left( {\mathcal E} \left( - n D _{\tau} \right) \cap _D \mathcal B ( {\mathcal E} ) \right)\\
& = \Gamma \left( \mathbb{A} _{\tau}, {\mathcal E} \left( - n D _{\tau} \right) \MID _{\mathbb{A} _{\tau}} \right) \cap B ( {\mathcal E} ) \otimes k t ^0 = {} ^{\prime} F ^{\tau} \left( n, {\mathcal E} \right)
\end{align*}
as vector subspaces of $B ( {\mathcal E} )$ by Definition \ref{main functor for trivial factor}. (Here we put $D \MID _{\mathbb{A} _{\tau}} = M D _{\tau} \MID _{\mathbb{A} _{\tau}}$.) Therefore, we have
$$\Gamma \left( \mathbb{A} _{\tau}, ( {\mathcal B} \left( {\mathcal E} \right) ^{N D _{\tau}} \cap _D {\mathcal E} ) \MID
_{\mathbb{A} _{\tau}} \right) = \bigoplus _{n = - \infty}
^{N} F ^{\tau} \left( n, {\mathcal E} \right) \otimes k t ^{- n} \subset B ( {\mathcal E} ) \otimes t ^{- N} k [ t ]$$
as compatible $( \mathbb{G} ^{\tau} _m, k [t] )$-modules.
By Lemma \ref{module}, each inclusion $F ^{\tau} \left( n, {\mathcal E} \right) \hookrightarrow B (
{\mathcal E} )$ is a $\triangle ^d ( \tilde{L} ^{\tau} )$-module
inclusion.
\end{proof}

For a $\tilde{G}$-module $V$, we define the formal loop space $V _{[t]} ^{\tau}$ of $V$ with respect to $\tau \in \Sigma ( 1 )$ as follows:
\begin{eqnarray}
V _{[t]} ^{\tau} := \mathop{\varinjlim} _N \Gamma \left( \mathbb{A} _{\tau}, V \otimes {\mathcal O} _X \left( N D _{\tau} \right) \MID _{\mathbb{A} _{\tau}} \right) = V \otimes \mathop{\varinjlim} _N t ^{- N} k [ t ] \cong V \otimes k [ t, t ^{- 1} ] \label{loopspace}
\end{eqnarray}
Here $t, t ^{-1}$, and $V \otimes k$ generates $V _{[t]} ^{\tau}$.

\begin{definition}
We define the category ${\mathfrak B} \left( \Sigma \right) ^{\tau}$ as follows:
\begin{itemize}
\item $\bf{(Objects)}$ Triples $\left( R, V, \iota \right)$ such that:
\begin{itemize}
\item $R$ is a compatible $( \triangle ^d ( \tilde{L} ^{\tau} ) \mathbb{G} ^{\tau} _m, k [t] )$-module;
\item $V$ is a $\tilde{G}$-module;
\item $\iota : R \hookrightarrow V _{[t]} ^{\tau}$ is an injective morphism of compatible $( \triangle ^d ( \tilde{L} ^{\tau} ) \mathbb{G} ^{\tau} _m, k [t] )$-modules.
\end{itemize}
For simplicity, we may also denote it by $( R \subset V _{[t]} ^{\tau})$.
\item $\bf{(Morphisms)}$ We define the morphism by the following
commutative diagram of compatible $( \triangle ^d ( \tilde{L} ^{\tau} ) \mathbb{G} ^{\tau} _m, k [t] )$-modules:
$$
\begin{matrix}
( R _1 & \hookrightarrow & ( V _1 ) _{[t]} ^{\tau} ) & \in & {\bf Ob} {\mathfrak B} \left( \Sigma \right) ^{\tau}\\
\downarrow & & \downarrow & & \\
( R _2 & \hookrightarrow & ( V _2 ) _{[t]} ^{\tau} ) & \in & {\bf Ob} {\mathfrak B} \left( \Sigma \right) ^{\tau}
\end{matrix}
$$
Here $( V _1 ) _{[t]} ^{\tau} \rightarrow ( V _2 ) _{[t]} ^{\tau}$ is the morphism induced by a $\tilde{G}$-module map $V _1 \rightarrow V _2$.
\end{itemize}
\end{definition}

For each ${\mathcal E} \in {\bf Ob} EV \left( \Sigma, 1 \right) _c$, we define an object $\rho ^{\tau} _{\infty} \left( {\mathcal E}
\right) \in {\bf Ob} {\mathfrak B} \left( \Sigma \right) ^{\tau}$ as follows:
$$\rho ^{\tau} _{\infty} \left( {\mathcal E} \right) := \left( \Gamma \left(
\mathbb{A} _{\tau},{\mathcal E} \MID _{ \mathbb{A} _{\tau}} \right) = \mathop{\varinjlim} _N \Gamma \left( \mathbb{A} _{\tau}, ( {\mathcal E} \cap _D {\mathcal B} \left( {\mathcal E} \right) ^{N D _{\tau}} ) \MID _{ \mathbb{A} _{\tau}} \right) \subset B \left( {\mathcal E} \right) _{[t]} ^{\tau}\right)$$
Here $D$ is a sufficiently large equivariant divisor which depends on
${\mathcal E}$ (cf. Lemma \ref{intersection of equivariant vector bundles}).

\begin{remark}
It is important NOT to forget the ambient space $B \left( {\mathcal E} \right) _{[t]} ^{\tau}$.
\end{remark}

\begin{corollary}\label{def of rho}
For each ${\mathcal E} \in {\bf Ob} EV \left( \Sigma, 1 \right) _{c}$, we have
$$\rho ^{\tau} _{\infty} \left( {\mathcal E} \right) \cong \left( \left( 
\bigoplus _{n \in \Z} F ^{\tau} \left( n, {\mathcal E} \right) \otimes k t ^{n} \right) \subset B \left( {\mathcal E} \right) _{[t]} ^{\tau} \right)$$
as $\triangle ^d ( \tilde{L} ^{\tau} ) \tilde{\mathbb{G}} _m^{\tau}$-modules. In
particular, the above direct sum decomposition is the $\mathbb{G}
^{\tau} _m$-isotypical decomposition.
\end{corollary}

\begin{proof}
The assertion follows from the definition of $\rho ^{\tau} _{\infty}$ and Corollary \ref{cutting down}.
\end{proof}

\begin{corollary}\label{prenaive functor}
$\rho ^{\tau} _{\infty}$ defines a functor $EV \left( \Sigma, 1 \right) _{c} \rightarrow {\mathfrak B} \left( \Sigma \right) ^{\tau}$.
\end{corollary}

\begin{proof}
This follows immediately from Proposition \ref{naive functor} and Corollary \ref{def of rho}.
\end{proof}

\subsubsection{Basic computation}
We work under the same settings as in \S \ref{2red}. In particular, $\tau \in \Sigma ( 1 )$. Here, we compute $\rho ^{\tau} _{\infty}$ for two kinds of
$\tilde{G} \times G$-equivariant vector bundles. The description of the first series of
$\tilde{G} \times G$-equivariant vector bundles is a direct consequence of the definition.

\begin{corollary}\label{def of V}
For a $\tilde{G}$-module $V$, we have $\rho ^{\tau} _{\infty} \left( V \otimes {\mathcal O} _X \right) = \left( V \otimes k[ t ] \subset V ^{\tau} _{[t]} \right)$.
\end{corollary}

\begin{proof}
It is a direct consequence of Corollary \ref{def of rho}.
\end{proof}

To describe our second series of $\tilde{G} \times G$-equivariant vector bundles, we need a preparation.

\begin{definition}\label{max}
Let $\lambda$ be a dominant weight. Let $v _{w _0 \lambda} ^*$ be a highest weight vector of $V _{- w _0 \lambda}$ (cf. \S \ref{NoteAlg}). We define
$$V _{\lambda} ^{\tau} \left( n \right) := \{v \in V _{\lambda} ; \tau ( t ) ( v
\otimes v _{w _0 \lambda} ^* ) = t ^n v
\otimes v _{w _0 \lambda} ^* \mbox{ for all } t \in \mathbb{G} _m ( k ) \subset k \}$$
for every $n \in \Z$. We have a direct sum decomposition $V _{\lambda} = \oplus _{n \ge 0} V _{\lambda} ^{\tau} \left( n \right)$. Moreover, we define 
$$F ^{\tau} _{max} \left( m, V _{\lambda} \right) := \bigoplus _{n \ge m} V _{\lambda} ^{\tau} \left( n \right) \text{ for all } m \in \Z.$$
We call this filtration the maximal filtration of the $\tilde{G} \times \tilde{Z} ( G )$-module $V _{\lambda} \boxtimes k$.

Let $V$ be a finite dimensional rational $\tilde{G}$-module. We denote its ($\tilde{G}$-) irreducible decomposition by $\oplus _{\mu} V _{\mu} ^{\oplus n _{\mu}}$. Then, we define
$$F ^{\tau} _{max} \left( m, V \right) := \bigoplus _{\mu} F ^{\tau}
_{max} \left( m, V _{\mu} \right) ^{\oplus n _{\mu}} \text{ for all } m \in \Z.$$
We call this filtration the maximal filtration of the $\tilde{G} \times \tilde{Z} ( G )$-module $V \boxtimes k$. It is independent of the choice of an irreducible decomposition of $V$.
\end{definition}

Taking account the weight decomposition of $V$ and Definition
\ref{distributive lattice}, we see that $F ^{\tau} _{max} \left( 
\bullet, V \right)$ is an $\tau$-transversal filtration.

Now we can state and prove the description of our second series of
$\tilde{G} \times G$-equivariant vector bundles.

\begin{proposition}\label{maximal}
We have $\Xi \left( {\mathcal L} _{- w _0 \lambda}
\otimes V _{\lambda} \right) = \left( V _{\lambda}, \{ F ^{\tau}
_{max} \left( \bullet, V _{\lambda} \right)\} _{\tau \in \Sigma ( 1 )} \right)$ for every dominant weight $\lambda$. $($Here ${\mathcal L} _{\lambda}$ is defined before Lemma \ref{divisor}.$)$
\end{proposition}

\begin{proof}
We have ${\mathcal L} _{- w _0 \lambda} \in {\bf Ob} EV \left( \Sigma, - w _0 \bar{\lambda} \right) _c$ by Theorem \ref{Strickland} and the definition. Thus, we have ${\mathcal L} _{- w _0 \lambda} \otimes V _{\lambda} \in {\bf Ob} EV \left( \Sigma, 1 \right) _c$.
We compute the boundary behavior as in \S \ref{CatEV} via Proposition \ref{identification} for each direction. We choose an arbitrary $\tau \in \Sigma ( 1 )$. By the remarks at the beginning of \S \ref{infinitesimal}.1, we restrict ourselves to $\left( {\mathcal L} _{- w _0 \lambda} \otimes V _{\lambda} \right) \MID _{\mathbb{A} _{\tau}}$.
By Theorem \ref{Klyachko's Proposition} 1), we have
$$\left( {\mathcal L} _{- w _0 \lambda} \otimes V _{\lambda} \right) \MID _{\mathbb{A} _{\tau}} \cong \mathbb{A} _{\tau} \times \left( ( w _0 \lambda ) ^{- 1} \otimes V _{\lambda} \right)$$
as $\mathbb{G} ^{\tau} _m$-equivariant vector bundles on $\mathbb{A} _{\tau}$. Therefore, we have 
$$F ^{\tau} \left( m, \mathbb{A} _{\tau} \times \left( ( w _0 \lambda
) ^{- 1} \otimes V _{\lambda} \right) \right) = \bigoplus _{n \ge m} V
_{\lambda} ^{\tau} \left( n \right) = F ^{\tau} _{max} \left( m, V _{\lambda} \right)$$
for every $m \in \Z$ by Theorem \ref{Klyachko's Proposition} 3).
\end{proof}

\subsubsection{Injectivity of $\Xi$}
The rest of this subsection is devoted to the proof of Proposition \ref{last}.

\begin{proposition}\label{last}
$\Xi$ induces an injective map ${\bf Ob} EV \left( \Sigma
\right) _{c} \hookrightarrow {\bf Ob} {\mathfrak C}  \left( \Sigma
\right) _{c} ^l$.
\end{proposition}

First, we prove a simple special case.

\begin{lemma}\label{plast}
Let $\tau \in \Sigma ( 1 )$. Let ${\mathcal E} _1, {\mathcal E} _2 \in {\bf Ob} EV \left( \tau, 1
\right) _{c}$ be such that $\Xi \left( {\mathcal E}
_1 \right) \cong \Xi \left( {\mathcal E} _2 \right)$. For a
sufficiently large equivariant divisor $D$, we fix inclusions ${\mathcal B} \left(
{\mathcal E} _1 \right) ^{- D} \subset {\mathcal E} _1, {\mathcal E} _2
\subset {\mathcal B} \left(
{\mathcal E} _1 \right) ^D$ which gives rise
to an isomorphism $\Xi \left( {\mathcal E}
_1 \right) \cong \Xi \left( {\mathcal E} _2 \right)$ $($see Corollary
\ref{sandwich}$)$. Then, we have ${\mathcal E} _1 = {\mathcal E} _2$ as a $\tilde{G} \times G$-equivariant coherent subsheaves of ${\mathcal B} \left(
{\mathcal E} _1 \right) ^D$.
\end{lemma}

\begin{proof}
We prove by induction on $p$. We have ${\mathcal
E} _1 \cap _{H} {\mathcal B} \left(
{\mathcal E} _1 \right) ^{p D _{\tau}}
 = {\mathcal
E} _2 \cap _{H} {\mathcal B} \left(
{\mathcal E} _1 \right) ^{p D _{\tau}}
$ for every sufficiently
small integer $p$. By assumption, we have $\overline{\mathbf O} _{\tau} = \mathbf O _{\tau}$. Hence, every $\tilde{G} \times G$-equivariant ${\mathcal O} _{X ( \tau )}$-module annihilated by ${\mathcal O} _{X ( \tau )} \left( - D _{\tau} \right)$ is a $\tilde{G} \times G$-equivariant vector bundle on $\mathbf O _{\tau}$. In particular, it is determined by its fiber at $x _{\tau}$. Suppose
$${\mathcal E} _1 \cap _D {\mathcal B} \left(
{\mathcal E} _1 \right) ^{p D _{\tau}}
= {\mathcal
E} _2 \cap _D {\mathcal B} \left(
{\mathcal E} _1 \right) ^{p D _{\tau}}
$$
and $F ^{\tau} \left( p + 1, {\mathcal E} _1 \right) = F
^{\tau} \left( p + 1, {\mathcal E} _2 \right)$ as a subspace of $B
\left( {\mathcal E} _1 \right)$ for some integer $p$. Let ${\mathcal
Q} _{p + 1} \left( {\mathcal E} _1 \right) \tilde{}$ be a $\tilde{G} \times
G$-equivariant ${\mathcal O} _{\mathbf O _{\tau}}$-submodule of $B ( {\mathcal
E} _1 ) \otimes {\mathcal O} _{\mathbf O _{\tau}}$
such that ${\mathcal
Q} _{p + 1} \left( {\mathcal E} _1 \right) \tilde{} \otimes _X k ( x _{\tau} ) \cong F ^{\tau} \left( p + 1, {\mathcal E} _1 \right)$ as $\tilde{G} ^{\tau}$-modules. We put ${\mathcal
Q} _{p + 1} \left( {\mathcal E} _1 \right) := {\mathcal
Q} _{p + 1} \left( {\mathcal E} _1 \right) \tilde{} \otimes _X {\mathcal O} _X \left( ( p + 1 ) D _{\tau} \right)$. Then, we have the following commutative diagrams
\begin{eqnarray}
\begin{matrix}
0 & \rightarrow & {\mathcal B} \left(
{\mathcal E} _1 \right) ^{p D _{\tau}}
 & \rightarrow & {\mathcal B} \left(
{\mathcal E} _1 \right) ^{( p + 1 ) D _{\tau}}
& \stackrel{\varphi ^p _i}{\rightarrow} & {\mathcal B} \left(
{\mathcal E} _1 \right) ^{( p + 1 ) D _{\tau}}
 \MID _{D _{\tau}} & \rightarrow & 0 &
\text{(exact)}\\
& & \cup s ^p & & \cup & & \cup u & \\
0 & \rightarrow & {\mathcal E} _i \cap _D {\mathcal B} \left(
{\mathcal E} _1 \right) ^{p D _{\tau}}
 & \rightarrow & {\mathcal E} _i \cap _D {\mathcal B} \left(
{\mathcal E} _1 \right) ^{( p + 1 ) D _{\tau}}
 & \rightarrow & {\mathcal Q} _{p + 1} \left( {\mathcal E} _1
\right) & \rightarrow & 0 &
\text{(exact)}
\end{matrix}\label{miniext}
\end{eqnarray}
for $i = 1, 2$.
\begin{claim}
In the diagram $(\ref{miniext})$, the extension ${\mathcal E} _i \cap _D {\mathcal B} \left(
{\mathcal E} _1 \right) ^{( p + 1 ) D _{\tau}}
$ of ${\mathcal Q} _{p + 1} \left(
{\mathcal E} _1
\right)$ by ${\mathcal E} _i \cap _D {\mathcal B} \left(
{\mathcal E} _1 \right) ^{p D _{\tau}}$ $($as $\tilde{G} \times G$-equivariant coherent ${\mathcal O} _X$-modules$)$ is uniquely determined by the inclusions $s ^p$ and $u$.
\end{claim}

\begin{proof}[Proof of Claim]
The preimage of ${\mathcal Q} _{p + 1} \left( {\mathcal E} _1 \right)$ under $\varphi ^p _i$ is uniquely determined in ${\mathcal B} \left( {\mathcal E} _1 \right) ^{( p + 1 ) D _{\tau}}$. By standard homological algebra (cf. Rotman \cite{Ro} Theorem 7.19 and Ex. 7.20-7.22), the ambiguity of extension is given by a morphism
$$h : {\mathcal Q} _{p + 1} \left( {\mathcal E} _1 \right) \rightarrow {\mathcal B} \left( {\mathcal E} _1 \right) ^{p D _{\tau}} / ( {\mathcal E} _i \cap _D {\mathcal B} \left( {\mathcal E} _1 \right) ^{p D _{\tau}} )$$
as $\tilde{G} \times G$-equivariant coherent ${\mathcal O} _{X ( \tau )}$-modules. ${\mathcal B} \left(
{\mathcal E} _1 \right) ^{p D _{\tau}} / ( {\mathcal E} _i \cap _D {\mathcal B} \left(
{\mathcal E} _1 \right) ^{p D _{\tau}} )$ is a successive extension of elements in $\{ {\mathcal B} \left( {\mathcal E} _1 \right) ^{q D _{\tau}} / ( \varphi ^{q - 1} _i ) ^{- 1} \left( {\mathcal Q} _q \left( {\mathcal E} _1 \right) \right) ; q \le p\}$ as $\tilde{G} \times G$-equivariant coherent ${\mathcal O} _{X ( \tau )}$-modules by construction. For each $q \le p$, ${\mathcal B} \left( {\mathcal E} _1 \right) ^{q D _{\tau}} / ( \varphi ^{q - 1} _i ) ^{- 1} \left( {\mathcal Q} _q \left( {\mathcal E} _1 \right) \right)$ is an ${\mathcal O} _{D _{\tau}}$-module. Moreover, $\mathbb{G} ^{\tau} _m$ acts on its fiber at $x _{\tau}$ by weight $- q$ (cf. Theorem \ref{embed ed} 3)). Similarly, $\mathbb{G} ^{\tau} _m$ acts on ${\mathcal Q} \left( {\mathcal E} _1, p + 1 \right) \otimes _X k ( x _{\tau} )$ by weight $- ( p + 1 )$. Therefore, we have
$$\Hom _{( \tilde{G} \times G, {\mathcal O} _{\mathbf O _{\tau}} )} \left( {\mathcal Q} _{p + 1} \left( {\mathcal E} _1 \right), {\mathcal B} \left( {\mathcal E} _1 \right) ^{q D _{\tau}} / ( \varphi ^{q - 1} _i ) ^{- 1} \left( {\mathcal Q} _q \left( {\mathcal E} _1 \right) \right) \right) = 0.$$
Every $\mathcal O  _{X ( \tau )}$-morphism from ${\mathcal Q} _{p + 1} \left( {\mathcal E} _1 \right)$ descends to a $\mathcal O  _{\mathbf O _{\tau}}$-morphism. Hence, we have $h = 0$ by induction on $q$.
\end{proof}

By the above claim, we have ${\mathcal E} _1 \cap _D {\mathcal B} \left(
{\mathcal E} _1 \right) ^{( p + 1 ) D _{\tau}}
 = {\mathcal
E} _2 \cap _D {\mathcal B} \left(
{\mathcal E} _1 \right) ^{( p + 1 ) D _{\tau}}$. Thus, the induction on $p$ proceeds.
\end{proof}

\begin{proof}[Proof of Proposition \ref{last}]
Let ${\mathcal E} _1, {\mathcal E} _2 \in {\bf Ob} EV \left( \Sigma \right) _{c}$
such that $\Xi \left( {\mathcal E} _1 \right) \cong \Xi \left(
{\mathcal E} _2 \right)$. By Corollary \ref{isotypical decomposition of category of vector bundles}
and Corollary \ref{reduction to trivial character}, we need only to prove the case
${\mathcal E} _1, {\mathcal E} _2 \in {\bf Ob} EV \left( \Sigma, 1
\right) _{c}$. 
By Lemma \ref{plast}, we have ${\mathcal E} _1 \MID _{X \left( \tau \right)} = {\mathcal
E} _2 \MID _{X \left( \tau \right)}$ as a $\tilde{G} \times \tilde{G}$-equivariant coherent
subsheaf of ${\mathcal B} \left( {\mathcal E}_1 \right) ^D \MID _{X \left( \tau \right)}$ for a sufficiently large equivariant divisor $D$. 
Collecting these isomorphisms (as coherent ${\mathcal O} _{X
\left( \tau \right)}$-submodules of ${\mathcal B} \left( {\mathcal E}_1 \right) ^D
 \MID _{X \left( \tau \right)}$) for all $\tau \in \Sigma ( 1 )$, we have ${\mathcal E} _1
= {\mathcal E} _2$ on $\cup
_{\tau \in \Sigma ( 1 )} X \left( \tau \right)$. Thus, we have
${\mathcal E} _1 = {\mathcal E} _2$ as vector bundles because a vector bundle on a normal variety is uniquely determined by its restriction to a complement of a codimension two locus. Hence, we have ${\mathcal
E} _1 = {\mathcal E} _2$ as $\tilde{G} \times \tilde{G}$-equivariant vector subbundles of ${\mathcal B} \left( {\mathcal E}_1 \right) ^D
$ by Corollary \ref{uniqueness}.
\end{proof}

\subsection{Proof of Transversality}\label{transversality}
This subsection is devoted to prove the following proposition.

\begin{proposition}\label{induced action}
For each ${\mathcal E} \in {\bf Ob} EV \left( \Sigma \right) _{c}$, we have $\Xi \left( {\mathcal E} \right) \in {\bf Ob} {\mathfrak C} \left( \Sigma \right) _{c}$.
\end{proposition}

Since the transversality condition (Definition \ref{adm}) is closed under direct
sums, we can reduce the problem to ${\bf Ob} EV \left( \Sigma, 1 \right) _{c}$ by Corollary \ref{isotypical decomposition of category of vector bundles} and Corollary \ref{reduction to trivial character}.

We fix an arbitrary $\tau \in \Sigma ( 1 )$.

\subsubsection{Left actions}\label{left}
\begin{definition}
For each $\tilde{G}$-module $V$, we introduce a left $\mathfrak u ^{\tau} _+$-action $L _{\tau}$ on $V ^{\tau} _{[t]}$ as follows:
$$L _{\tau} ( X ) : V ^{\tau} _{[t]} \ni \sum _{n \in \Z} v _n \otimes t ^n \mapsto \sum _{n \in \Z} ( X v _n) \otimes t ^n \in V ^{\tau} _{[t]}$$
for every $X \in \mathfrak u ^{\tau} _+$. We extend it to a left $\mathfrak g$-action $\tilde{L} _{\tau}$ on $V ^{\tau} _{[t]}$ as follows:
$$\tilde{L} _{\tau} ( X ) : V ^{\tau} _{[t]} \ni \sum _{n \in \Z} v _n \otimes t ^n \mapsto \sum _{n \in \Z} ( X v _n) \otimes t ^n \in V ^{\tau} _{[t]}$$
for every $X \in \mathfrak g$.
\end{definition}

\begin{lemma}\label{left action preserving}
For every ${\mathcal E} \in {\bf Ob} EV \left( \Sigma , 1 \right) _{c}$, the $\mathbb G ^{\tau} _m$-isotypical decomposition
$$\rho ^{\tau} _{\infty} ( {\mathcal E} ) = \left( \left( 
\bigoplus _{n \in \Z} F ^{\tau} \left( n, {\mathcal E} \right) \otimes k t ^{n} \right) \subset B \left( {\mathcal E} \right) _{[t]} ^{\tau} \right)$$
identifies the differential of the unipotent radical of the $P ^{\tau}$-action on $F ^{\tau} (\bullet)$ (coming from Proposition \ref{module}) with the induced action of the left $\mathfrak u ^{\tau} _+$-action. In particular, $\rho ^{\tau} _{\infty} ( {\mathcal E} )$ is preserved by the left $\mathfrak u ^{\tau} _+$-action.
\end{lemma}

\begin{proof}
Since the construction of $F ^{\tau} (\bullet)$ factors through $\mathbb A _{\tau}$, this is a rephrasement of Lemma \ref{module} and Corollary \ref{def of rho}.
\end{proof}

\subsubsection{Right action}\label{right}
Here we introduce a $\mathfrak{u} ^{\tau} _-$-action on $B ( \bullet ) ^{\tau} _{[t]}$ which is compatible with the natural $\mathfrak{u} ^{\tau} _-$-action on the fiber of ${\mathcal O} _{X} \otimes B ( \bullet )$ at $x _{\tau}$. What we do here is only to switch from left to right in the definition of $L _{\tau}$. However, this is little complicated because our description heavily relies on Definition \ref{main functor for trivial factor}, where we choose $B ( \bullet ) \otimes {\mathcal O} _X$ as the standard $\tilde{G} \times G$-equivariant vector bundle (its alternative was ${\mathcal O} _X \otimes B ( \bullet )$).

We define ${}^{\tau} \mathbb{G} _m$ to be the image of  $\left( \tau
\times 1 \right) : \mathbb{G} _m \rightarrow T \times T$. Let $V$ be a $\tilde{G}$-module. Then,
${}^{\tau} \tilde{\mathbb{G}} _m$ acts on a $\tilde{G} \times G$-equivariant vector bundle
$V \otimes {\mathcal O} _X$. We have ${}^{\tau} \tilde{\mathbb{G}} _m
\subset \triangle ^d ( \tilde{L} ^{\tau} ) \tilde{\mathbb{G}} ^{\tau}
_m$. In particular, ${}^{\tau} \tilde{\mathbb{G}}
_m$ acts on $V ^{\tau} _{[t]}$. Fix a set of representative $\{ \lambda _1, \lambda _2,
\ldots, \lambda _{h}\}$ of $\tilde{Z} ( G ) ^{\vee}$
in $X ^* ( \tilde{T} )$. Then, we define ${\mathcal R} \left( {\mathcal E} \right)$ for each ${\mathcal E} \in {\bf Ob} EV \left( \Sigma, 1 \right) _{c}$ as follows: 
$${\mathcal R} \left( {\mathcal E} \right) := \bigoplus _{i = 1} ^{h} \bigoplus _{\lambda - \lambda _i \in X ^* ( T )} {\mathcal L} _{- \lambda _i} \otimes V _{\lambda} ^{\oplus n _{\lambda}}.$$
Here $\bigoplus _{\lambda \in X ^* ( \tilde{T} )} V _{\lambda}
^{\oplus n _{\lambda}}$ is the $\tilde{G}$-isotypical decomposition of
$B ( \mathcal E )$. For an equivariant divisor $D$, let ${\mathcal R} \left( {\mathcal E} \right) ^D$ denote $\mathcal R ( \mathcal E ) \otimes _X \mathcal O _X ( D )$. By Lemma \ref{switch}, we have $B ( \mathcal E ) \otimes
{\mathcal O} _G \cong {\mathcal R} \left( {\mathcal E} \right)
\MID _{G}$ as $\tilde{G} \times \tilde{Z} ( G )$-modules. Then, by using Corollary \ref{sandwich} repeatedly, we have
$$\cdots \subset \mathcal B ( \mathcal E ) ^{- D} \subset {\mathcal R} \left( {\mathcal E} \right) \subset \mathcal B ( \mathcal E ) ^D \subset {\mathcal R} \left( {\mathcal E} \right) ^{2 D} \subset \cdots$$
for a sufficiently large equivariant divisor $H$. As a result, we have the following isomorphism of $( \triangle ^d ( \tilde{L} ^{\tau} ) \mathbb{G} ^{\tau} _m, k [t, t ^{-1}] )$-modules.
$$B ( \mathcal E ) ^{\tau} _{[t]} = \mathop{\varinjlim} _N \Gamma \left( \mathbb{A}
_{\tau}, \mathcal B ( \mathcal E ) ^{N D _{\tau}} \MID _{\mathbb{A} _{\tau}} \right) \cong \mathop{\varinjlim} _N \Gamma ( \mathbb{A} _{\tau}, {\mathcal R} ( {\mathcal E} ) ^{N D _{\tau}} \MID _{\mathbb{A} _{\tau}} )$$

\begin{lemma}
For every $\mathcal E \in {\bf Ob} EV ( \Sigma, 1 ) _c$, there exist a section $s : {\mathcal R} ( \mathcal E ) \MID _{x _{\tau}} \hookrightarrow B \left( {\mathcal E} \right) _{[t]} ^{\tau}$ which yields the following ${} ^{\tau} \tilde{\mathbb G} _m$-isotypical decomposition:
$$B \left( {\mathcal E} \right) _{[t]} ^{\tau} = \bigoplus _{n \in \mathbb Z} t ^n. s ( {\mathcal R} ( \mathcal E ) \MID _{x _{\tau}} ).$$
\end{lemma}

\begin{proof}
We have the residual map $\psi ^R _{\tau} ( \lambda ) : \left( {\mathcal L} _{- w _0 \lambda} \otimes V _{\lambda} \right) \MID _{\mathbb{A} _{\tau}} \rightarrow \left( {\mathcal L} _{- w _0 \lambda} \otimes V _{\lambda} \right) \MID _{x _{\tau} }$.
By Proposition \ref{maximal}, we have an inclusion (of ${}^{\tau}
\tilde{\mathbb{G}} _m$-modules)
$$\left( {\mathcal L} _{- w _0 \lambda} \otimes V _{\lambda} \right) \MID _{ x _{\tau} } \stackrel{s _{\lambda}}{\hookrightarrow} \bigoplus _{n = 0} ^{\left<  \tau, \lambda - w _0 \lambda \right>} V ^{\tau} _{\lambda} ( n ) \otimes k t ^{- n} \subset \bigoplus _{n = 0} ^{\left< \tau, \lambda - w _0 \lambda \right>} \bigoplus _{m \ge - n} V _{\lambda} ^{\tau} ( n ) \otimes k t ^{m} \subset ( V _{\lambda} ) ^{\tau} _{[t]}$$
 induced from a ${} ^{\tau} \tilde{\mathbb{G}} _m$-module
section of $\psi ^R _{\tau} ( \lambda )$. Here the image of $s _{\lambda}$ is contained in a ${}^{\tau} \tilde{\mathbb{G}} _m$-isotypical component since $\left( {\mathcal L} _{- \lambda _i} \otimes V _{\lambda} \right) \MID _{ x _{\tau} }$ is of the form $\lambda _i ^{\oplus \dim V _{\lambda}}$ regarded as a representation of ${}^{\tau} \tilde{\mathbb{G}} _m$. Moreover, the image of $s _{\lambda}$ generates $( V _{\lambda} ) ^{\tau} _{[t]}$ as a $k [ t, t ^{- 1} ]$-module.

Let $\phi ^R _{\tau} : {\mathcal R} ( {\mathcal E} ) \rightarrow {\mathcal R} ( {\mathcal E} ) \MID _{ x _{\tau} }$ be the residual map. It factors through $\psi ^R _{\tau} : {\mathcal R} ( {\mathcal E} ) \MID _{\mathbb{A} _{\tau}} \rightarrow {\mathcal R} ( {\mathcal E} ) \MID _{ x _{\tau} }$. We define a ${} ^{\tau} \tilde{\mathbb{G}} _m$-module inclusion ${\mathcal R} \left( \mathcal E \right) \MID _{ x _{\tau} } \stackrel{s}{\hookrightarrow} V ^{\tau} _{[t]}$ as
$$s := \sum _{i = 1} ^{h} \bigoplus _{\lambda - \lambda _i \in X ^* ( T )} t ^{\left< \tau, w _0 \lambda - \lambda _ i \right>} . s _{\lambda}  = \sum _{i = 1} ^{h} s _i.$$
Since we have $\left< \tau, w _0 \lambda - \lambda _ i \right> \in \Z$ for every $\tau \in \Sigma ( 1 )$, the above construction is well-defined. $s$ is a sum of the inclusions of the form $s _{\lambda}$ with $t$-twist. Hence, the image of $s$ generates $V ^{\tau} _{[t]}$ as a $k [ t, t ^{- 1} ]$-module. By the above construction, the image of $s _i$ is contained in a ${}^{\tau} \tilde{\mathbb{G}} _m$-isotypical component. Therefore, the image of $s _i$ form a ${}^{\tau} \tilde{\mathbb{G}} _m$-isotypical component since $\tilde{Z} ( G ) \times 1$ acts on the various $\mathrm{Im} s _i$ by distinct characters.
\end{proof}

Since ${\mathcal R} ( {\mathcal E} ) \MID _{x _{\tau}}$
has a natural $\mathfrak{u} ^{\tau} _-$-action induced from $( 1
\times U ^{\tau} _-) \subset \tilde{G} ^{\tau}$, we can introduce a $\mathfrak{u} ^{\tau} _-$-action on the image of $s$.

\begin{definition}
We define a right $\mathfrak u ^{\tau} _-$-action $R _{\tau}$ on $B \left( {\mathcal E} \right) _{[t]} ^{\tau}$ by
$$R _{\tau} ( X ) : B \left( {\mathcal E} \right) _{[t]} ^{\tau} \ni \oplus _{n \in \mathbb Z} t ^n. v _n \mapsto \oplus _{n \in \mathbb Z} t ^n. ( X v _n ) \in B \left( {\mathcal E} \right) _{[t]} ^{\tau}$$
for every $X \in \mathfrak u ^{\tau} _-$.
\end{definition}

\begin{corollary}\label{right action preserving}
$\rho ^{\tau} _{\infty} ( {\mathcal E} )$ is stable under the right $\mathfrak{u} ^{\tau} _-$-action for every ${\mathcal E} \in {\bf Ob} EV \left( \Sigma , 1 \right) _{c}$.
\end{corollary}

\begin{proof}
By using ${\mathcal R} ( {\mathcal E} ) ^{- D} \subset {\mathcal E} \subset {\mathcal R} ( {\mathcal E} ) ^D$ instead of $\mathcal B ( \mathcal E ) ^{- D} \subset {\mathcal E} \subset \mathcal B ( \mathcal E ) ^D$ in Definition \ref{main functor for trivial factor}, we define
$${} ^R F ^{\tau} \left( n, {\mathcal E} \right) := \phi _{\tau} ^R \left( {\mathcal E} ( - n D _{\tau} ) \cap _D {\mathcal R} \left( {\mathcal E} \right) \right) \subset {\mathcal R} \left( {\mathcal E} \right) \MID _{ x _{\tau} }$$
for every $n \in \Z$. This is a $\tilde{P} ^{- \tau}$-submodule of $B ( \mathcal E )$ by a similar argument as in the proof of Lemma \ref{module} via the same isomorphism $k ( e ) \cong k ( x _{\tau} )$ as in Definition \ref{main functor for trivial factor}. Then, we have
$$\Gamma ( \mathbb{A} _{\tau}, {\mathcal E} \MID _{\mathbb{A} _{\tau}} ) = \bigoplus _{m \in \Z} F ^{\tau} ( n, {\mathcal E} ) \otimes k t ^{- n} = \bigoplus _{m \in \Z} t ^{- n} . {}^R F ^{\tau} ( n, {\mathcal E} ) .$$
as compatible $( k [ t ], {} ^{\tau} \tilde{\mathbb{G}} \triangle ^d ( \tilde{L} ^{\tau} ) )$-submodules of $B ( {\mathcal E} ) ^{\tau} _{[t]}$. Here the latter direct sum gives the ${} ^{\tau} \tilde{\mathbb{G}}$-isotypical decomposition. Hence, the result follows from the definition of  $\rho ^{\tau} _{\infty} ( {\mathcal E} )$.
\end{proof}

\subsubsection{Compatibility}
We work under the same assumptions as in \S \ref{left}. The main technical result of this subsection is the following.

\begin{theorem}\label{compatibility of action}
Let $\alpha$ be a root of $\mathfrak{g}$ such that $e _{\alpha} \in \mathfrak{u} ^{\tau} _-$ $($i.e. $\left< \tau, \alpha \right> < 0)$. Then, the right $\mathfrak{u} ^{\tau} _-$-action $R _{\tau} \left( e _{\alpha} \right)$ of $e _{\alpha}$ is written as
$$R _{\tau} \left( e _{\alpha} \right) = c \tilde{L} _{\tau} \left( e _{\alpha} \right) t ^{- \left< \tau, \alpha \right>},$$
where $\tilde{L} _{\tau}$ is the left $\mathfrak g$-action and $c$ is a non-zero constant.
\end{theorem}

We first prove our main result (Proposition \ref{induced action}) by assuming Theorem \ref{compatibility of action}.

\begin{proof}[Final step of the proof of Proposition \ref{induced action}]
By Corollary \ref{def of rho}, we have
$$\rho ^{\tau} _{\infty} \left( {\mathcal E} \right) \cong \left(
\left( \bigoplus _{n \in \Z} F ^{\tau} \left( n, {\mathcal E} \right) \otimes k t ^{- n} \right) \subset B \left( {\mathcal E} \right) ^{\tau} _{[t]} \right).$$
Let $\alpha$ be a root such that $e _{\alpha} \in \mathfrak{u} ^{\tau} _-$. By Corollary \ref{right action preserving}, $R _{\tau} \left( e _{\alpha} \right)$ must preserve $\rho ^{\tau} _{\infty} \left( {\mathcal E} \right)$ in $B \left( {\mathcal E} \right) ^{\tau} _{[t]}$. By Theorem \ref{compatibility of action}, it suffices to check whether the composition of the $\tilde{L} _{\tau} \left( e _{\alpha} \right)$-action and multiplying $t ^{- \left< \tau, \alpha \right>}$ preserves $\rho ^{\tau} _{\infty} \left( {\mathcal E} \right)$ or not. Then, to preserve $\rho ^{\tau} _{\infty} \left( {\mathcal E} \right)$, we must have
$$(1 \otimes t ^{- \left< \tau, \alpha \right>}) .e _{\alpha} F
^{\tau} \left( n, {\mathcal E} \right) \otimes k t ^{- n} \subset F
^{\tau} \left( n + \left< \tau, \alpha \right>, {\mathcal E} \right) \otimes k t ^{- n - \left< \tau, \alpha \right>}.$$
Hence, we have
$$U \left( \mathfrak{u} ^{\tau} _- \right) ^{\tau} _m  F ^{\tau}
\left( n, {\mathcal E} \right) \subset F ^{\tau} \left( n + m, {\mathcal E} \right)$$
This is the $\tau$-transversality condition. Hence, we have proved the result for $\tau$. Since $\tau$ is an arbitrary one-dimensional cone, we obtain the result.
\end{proof}

Before the proof of Theorem \ref{compatibility of action}, we need some preparations. We define two full-subcategories $\mathfrak{L}$ and $\mathfrak{R}$ of $EV \left( \Sigma, 1 \right) _{c}$ as follows:

\begin{itemize}
\item ${\bf Ob} \mathfrak{L} := \{ V \otimes {\mathcal O} _X ; V \in \mathrm{Rep} \tilde{G} \}$;
\item ${\bf Ob} \mathfrak{R} := \{ \oplus _{\lambda \in X ^* ( \tilde{T} )} {\mathcal L} _{- w _0 \lambda} \otimes V _{\lambda} ^{\oplus n _{\lambda}} ; n _{\lambda} \in \Z _{\ge 0} \}$.
\end{itemize}

We put ${\mathcal V} ^R _{\lambda} := {\mathcal L} _{- w _0 \lambda} \otimes V _{\lambda}$ for each $\lambda \in X ^* ( \tilde{T} )$. $\mathfrak{L}$ is a rigid tensor category (cf. Deligne-Milne \cite{DMOS}) by
the usual tensor product over ${\mathcal O} _X$.
\begin{lemma} There exists a tensor structure $\otimes^{\prime}$ on $\mathfrak{R}$ which makes $\mathfrak{R}$ a rigid tensor category which is equivalent to $\mathrm{Rep} \tilde{G}$.
\end{lemma}
\begin{proof}
Since $\mathfrak{R}$ is semisimple as an
abelian category, it suffices to construct tensor products for each
pair of irreducible objects. For each dominant weights $\lambda,
\gamma$, we write an irreducible decomposition of $V _{\lambda}
\otimes V _{\gamma}$ by
$$V _{\lambda} \otimes V _{\gamma} = \bigoplus_{\epsilon \text{ : dominant}} V _{\epsilon} ^{\oplus m _{\lambda, \gamma} ^{\epsilon}}.$$
Again by semisimplicity, we can twist the above direct sum decomposition isotypical componentwise. We have
$${\mathcal V} ^R _{\lambda} \otimes _X {\mathcal V} ^R _{\gamma} = {\mathcal L} _{w _0 (\lambda + \gamma )} ^{\vee} \otimes \bigoplus_{\epsilon \text{ : dominant}} V _{\epsilon} ^{\oplus m _{\lambda, \gamma} ^{\epsilon}}.$$
We twist the component ${\mathcal L} _{w _0 (\lambda + \gamma)}
^{\vee} \otimes V _{\epsilon} ^{\oplus m _{\lambda, \gamma}
^{\epsilon}}$ by ${\mathcal L} _{w _0 (\epsilon - \lambda - \gamma)}
^{\vee}$. Notice that we have $\left< \tau, w _0 ( \epsilon - \lambda -
\gamma ) \right> \in \Z$ for every $\tau \in \Sigma ( 1 )$. In
particular, $D _{- w _0 (\epsilon - \lambda - \gamma)}$ is an
equivariant divisor. Hence, defining the twisted tensor product
$\otimes ^{\prime}$ by
$${\mathcal V} ^R _{\lambda} \otimes ^{\prime} {\mathcal V} ^R _{\gamma} := \bigoplus_{\epsilon \text{ : dominant}} \left( {\mathcal V} ^R _{\epsilon} \right) ^{\oplus m _{\lambda, \gamma} ^{\epsilon}}$$
yields the result.
\end{proof}
Hence, we obtain two different rigid tensor categories $\mathfrak{L}$ and $\mathfrak{R}$ which are equivalent to $\mathrm{Rep} \tilde{G}$. We denote their images under $\rho ^{\tau} _{\infty}$ by $\hat{\mathfrak{L}} ^{\tau}$ and $\hat{\mathfrak{R}} ^{\tau}$, respectively. Then, we have the following:
$$\rho ^{\tau} _{\infty} \left( {\mathcal E} \otimes _X {\mathcal F} \right) \cong \left( \left( \Gamma \left( \mathbb{A} _{\tau}, {\mathcal E} \MID _{\mathbb{A} _{\tau}} \right) \otimes _{k [ t ]} \Gamma \left( \mathbb{A} _{\tau}, {\mathcal F} \right) \MID _{\mathbb{A} _{\tau}} \right) \subset \left( B \left( {\mathcal E} \right) \otimes _k B \left( {\mathcal F} \right) \right) ^{\tau} _{[t]} \right)$$
for every ${\mathcal E}, {\mathcal F} \in {\bf Ob} EV \left( \Sigma, 1 \right) _{c}$, and
$$\rho ^{\tau} _{\infty} \left( {\mathcal V} ^R _{\lambda} \otimes ^{\prime} {\mathcal V} ^R _{\gamma} \right) \cong \bigoplus_{\epsilon \text{ : dominant}} t ^{\left<\tau, w _0 ( \epsilon - \lambda - \gamma ) \right>} \rho ^{\tau} _{\infty} \left( {\mathcal L} _{w _0 (\lambda + \gamma )} ^{\vee} \otimes V _{\epsilon} ^{\oplus m _{\lambda, \gamma} ^{\epsilon}} \right).$$

By the above formulas, we equip $\hat{\mathfrak{L}} ^{\tau}$ and $\hat{\mathfrak{R}} ^{\tau}$ with tensor products $\otimes _0$ and $\otimes ^{\prime} _0$ that arise from $\otimes _X$ and $\otimes ^{\prime}$, respectively. We also define a tensor product $\otimes _0$ of ${\mathfrak B} \left( \Sigma \right) ^{\tau}$ by
$$\left( A _1 \subset ( V _1 ) ^{\tau} _{[t]} \right) \otimes _0 \left( A _2 \subset ( V _2 ) ^{\tau} _{[t]} \right) := \left( \left( A _1 \otimes _{k [t]} A _2 \right) \subset ( V _1 \otimes V _2 ) ^{\tau} _{[t]} \right)$$
for each $( A _1 \subset ( V _1 ) ^{\tau} ), ( A _2 \subset ( V _2 ) ^{\tau} _{[t]} ) \in {\bf Ob} {\mathfrak B} \left( \Sigma \right) ^{\tau}$. This tensor product coincides with the tensor product of $EV \left( \Sigma, 1 \right) _{c}$ via $\rho ^{\tau} _{\infty}$. Therefore, it is an extension of $\otimes _0$ in $\hat{\mathfrak{L}}$.

\begin{proof}[Proof of Theorem \ref{compatibility of action}]
Exponential $\exp R _{\tau} \left( e _{\alpha} \right)$ of $R _{\tau}
\left( e _{\alpha} \right)$ acts on any object of $\hat{\mathfrak{R}}
^{\tau}$ and commutes with the tensor product $\otimes _0$ of two
elements of $\hat{\mathfrak{R}}$ (in ${\mathfrak B} \left( \Sigma
\right) ^{\tau}$). By construction, the right $\mathfrak{u} _-$-action
commutes with $t$-twists. Hence, $\exp R _{\tau} \left( e _{\alpha}
\right)$-action commutes with the tensor product $\otimes ^{\prime}
_0$ of $\hat{\mathfrak{L}}$. We define a twisted right $\mathfrak{u} ^{\tau} _-$-action on formal loop spaces by the following:
$$R \left( e _{\alpha} \right) t ^{\left< \tau, \alpha \right>} : V ^{\tau} _{[t]} \ni v \mapsto R \left( e _{\alpha} \right) v \mapsto \left( R \left( e _{\alpha} \right) t ^{\left< \tau, \alpha \right>} \right) v \in V ^{\tau} _{[t]}$$
for each $V \in \Rep \tilde{G}$. Here, the weight of this endomorphism of $V ^{\tau} _{[t]}$
with respect to $\mathbb{G} ^{\tau} _m$-action is zero. Thus, $R
\left( e _{\alpha} \right) t ^{\left< \tau, \alpha \right>}$
preserves $\rho ^{\tau} _{\infty} \left( V \otimes {\mathcal O} _X
\right)$. Hence, it operates on objects of $\hat{\mathfrak{L}}$ compatibly
with $\otimes _0$. Moreover, its image under the residual map
$\psi _{\tau} : V \otimes {\mathcal O} _{\mathbb{A} _{\tau}}
\rightarrow V \otimes k ( x _{\tau} )$ is non-zero if $R _{\tau} \left( e
_{\alpha} \right)$ operates $V ^{\tau} _{[t]}$ nontrivially. $\psi _{\tau}$ commutes with the decomposition rule of the tensor category $\Rep \tilde{G}$. Hence, $\psi _{\tau}$ is a fiber functor of $\hat{\mathfrak{L}}$. Here $\exp \left( s R \left( e _{\alpha} \right) t ^{\left< \tau, \alpha \right>} \right)$ commutes with $\psi _{\tau}$ for every $s \in k$. Thus, $\exp \left( s R \left( e _{\alpha} \right) t ^{\left< \tau, \alpha \right>} \right)$ defines some element of $\tilde{G}$ (as an automorphism of $V \otimes k ( x _{\tau} )$ for every $V \in \Rep \tilde{G}$) by \cite{DMOS} II Prop. 2.8. Hence, its derivative $R \left( e _{\alpha} \right) t ^{\left< \tau, \alpha \right>}$ defines some element of $\mathfrak{g}$. The weight of $R \left( e _{\alpha} \right) t ^{\left< \tau, \alpha \right>}$ is $\alpha$ by $\triangle ^d ( \tilde{T} )$-action. Moreover, every root space of a reductive Lie algebra is one-dimensional. Therefore, $R \left( e _{\alpha} \right) t ^{\left< \tau, \alpha \right>}$ must operate as a nonzero constant multiple of $\tilde{L} _{\tau} \left( e _{\alpha} \right)$.
\end{proof}

\begin{remark}
The term "tensor category" in \cite{DMOS} (and this paper) is the same as "symmetric tensor category" in the standard notation after the invention of quantum groups (cf. Bakalov and Kirillov \cite{BK}).
\end{remark}
\section{Main theorem}\label{MT}
Now we can state our main result as follows:

\begin{theorem}\label{Main Theorem}
We have an equivalence of categories
$$\Xi : EV \left( \Sigma \right) _{c}
\stackrel{\cong}{\rightarrow} {\mathfrak C} \left( \Sigma \right)
_{c}.$$
Moreover, $\Xi$ induces a category equivalence $EV \left( \Sigma \right) \cong {\mathfrak C} \left( \Sigma \right)$.
\end{theorem}

\begin{example}[Tangent bundles]
Let $G$ be a semisimple adjoint group. Let $X = X ( \Sigma _0 )$ be its wonderful compactification. Then, the set of one-dimensional cones of $\Sigma _0$ is represented by the set of fundamental co-weights $\omega _1 ^{\vee}, \ldots, \omega _r ^{\vee}$ (cf. \S \ref{NoteVar}). In this setting, we have
$$\Xi ( T X ) \cong ( \mathfrak g, \{ F ^{\omega _i ^{\vee}} \left( \bullet \right) \} ),$$
where
$$F ^{\omega _i ^{\vee}} \left( n \right) = \left\{ \begin{array}{ll}
\bigoplus _{\left< \omega _i ^{\vee}, \alpha \right> \ge n} k e_{\alpha} & (n \ge 2)\\
k H _i \oplus \bigoplus _{\left< \omega _i ^{\vee}, \alpha \right> \ge n} k e _{\alpha} & (n = 1)\\
\mathfrak g & (n \le 0)
\end{array} \right. .$$
Here $e _{\alpha}$ is a root vector of $\mathfrak g$ and we put $H _i := [ e _{\alpha _i}, e _{- \alpha _i} ] \in \mathfrak t$. If we replace $F ^{\omega _i ^{\vee}} \left( 1 \right)$ by $\oplus _{\left< \omega _i ^{\vee}, \alpha \right> \ge 1} k e _{\alpha}$, then the corresponding vector bundle is the ``logarithmic" tangent bundle of $X$ along the reduced union of boundary divisors. (cf. Bien-Brion \cite{BB}.)
\end{example}

\begin{example}[Case $G = PGL _2$]\label{example}
In this case, the only nontrivial (partial) compactification of
$\mathop{PGL} _2$ is the projecfication $\mathbb{P} \left( M _2 \right)$ of the set of two by two matrices $M _2$, which is 
the wonderful compactification. Here, we have
$\tilde{G} = \mathop{SL} _2$. Moreover, $\mathop{SL} _2 \times \mathop{SL}
_2$-action on $\mathbb{P} \left( M _2
\right)$ factors through the following action:
$$\mathop{SL _2} \times \mathop{SL _2} \times M _2 \ni ( g _1, g _2, A ) \mapsto g _1 A g _2 ^{-1} \in M _2.$$
We have a unique boundary divisor $D$ of $\mathbb{P} \left( M _2
\right)$ described as the zeros of the determinant. Let $\alpha$ be a unique positive root of $\mathop{PGL} _2$ and let $\varpi$ be the unique fundamental coweight of $\mathop{PGL} _2$.

For a $\mathop{SL} _2$-module $V$, we have the $\tilde{T}$-isotypical decomposition $V = \oplus _{\lambda \in X ^{*} ( \tilde{T} )} V ^{\lambda}$. In this case, we define the character $\mathrm{ch} V$ of $V$ as $\mathrm{ch} V := \sum _{\lambda \in X ^{*} ( \tilde{T} )} ( \dim V ^{\lambda} ) e ^{\lambda}$.

Let ${\mathcal E}$ be a $\mathop{SL} _2 \times \mathop{SL} _2$-equivariant vector bundle on $X$ such that ${\mathcal E} \otimes _X k ( e ) \cong \mathfrak{sl} _2 \boxtimes k$ as $\mathop{SL} _2 \times \Z / 2 \Z$-modules. By Theorem \ref{Main Theorem}, we can classify such $\mathop{SL} _2 \times \mathop{SL} _2$-equivariant vector bundles in ${\bf Ob} {\mathfrak C} \left( \Sigma _0 \right)$. We have $\Xi \left( {\mathcal E} \right) \cong \left( \mathfrak{sl} _2 \boxtimes k, F ^{\varpi} \left( \bullet \right) \right)$. Then, we have the following classification by the transversality condition of Definition \ref{adm}.

Every object in ${\mathfrak C} \left( \Sigma _0 \right)$ with the form $\left( \mathfrak{sl} _2 \boxtimes k, F ^{\varpi} \left( \bullet \right) \right)$ is given by one of the following four objects up to degree shift:
\begin{eqnarray*}
\mathop{ch} F ^{\varpi} \left( n, \mathfrak{sl} _2 \otimes {\mathcal O} _X \right) & = & \left\{ \begin{array}{ll}
0 & (n \ge 1)\\
e ^{\alpha} + e ^0 + e ^{- \alpha} & (n \le 0)
\end{array} \right.\\
\mathop{ch} F ^{\varpi} \left( n, {\mathcal O} _X \otimes \mathfrak{sl} _2 \right) & = & \left\{ \begin{array}{llll}
0 & (n \ge 2)\\
e ^{\alpha} & (n = 1)\\
e ^{\alpha} + e ^0 & (n = 0)\\
e ^{\alpha} + e ^0 + e ^{- \alpha} & (n \le - 1)
\end{array} \right.\\
\mathop{ch} F ^{\varpi} \left( n, T X \right) & = & \left\{ \begin{array}{lll}
0 & (n \ge 2)\\
e ^{\alpha} + e ^0 & (n = 1)\\
e ^{\alpha} + e ^0 + e ^{- \alpha} & (n \le 0)\\
\end{array} \right.\\
\mathop{ch} F ^{\varpi} \left( n, T ^* X \right) & = & \left\{ \begin{array}{lll}
0 & (n \ge 1)\\
e ^{\alpha} & (n = 0)\\
e ^{\alpha} + e ^0 + e ^{- \alpha} & (n \le - 1)\\
\end{array} \right.
\end{eqnarray*}

We say a $\mathop{SL} _2 \times \mathop{SL} _2$-equivariant vector bundle ${\mathcal E}$ on $\mathbb{P} \left( M _2
\right)$ irreducible, if there exists no proper $\mathop{SL} _2 \times \mathop{SL} _2$-equivariant vector subbundle in $EV \left( \Sigma _0 \right)$. Applying the inverse functor of $\Xi$, we obtain:

{\it Every irreducible rank three $\mathop{SL _2} \times \mathop{SL _2}$-equivariant vector bundle on $\mathbb{P} ^3$ is isomorphic to $\mathfrak{sl} _2 \otimes {\mathcal O} _X$, ${\mathcal O} _X \otimes \mathfrak{sl} _2$, $T X$, or $T ^* X$ up to line bundle twist.}

Similarly, we can prove that the number of isomorphism classes of irreducible $\mathop{SL} _2 \times \mathop{SL} _2$-equivariant vector bundles of rank $n$ on $\mathbb{P} ^3$ is just $2 ^{n - 1}$ up to line bundle twist.
\end{example}

\begin{remark}[Tensor structures]
All categories in Theorem \ref{Main Theorem} have natural tensor products. In fact, $\Xi$ is a tensor functor when we restrict to $EV \left( \Sigma, 1 \right) _c$ or $EV \left( \Sigma, 1 \right)$. However, due to our presentation of the category $\mathfrak C$, $\Xi$ does not preserve natural tensor products in general.
\end{remark}

We postpone the proof of our main theorem until \S \ref{finish} since we need the inverse functor $\Phi : {\mathfrak C} \left( \Sigma \right) _{c} \rightarrow EV \left( \Sigma \right) _{c}$ in order to prove Theorem \ref{Main Theorem}.

\subsection{Construction of the inverse functor}\label{reverse}
In this subsection, we construct an inverse functor $\Phi$ of $\Xi$. To achieve this, we construct some $\tilde{G} \times G$-equivariant vector bundle from a given object $\left( V , \{ F ^{\tau} \left( \bullet \right) \} \right)$ of ${\mathfrak C} \left( \Sigma, 1 \right) _{c}$.

\subsubsection{The functor $\Phi$}
Here we devote ourselves to construct a functor $\Phi$ by assuming the following Proposition \ref{simple makeup}. This Proposition is proved in \S \ref{sc}.

\begin{proposition}\label{simple makeup}
Let ${\mathcal A} := \left( V , \{ F ^{\tau} \left( \bullet \right) \}
\right) \in {\bf Ob} {\mathfrak C} \left( \Sigma, 1 \right) _{c}$ and choose
an arbitrary $\tau \in \Sigma ( 1 )$. Assume that the following condition holds:

\begin{itemize}
\item[$\diamond$] $F ^{\xi} \left( 0 \right) = V$ and $F ^{\xi} \left( 1 \right) = \{ 0 \}$ for every $\xi \in \Sigma ( 1 ) \backslash \{ \tau \}$. 
\end{itemize}

Then, there exist a $\tilde{G} \times G$-equivariant vector bundle $\Phi ^{\tau} ( {\mathcal A} )$ with the following properties:
\begin{enumerate}
\item[a.] $\Xi \left( \Phi ^{\tau} ( {\mathcal A} ) \right) = \left( V , \{ F ^{\tau} \left( \bullet \right) \} \right)$;
\item[b.] There exist $\tilde{G} \times G$-equivariant inclusions $V \otimes {\mathcal O} _X
\left( - N D _{\tau} \right)
\hookrightarrow \Phi ^{\tau} ( {\mathcal A} ) \hookrightarrow V
\otimes {\mathcal O} _X \left( N D _{\tau} \right)$ for a sufficiently
large integer $N$.
\end{enumerate}
Moreover, for each ${\mathcal B} := \left( V _1 , \{ F ^{\tau} _1 \left( \bullet \right) \} \right) \in {\bf Ob} {\mathfrak C} \left( \Sigma, 1
\right) _{c}$ which satisfies $\diamond$, we have a natural inclusion
$$\Hom _{{\mathfrak C} \left( \Sigma \right) _{c}} \left( {\mathcal
A}, {\mathcal B} \right) \hookrightarrow \Hom _{EV \left( \Sigma
\right) _{c}} \left( \Phi ^{\tau} \left( {\mathcal
A} \right), \Phi ^{\tau} \left( {\mathcal B} \right) \right).$$
\end{proposition}
For a $\tilde{G}$-module $V$, we define a filtration $F _{null}
\left( \bullet, V \right)$ of
$V$ indexed by $\Z$ as follows:

$$F _{null} \left( n, V \right) = 0 \text{ if } n \ge 1 \text{ and } V
\text{ if } n \le 0.$$
This is clearly an $\tau$-transversal filtration for every $\tau \in
\Sigma ( 1 )$. 

For every ${\mathcal A} = \left( V , \{ F ^{\tau} \left( \bullet \right) \}
\right) \in {\bf Ob} {\mathfrak C} \left( \Sigma, 1 \right) _{c}$ and every
$\tau \in \Sigma ( 1 )$, we define
$${\mathcal A} _{\tau} := \left( V, \{ F _{null} \left( \bullet, V \right)
\} _{\eta \in \Sigma ( 1 ) \backslash \{ \tau \}} \cup \{ F ^{\tau}
\left( \bullet \right)\} \right) \in {\bf Ob} {\mathfrak C} \left( \Sigma, 1 \right) _{c}.$$
${\mathcal A} _{\tau}$ satisfies the assumption of Proposition
\ref{simple makeup}. We put $D _0 := \sum _{\tau \in \Sigma ( 1 )} D _{\tau}$ and $D ^c _{\eta} :=  \sum _{\tau \in \Sigma ( 1 ) \backslash \eta} D _{\tau}$ for each $\eta \in \Sigma ( 1 )$. For each $\mathcal A \in {\bf Ob} \mathfrak C ( \Sigma, 1 ) _c$, we define
$$\Phi _1 \left( {\mathcal A} \right) := \bigcap _{\tau \in \Sigma ( 1 )}
\Phi ^{\tau} \left( {\mathcal A} _{\tau} \right) \otimes _X {\mathcal O}
_X ( M D _{\tau} ^c ) \subset V \otimes {\mathcal O} _X ( M D _0 )$$
Here we assume that $M$ is a sufficiently large integer and drop the sufficiently large equivariant divisor $D$ needed to define the intersection (Lemma \ref{intersection of equivariant vector bundles}) by letting $D = M D _0$. We use this convention throughout this subsection. It is clear that $\Phi _1 \left( {\mathcal A} \right)$ does not depend on the choice of $M$. We have $\Phi _1 \left(
{\mathcal A} \right) \MID _{X \left( \tau \right)} \cong \Phi ^{\tau} \left(
{\mathcal A} _{\tau} \right) \MID _{X
\left( \tau \right)}$ for each $\tau \in \Sigma ( 1 )$.

\begin{lemma}\label{vector bundleness of Phi}
For every $( V , \{ F ^{\tau} \left( \bullet \right) \}
) \in {\bf Ob} {\mathfrak C} \left( \Sigma, 1 \right) _{c}$, $\Phi _1
\left( ( V , \{ F ^{\tau} \left( \bullet \right) \}
) \right)$ is a $\tilde{G} \times G$-equivariant vector bundle.
\end{lemma}

\begin{proof}
What we need to show is that $\Phi _1
\left( ( V , \{ F ^{\tau} \left( \bullet \right) \}
) \right)$ is a vector bundle. We restrict our construction from $X \left( \Sigma \right)$ to $T \left( \Sigma \right)$ to obtain the corresponding object in Klyachko's category. We fix an inclusion $\Phi _1 \left( ( V , \{ F ^{\tau} \left( \bullet \right) \}
) \right) \otimes _X {\mathcal O} _X \left( M D _{\tau} ^c \right) \subset V
\otimes {\mathcal O} _X \left( M D
_0 \right)$. It induces the corresponding inclusion between the identity
fibers. For each $\sigma \in \Sigma$, the family of filtrations $F _{null} \left( \bullet, V \right)$ and $\{ F ^{\tau} \left(
\bullet \right) \} _{\tau \in \sigma ( 1 )}$ forms a distributive
lattice. Hence, they satisfies the condition of Klyachko's
category \cite{Kl} 1.1. As a result,
$$\bigcap _{\tau \in \Sigma ( 1 )}
\left( \Phi ^{\tau} \left( ( V , \{ F ^{\tau} \left( \bullet \right) \}
)  _{\tau} \right) \otimes _X {\mathcal O}
_X ( M D _{\tau} ^c )
\MID _{T \left( \Sigma  \right)} \right) \subset V
\otimes {\mathcal O} _X ( M D _0 ) \MID _{T \left( \Sigma \right)}$$
is a vector bundle for a sufficiently large integer $M$. Since $\Phi _1 \left( ( V , \{ F ^{\tau} \left( \bullet \right) \} ) \right)$ is a $\tilde{G} \times G$-equivariant coherent sheaf, the local structure theorem (Theorem \ref{structure theorem}) asserts that its restriction to the open subset
$$T ( \Sigma ) \times \mathbb A ^{2 \dim U} \cong U ^- .T ( \Sigma ). U\subset X ( \Sigma )$$
is still a vector bundle. Further, this open subset meets all $( \tilde{G} \times \tilde{G} )$-orbits, which yields the result.
\end{proof}

It follows that the map $\Phi _1 :
{\bf Ob} {\mathfrak C} \left( \Sigma, 1 \right) _{c} \rightarrow {\bf Ob} EV \left( \Sigma, 1
\right) _{c}$ satisfies $\Xi \circ \Phi _1 = \id$ on ${\bf Ob} {\mathfrak C} \left( \Sigma, 1 \right) _{c}$. 

\begin{lemma}\label{faithfulness of Phi}
For every ${\mathcal A}, {\mathcal B} \in {\bf Ob} {\mathfrak C}
\left( \Sigma, 1 \right) _{c}$, we have a natural inclusion
$$\Hom _{{\mathfrak C} \left( \Sigma \right) _{c}} \left( {\mathcal A}, {\mathcal B} \right) \hookrightarrow \Hom _{EV \left( \Sigma
\right) _{c}} \left( \Phi _1 \left( {\mathcal A} \right), \Phi _1 \left( {\mathcal B} \right) \right)$$
\end{lemma}

\begin{proof}
Put ${\mathcal A} := ( V , \{ F ^{\tau} \left( \bullet \right) \} )$ and ${\mathcal B} := ( V _1, \{ F ^{\tau} _1 \left( \bullet \right) \} )$.
For each morphism $f : ( V , \{ F ^{\tau} \left( \bullet \right) \} ) \rightarrow ( V _1 , \{
F _1 ^{\tau} \left( \bullet \right) \} )$ in ${\mathfrak C} \left(
\Sigma, 1 \right) _{c}$, we have a base morphism $V \rightarrow V
_1$. Hence, we have a corresponding morphism $f _{\tau} : ( V , \{ F ^{\tau}
\left( \bullet \right) \} ) _{\tau} \rightarrow ( V _1 , \{
F _1 ^{\tau} \left( \bullet \right) \} ) _{\tau}$ for every $\tau \in
\Sigma ( 1 )$. By Proposition \ref{simple makeup}, we have a
corresponding morphism $f ^{\prime} _{\tau} : \Phi ^{\tau} \left( ( V , \{ F ^{\tau}
\left( \bullet \right) \} ) _{\tau} \right) \rightarrow \Phi ^{\tau} \left( ( V _1 , \{
F _1 ^{\tau} \left( \bullet \right) \} ) _{\tau} \right)$. Then, by taking the intersection of the
both sides compatibly with $f \otimes \id : V \otimes {\mathcal O}
_X ( M D _0 ) \rightarrow V _1 \otimes {\mathcal O}
_X ( M D _0 )$, we have a
corresponding morphism
\begin{align*}
&f ^{\prime \prime} : \Phi _1 \left( ( V , \{ F ^{\tau} \left( \bullet
\right) \} ) \right) = \bigcap _{\tau \in \Sigma ( 1 )}
\left( \Phi ^{\tau} \left( ( V , \{ F ^{\tau} \left( \bullet
\right) \} ) _{\tau} \right) \otimes _X {\mathcal O}
_X ( M D _{\tau} ^c ) \right) \\
& \stackrel{\cap f ^{\prime} _{\tau} \otimes \mathrm{id}}{\rightarrow} \bigcap _{\tau \in \Sigma ( 1 )}
\left( \Phi ^{\tau} \left( ( V _1 , \{ F ^{\tau} _1 \left( \bullet \right) \} ) _{\tau} \right) \otimes _X {\mathcal O} _X ( M D _{\tau} ^c ) \right) = \Phi _1 \left( ( V _1 , \{ F _1 ^{\tau} \left(
\bullet \right) \} ) \right).
\end{align*}
If $f$ is nontrivial, then its base morphism $V \rightarrow V _1$ is also
nontrivial. Hence, $f ^{\prime \prime}$ induces nontrivial
morphism.
\end{proof}

For each $\lambda \in X ^* ( \tilde{T} )$, we have a category
equivalence $\lambda : {\mathfrak C} \left( \Sigma, 1
\right) _c \rightarrow {\mathfrak C} \left( \Sigma, \bar\lambda
\right) _c$ induced from
$${\bf Ob} {\mathfrak C} \left( \Sigma, 1
\right) _c \ni \left( V \boxtimes k, \{ F ^{\tau} \left( \bullet \right)
\}\right) \mapsto \left( V \boxtimes \bar{\lambda}, \{ F ^{\tau} \left( \bullet - \lfloor \left< \tau,
\lambda \right>\rfloor  \right)
\}\right) \in {\bf Ob} {\mathfrak C} \left( \Sigma, \bar\lambda
\right) _c.$$

For $\lambda \in X ^* \left( T \right)$, $\lambda$ yields a category auto-equivalence which corresponds to $\otimes _X {\mathcal O} _X \left( D ^{\lambda} \right)$ via $\Xi$ (from Lemma \ref{divisor} and Definition \ref{general redef}). Hence, we have the following commutative diagram.
$$
\begin{matrix}
{\mathfrak C} \left( \Sigma, 1
\right) _{c} & \stackrel{\lambda}{\rightarrow} & {\mathfrak C} \left( \Sigma,  1 \right) _c &\\
\Phi _1 \downarrow \quad& &\quad \downarrow \Phi _1 & \text{for every $\lambda \in  X ^* \left( T \right)$}\\
EV \left( \Sigma, 1 \right) _{c} & \stackrel{\otimes _X {\mathcal O}
_{X} \left( D ^{\lambda} \right)}{\rightarrow} & EV \left( \Sigma, 1 \right) _{c} &
\end{matrix}
$$
Similarly, we have a category auto-equivalence $\lambda$ on ${\mathfrak C} \left( \Sigma \right) _{c}$ for every $\lambda \in X ^* ( \tilde{T} )$ as the direct sum of the composition of the inverse of $\mu$ and $\lambda + \mu$ in
$${\mathfrak C} \left( \Sigma, \bar\mu \right) _{c} \stackrel{\mu}{\leftarrow} {\mathfrak C} \left( \Sigma, 1 \right) _{c} \stackrel{\lambda + \mu}{\rightarrow} {\mathfrak C} \left( \Sigma, \bar\lambda + \bar\mu \right) _{c} \text{ for }\mu \in X ^* ( \tilde{T} ).$$
Notice that the above construction is independent of the representative $\mu$ of $\bar\mu$. By the category auto-equivalence $\lambda$ on ${\mathfrak C} \left( \Sigma \right) _{c}$, define the category equivalence $\Phi _{\bar\lambda}$ as
$$\Phi _{\bar\lambda} := ( \otimes _X \mathcal L _{\lambda}) \circ \Phi _1 \circ ( \lambda )^{- 1} : {\mathfrak C} \left( \Sigma, \bar\lambda \right) _{c} \rightarrow EV \left( \Sigma, \bar\lambda \right) _{c}.$$
Now, we extend the category equivalence from ${\mathfrak C} \left( \Sigma, 1 \right) _{c} \rightarrow EV
\left( \Sigma, 1 \right) _{c}$ to ${\mathfrak C} \left( \Sigma \right) _{c} \rightarrow EV
\left( \Sigma \right) _{c}$ by setting $\Phi := \oplus _{\chi \in \tilde{Z} ( G ) ^{\vee}} \Phi _{\chi}$.

Summarizing the above, we have the following result.

\begin{proposition}\label{reverse functor Phi}
$\Phi : {\mathfrak C} \left( \Sigma \right) _{c} \rightarrow EV
\left( \Sigma \right) _{c}$ is a faithful covariant functor. Moreover, $\Xi
\circ \Phi$ is the identity on ${\bf Ob} {\mathfrak C} \left( \Sigma \right)_{c}$.
\end{proposition}

\subsubsection{Simple case}\label{sc}
Here we prove Proposition \ref{simple makeup}. We fix $\tau \in \Sigma ( 1 )$ which is the same as in the statement of Proposition \ref{simple makeup}. First, we recall some notation and prove a preliminary result.

By the description of \S \ref{NoteVar}, we have $\mathbf O _{\tau} \cong ( G \times G ) / G ^{\tau}$. Moreover, the base point $1 \times 1 \mod G ^{\tau}$ corresponds to $x _{\tau}$. We denote the Lie algebra of $G ^{\tau}$ by $\mathfrak{g} ^{\tau}$.

\begin{lemma}\label{uniq stable action}
For every ${\mathcal E} \in {\bf Ob} EV \left( \Sigma, 1 \right) _{c}$, combining $\triangle ^d ( \tilde{L} ^{\tau} ) \tilde{\mathbb G} ^{\tau} _m$-action, the left $\mathfrak{u} ^{\tau} _+$-action, and the right $\mathfrak{u} ^{\tau} _-$-action, we can introduce a $\mathfrak{g} ^{\tau}$-module structure on
$$\Gamma \left( 
\mathbb{A} _{\tau}, {\mathcal E} \MID _{\mathbb{A} _{\tau}} \right) / t \Gamma \left( 
\mathbb{A} _{\tau}, {\mathcal E} \MID _{\mathbb{A} _{\tau}} \right) ( \cong {\mathcal E} \otimes _X k ( x _{\tau} ) \text{ as ${\mathcal O} _{X, x _{\tau}}$-modules} )$$
which coincides with the natural $\mathfrak{g} ^{\tau}$-module structure on ${\mathcal E} \otimes _X k ( x _{\tau} )$. 
\end{lemma}

\begin{proof}
We have $\mathfrak{g} ^{\tau} = \mathrm{Lie} ( \triangle ^d ( L ^{\tau} ) \mathbb{G} ^{\tau} _m ) + (\mathfrak{u} ^{\tau} _+ \oplus \mathfrak{u} ^{\tau} _-) \subset \mathfrak{g} \oplus \mathfrak{g}$ (see \S \ref{NoteVar}).
By construction, $\Gamma \left( 
\mathbb{A} _{\tau}, {\mathcal E} \MID _{\mathbb{A} _{\tau}} \right)$ inherits a natural $\triangle ^d ( \tilde{L} ^{\tau} ) \tilde{\mathbb G} ^{\tau} _m$-action from ${\mathcal E}$. Hence, the only nontrivial part is about the left $\mathfrak{u} ^{\tau} _+$-action and the right $\mathfrak{u} ^{\tau} _-$-action. We show that the left $\mathfrak{u} ^{\tau} _+$-action and the right $\mathfrak{u} ^{\tau}
_-$-action defined in \S \ref{transversality} indeed satisfies the above
property. For each $\alpha, \beta \in \triangle$ such that $e
_{\alpha} \in \mathfrak{u} ^{\tau} _+$ and $e _{\beta} \in
\mathfrak{u} ^{\tau} _-$, we have $[L _{\tau} ( e _{\alpha} ), R _{\tau} ( e _{\beta} )] = [L _{\tau} ( e _{\alpha} ), \tilde{L} _{\tau} ( e _{\beta} )] t ^{- \left< \tau, \beta \right>}$. $\rho ^{\tau} _{\infty} \left( {\mathcal E} \right)$ is preserved by both the left $\mathfrak{u} ^{\tau} _+$-action and the right $\mathfrak{u} ^{\tau} _-$-action. Hence, for each $v \in \Gamma \left( \mathbb{A} _{\tau}, {\mathcal E} \MID _{\mathbb{A} _{\tau}} \right)$, we have
$$\left( \tilde{L} _{\tau} ( [e _{\alpha}, e _{\beta} ] ) \otimes t ^{- \left< \tau, \alpha + \beta \right>} \right) v \in \Gamma \left( \mathbb{A} _{\tau}, {\mathcal E} \MID _{\mathbb{A} _{\tau}} \right).$$
We have $\left< \tau, \alpha + \beta \right> > \left< \tau, \beta \right>$. As a result, the left $\mathfrak{u} ^{\tau} _+$-action and the right $\mathfrak{u} ^{\tau}
_-$-action commute on $\Gamma \left( 
\mathbb{A} _{\tau}, {\mathcal E} \MID _{\mathbb{A} _{\tau}} \right) / t \Gamma \left( 
\mathbb{A} _{\tau}, {\mathcal E} \MID _{\mathbb{A} _{\tau}} \right)$. Therefore, we can restrict ourselves to $\triangle ^d ( \tilde{L} ^{\tau} ) \tilde{\mathbb{G}} ^{\tau} _m ( U ^{\tau} _+ \times 1)$. In this case, the assertion is a consequence of Lemma \ref{module} and the isotypical decomposition of Corollary \ref{def of rho}.
\end{proof}

\begin{proof}[Proof of Proposition \ref{simple makeup}]
Let $N$ be an integer such that $F ^{\tau} \left( N - 1 \right) = \{ 0 \}$ and $F ^{\tau} \left( - N \right) = V$.
By induction on $n$, we construct a $\tilde{G} \times G$-equivariant vector bundle
$\Phi ^{\tau} _{n}\left( {\mathcal A} \right)$ with the following
properties (which we call $( \sharp ) _n$):

\begin{enumerate}
\item $V \otimes {\mathcal O} _X \left( - N D _{\tau} \right)
\hookrightarrow \Phi ^{\tau} _{n} \left( {\mathcal A} \right)
\hookrightarrow V \otimes {\mathcal O} _X \left( n D _{\tau}
\right)$ as $\tilde{G} \times G$-equivariant coherent sheaves;
\item $F ^{\tau} \left( m,\Phi ^{\tau} _{n} \left( {\mathcal A} \right) \right) = F ^{\tau} \left( m + N - n \right)$ for $m > - N$;
\item $F ^{\tau} \left( m, \Phi ^{\tau} _{n} \left( {\mathcal A} \right)
 \right) = V$ for $m \le - N$.
\end{enumerate}

Here we put $\Phi ^{\tau} \left( {\mathcal A} \right) := \Phi ^{\tau} _{N} \left( {\mathcal A} \right)$. Since $( \sharp ) _N$ is the same as the conditions of Proposition \ref{simple makeup}, $\Phi ^{\tau} \left( {\mathcal A} \right)$ has the desired property if such an induction proceeds. We have $F ^{\tau} \left( m, \Phi ^{\tau} _{- N} \left( {\mathcal A} \right) \right) = F ^{\tau} \left( m + 2 N \right) = 0$ for $m > - N$. Thus,
$\Phi ^{\tau} _{- N} \left( {\mathcal A} \right) := V \otimes
{\mathcal O} _X \left( - N D _{\tau} \right)$ satisfies the property
$( \sharp ) _{- N}$.

Then, assuming $( \sharp ) _n$, we construct $\Phi ^{\tau} _{n
+ 1} \left( {\mathcal A} \right)$ with the property $( \sharp ) _{n + 1}$.

We denote $\mathrm{coker} \left[ \Phi ^{\tau} _n \left( {\mathcal A} \right) \rightarrow \Phi ^{\tau} _n \left( {\mathcal A} \right)
\otimes _X {\mathcal O} _X \left( D _{\tau} \right) \right]$ by
${\mathcal Q} _{n + 1} \left( {\mathcal A} \right)$.
Consider the following short exact sequence of $\tilde{G} \times G$-equivariant coherent sheaves:
\begin{eqnarray}
0 \rightarrow \Phi ^{\tau} _{n} \left( {\mathcal A} \right) \rightarrow \Phi ^{\tau} _{n} \left( {\mathcal A} \right) \otimes _X {\mathcal O} _X \left( D _{\tau} \right) \rightarrow
{\mathcal Q} _{n + 1} \left( {\mathcal A} \right) \rightarrow 0 \label{def of Q}
\end{eqnarray}

Since ${\mathcal O} _X \left( - D _{\tau} \right)$ annihilates
${\mathcal Q} _{n + 1} \left( {\mathcal A} \right)$, we can regard ${\mathcal Q} _{n + 1} \left( {\mathcal A} \right)$
as a $\tilde{G} \times G$-equivariant vector bundle on $D _{\tau}$. To proceed the proof, we need a lemma.

\begin{lemma}\label{slice}
We have a natural inclusion of $\mathfrak{g} ^{\tau}$-modules
$$\frac{ V \otimes
t ^{N} k[t] + \oplus _{m \in \Z} F ^{\tau} \left( m + N - n\right) \otimes k t ^{- m - 1}} { V \otimes
t ^{N} k[t] + \oplus _{m \in \Z} F ^{\tau} \left( m + N - n\right) \otimes k t ^{- m}} \hookrightarrow {\mathcal Q} _{n + 1} \left( {\mathcal A} \right)
\otimes _X k ( x _{\tau} ).$$
\end{lemma}

\begin{proof}
By Corollary \ref{def of rho}, we have
$$\frac{ V \otimes
t ^{N - 1} k[t] + \oplus _{m \in \Z} F ^{\tau} \left( m + N - n\right) \otimes k t ^{- m - 1}} { V \otimes
t ^{N} k[t] + \oplus _{m \in \Z} F ^{\tau} \left( m + N - n\right) \otimes k t ^{- m}} \cong {\mathcal Q} _{n + 1} \left( {\mathcal A} \right) \otimes _X k ( x _{\tau} ).$$
Since each term is $\left( \triangle ^d ( \tilde{L} ^{\tau} ) \mathbb{G} ^{\tau} _m, L _{\tau}, R _{\tau}, k [t] \right)$-stable by the induced action from $V ^{\tau} _{[t]}$, we have the result by Lemma \ref{uniq stable action}.
\end{proof}
\begin{flushleft}
{\it Continuation of the proof of Proposition \ref{simple makeup}.}
\end{flushleft}
We denote by ${\mathcal S} _{n + 1} ^0 \left( {\mathcal A} \right)$ the $\tilde{G} \times G$-equivariant vector bundle on $\mathbf O _{\tau}$ whose fiber at $x _{\tau}$ is isomorphic to the LHS of the Lemma
\ref{slice} as a $\mathfrak{g} ^{\tau}$-module. Hence, we have an inclusion
of $\tilde{G} \times G$-equivariant vector bundles ${\mathcal S} _{n + 1} ^0 \left( {\mathcal A} \right) \subset
{\mathcal Q} _{n + 1} \left( {\mathcal A} \right) \MID _{\mathbf O _{\tau}}$. We define a $\tilde{G} \times G$-equivariant coherent subsheaf ${\mathcal S} _{n + 1} \left( {\mathcal A} \right)$ of ${\mathcal Q} _{n + 1} \left( {\mathcal A} \right)$ as follows:
\begin{eqnarray}
{\mathcal S} _{n + 1} \left( {\mathcal A} \right) \left( {\mathcal U} \right) := \{ s \in \Gamma ( {\mathcal U}, {\mathcal Q} _{n + 1} \left( {\mathcal A} \right)  ) ; s
\MID _{\mathbf O _{\tau}} \in \Gamma ( {\mathcal U} \cap
\mathbf O _{\tau}, {\mathcal S} _{n + 1} ^0 \left( {\mathcal A}
\right) )\} \label{defofS}\\\nonumber
 \text{ for every Zariski open
set } {\mathcal U} \subset D _{\tau}.
\end{eqnarray}

The preimage ${} ^{\prime} \Phi ^{\tau} _{n + 1} \left( {\mathcal A} \right)$ of ${\mathcal S} _{n + 1} \left( {\mathcal A} \right)$ by $\Phi ^{\tau} _{n} \left( {\mathcal A} \right)
\otimes _X {\mathcal O} _X \left( D _{\tau} \right)
\rightarrow {\mathcal Q} _{n + 1} \left( {\mathcal A} \right)$ is a $\tilde{G} \times G$-equivariant ${\mathcal O}
_X$-module by the assumption (on $\Phi ^{\tau} _{n} \left( {\mathcal A} \right)$) and the construction (of ${\mathcal Q} _{n} \left( {\mathcal A} \right)$
and (\ref{def of Q})).

Since we have ${} ^{\prime} \Phi ^{\tau} _{n + 1} \left( {\mathcal A} \right) \subset {} ^{\prime} \Phi ^{\tau} _{n} \left( {\mathcal A}
\right) \otimes _X {\mathcal O} _X ( D _{\tau} )$, it follows
that ${} ^{\prime} \Phi ^{\tau} _{n + 1}\left( {\mathcal A} \right)$ satisfies the condition 1) of $( \sharp ) _{n + 1}$.

We show that ${} ^{\prime} \Phi ^{\tau} _{n + 1} \left( {\mathcal A} \right)$
is a vector bundle in order to check $( \sharp ) _{n + 1}$. We have a composition
$$\mathrm{Res} : \mathfrak{Coh} ^{\tilde{G} \times G} X ( \Sigma ) \rightarrow \mathfrak{Coh} ^{\tilde{B} \times B} U ^-. T ( \Sigma ). U \stackrel{\cong}{\longrightarrow} \mathfrak{Coh} ^{\tilde{T} \times T} T ( \Sigma )$$
of restriction functors. Here the second equivalence is a consequence of the local structure theorem (Theorem \ref{structure theorem}). This also implies that a $\tilde{G} \times G$-equivariant coherent sheaf on $X \left( \Sigma \right)$ is a vector bundle if, and only if, its restriction to $T \left( \Sigma \right)$ is a vector bundle. Hence, what has to be proved is:

\begin{lemma}
For each $\sigma \in \Sigma$ such that $\tau \in \sigma ( 1 )$,
${} ^{\prime} \Phi ^{\tau} _{n + 1} \left( {\mathcal A} \right) \MID _{T \left( \sigma \right)}$ is a vector bundle.
\end{lemma}

\begin{proof}
${\mathcal S} _{n + 1} \left( {\mathcal A} \right) \MID _{T (
\sigma )}$ is isomorphic to the sheaf obtained by replacing ${\mathcal
U}$ by ${\mathcal U} \cap T ( \sigma )$ and taking tensor products with
${\mathcal O} _{T (\sigma)}$ at the RHS of (\ref{defofS}). Since $T ( \sigma )$ is smooth, we can write $T \left( \sigma
\right) := \Spec R$, where $R = k [t , t _2, \ldots, t _{m}, t _{m + 1} ^{\pm 1},
\ldots, t _r ^{\pm 1}]$. Here we assume that $D _{\tau} \cap T \left(
\sigma \right)$ is defined by $t = 0$ and each of $t, t _2, 
\ldots, t _m$ spans a extremal ray of the dual cone $\sigma ^{\vee}$. By Theorem \ref{Klyachko's Proposition} 1), ${} ^{\prime} \Phi ^{\tau} _{n} \left( {\mathcal A} \right) \MID _{T \left( \sigma \right)}$ is written as $V _0 \otimes R$ by some
$1 \times T$-module $V _0$. We have
$$( {} ^{\prime} \Phi ^{\tau} _{n} \left( {\mathcal A} \right) \otimes _X {\mathcal O} _X ( D _{\tau} ) ) \MID _{T \left( \sigma \right)} \otimes _{{\mathcal O} _{T \left( \sigma \right)}} k ( x _{\tau} ) \cong {\mathcal Q} _{n + 1} \left( {\mathcal A} \right) \otimes _X k ( x _{\tau} ).$$
Hence, we have
$${\mathcal Q} _{n + 1} \left( {\mathcal A} \right) \MID _{T \left( \sigma \right)} \cong V _0 \otimes t ^{- 1} k[t _2, \ldots, t _m , t _{m + 1} ^{\pm 1}, \ldots, t _r ^{\pm 1}]$$
as a $R / t R$-module. For a vector subspace
$V _1 := {\mathcal S} _{n + 1} ^0 \left( {\mathcal A} \right) \otimes _X k ( x _{\tau} )$ of $V _0$, we have
$${\mathcal S} _{n + 1} ^0 \MID _{T ( \sigma )} \cong V _0 \otimes t ^{- 1} k[t _2, \ldots, t _m , t _{m + 1} ^{\pm 1},
\ldots, t _r ^{\pm 1}] \cap V _1 \otimes t ^{- 1} k[t _2 ^{\pm 1},
\ldots, t _r ^{\pm 1}] = V _1 \otimes t ^{- 1} R / R$$
Thus, we have
$${} ^{\prime} \Phi ^{\tau} _{n + 1} \left( {\mathcal A} \right) \MID _{T \left( \sigma
\right)} \cong V _0 \otimes R + V _1 \otimes t ^{- 1 } R \subset V _0
\otimes t ^{- 1}R.$$
Since this module is $R$-free, we obtain the result.
\end{proof}
\begin{flushleft}
{\it Continuation of the proof of Proposition \ref{simple makeup}.}
\end{flushleft}
The
definition of $F ^{\tau} ( \bullet, {}^{\prime} \Phi ^{\tau} _{n + 1} \left(
{\mathcal A} \right) )$ makes sense since ${}^{\prime} \Phi ^{\tau} _{n + 1} \left( {\mathcal A} \right)$ is a $\tilde{G} \times G$-equivariant vector bundle. We have
$$\rho ^{\tau} _{\infty} ( {}^{\prime} \Phi ^{\tau} _{n + 1} \left( {\mathcal A} \right) ) = \left( V \otimes t ^{N} k[t] + \oplus _{m \in \Z} F ^{\tau} \left( m + N - n\right) \otimes k t ^{- m - 1} \subset V ^{\tau} _{[t]} \right)$$
by the comparison with the corresponding extension over $T ( \tau )$ via $\mathrm{Res}$.

In other words, ${}^{\prime} \Phi ^{\tau} _{n + 1} \left( {\mathcal A} \right)$ satisfies 2) and 3) of $( \sharp ) _{n + 1}$. 
Hence, putting $\Phi ^{\tau} _{n + 1} \left( {\mathcal A} \right) := {}^{\prime} \Phi ^{\tau} _{n + 1} \left( {\mathcal A} \right)$ proceeds the induction on $n$.

Now we construct a morphism $\Phi ^{\tau} ( f ) : \Phi ^{\tau} \left( {\mathcal A} \right) \rightarrow \Phi ^{\tau} \left( {\mathcal B} \right)$ from $f : {\mathcal A} \rightarrow {\mathcal B}$. Enlarge $N$ if necessary to assume $F _1
^{\tau} \left( N - 1 \right) = \{ 0 \}$ and $F _1 ^{\tau} \left( - N \right) = V _1$.
We construct a morphism $f ^{\prime} _n : \Phi ^{\tau} _{n} \left( {\mathcal A} \right) \rightarrow \Phi ^{\tau} _{n} \left( {\mathcal B} \right)$ by induction on $n$. For $n = - N$, we put $f ^{\prime} _{- N} := f \otimes \mathrm{id} : V \otimes {\mathcal O} _X ( - N D _{\tau} ) \rightarrow V _1 \otimes {\mathcal O} _X ( - N D _{\tau} )$. Assume that we have $f ^{\prime} _n$. For each morphism $f : {\mathcal A} \rightarrow {\mathcal
B}$, we have the following associated commutative diagram
of compatible $\left( \triangle ^d ( \tilde{L} ^{\tau} ) \mathbb{G} ^{\tau} _m, L _{\tau}, R _{\tau}, k [t] \right)$-modules for all $n \in
\Z$.
$$
\begin{matrix}
V \otimes
t ^{N} k[t] + \oplus _{m \in \Z} F ^{\tau} \left( m + N - n\right) \otimes k t ^{- m - 1} & \hookrightarrow & V ^{\tau} _{[t]} \\
\downarrow & & \downarrow \\
V _1 \otimes
t ^{N} k[t] + \oplus _{m \in \Z} F _1 ^{\tau} \left( m + N - n\right)
\otimes k t ^{- m - 1} & \hookrightarrow & ( V _1 ) ^{\tau} _{[t]}
\end{matrix}
$$

Hence, we have the following commutative diagram arising from the $\mathfrak{g} ^{\tau}$-module homomorphism of Lemma \ref{uniq stable action}.
$$
\begin{matrix}
{\mathcal S} _{n + 1} \left( {\mathcal A} \right) &
\hookrightarrow & {\mathcal Q} _{n + 1} \left( {\mathcal A} \right) & \leftarrow & \Phi ^{\tau} _{n} \left( {\mathcal A} \right) \otimes _X
{\mathcal O} _X \left( D _{\tau} \right) \\
\downarrow & & \downarrow & & \qquad \downarrow f ^{\tau} _{n} \otimes \id\\
{\mathcal S} _{n + 1} \left( {\mathcal B} \right)
\quad &
\hookrightarrow & {\mathcal Q} _{n + 1} \left( {\mathcal B} \right) & \leftarrow &\Phi ^{\tau} _{n} \left( {\mathcal B} \right) \otimes _X
{\mathcal O} _X \left( D _{\tau} \right)
\end{matrix}
$$

Hence, taking preimages yield a morphism $f ^{\tau} _{n + 1} : {}^{\prime} \Phi
^{\tau} _{n + 1} \left( {\mathcal A} \right) \rightarrow
{}^{\prime} \Phi ^{\tau} _{n + 1} \left( {\mathcal B} \right)$ of
$\tilde{G} \times G$-equivariant ${\mathcal O} _X$-modules. Therefore, the induction on $n$ proceeds. Thus, setting $\Phi ^{\tau} ( f ) := f ^{\tau} _N$ yields a morphism $\Phi ^{\tau} _{hom} : \Hom _{{\mathfrak C} \left( \Sigma, 1 \right) _{c}} \left( {\mathcal A}, {\mathcal B} \right) \rightarrow \Hom _{EV \left( \Sigma, 1 \right) _{c}} \left( \Phi ^{\tau} ( {\mathcal A} ), \Phi ^{\tau} ( {\mathcal B} ) \right)$. Since $f : V \rightarrow V _1$ yields a morphism between the identity fibers of $\Phi ^{\tau} ( {\mathcal A} )$ and $\Phi ^{\tau} ( {\mathcal B} )$, the injectivity of $\Phi ^{\tau} _{hom}$ is clear.
\end{proof}

\subsection{Proof of Theorem \ref{Main Theorem}}\label{finish}
By Proposition \ref{last} and Proposition \ref{reverse functor Phi}, both $\Xi$
and $\Phi$ define one-to-one correspondences. By Proposition
\ref{naive functor} and Proposition \ref{reverse functor Phi}, both $\Xi$ and
$\Phi$ are faithful. Hence, we have the following category equivalence:
$$\Xi : EV \left(
\Sigma \right) _{c} \stackrel{\cong}{\rightarrow} {\mathfrak C}
\left( \Sigma \right) _{c}.$$
Next, we want to check what happens if we restrict our attention to $EV \left(
\Sigma \right)$.

\begin{lemma}\label{modulo codim two}
Let ${\mathcal E}, {\mathcal F} \in {\bf Ob} EV \left( \Sigma, 1 \right) _{c}$ and let $f
: {\mathcal E} \rightarrow {\mathcal F}$ be a morphism in $EV \left(
\Sigma, 1 \right)$.
 Then, for each $\tau \in \Sigma ( 1 )$, we have
$$\Xi ( f ) \left( F ^{\tau} ( n, {\mathcal E} )
\right) = \Xi ( f )
\left ( B ( {\mathcal E} ) \right) \cap F ^{\tau} ( n, {\mathcal F}
)\text{ for all } n \in \Z \text{ and}$$
$$\{ \Xi ( f ) ( B \left( {\mathcal E} \right) ), \{ F
^{\tau} \left( n, {\mathcal F} \right) \} _{n \in \Z} \}\text{ forms a
distributive lattice.}$$
\end{lemma}

\begin{proof}
Consider the restriction to $\mathbb A _{\tau}$ as in \S 3.1.1. Since a $\mathbb G ^{\tau} _m$-equivariant vector bundle on $\mathbb A _{\tau}$ splits, we have $$\mathcal F \MID _{\mathbb A _{\tau}} \cong \mathrm{Im} f \MID _{\mathbb A _{\tau}} \oplus \mathrm{coker} f \MID _{\mathbb A _{\tau}}$$
as $\mathbb G ^{\tau} _m$-equivariant vector bundles. Hence, we have a splitting
$$F ( n, \mathcal F ) = \left( F ( n, \mathcal F ) \cap \mathrm{Im} \Xi ( f ) \right) \oplus \left( F ( n, \mathcal F ) \cap \mathrm{coker} \Xi ( f ) \right)$$
for every $n \in \mathbb Z$.
\end{proof}

Let ${\mathcal E}, {\mathcal F} \in {\bf Ob} EV \left( \Sigma \right) _{c} (= {\bf Ob} EV \left( \Sigma \right))$ and let $f
: {\mathcal E} \rightarrow {\mathcal F}$ be a morphism of $EV \left(
\Sigma \right)$. Then $\Xi \left( f \right)$
satisfies the condition (L) of Definition \ref{cat} by Lemma \ref{modulo
codim two}. By definition, $\mathrm{coker}
f$ is a $\tilde{G} \times \tilde{G}$-equivariant vector bundle. Thus, $\Xi \left( f \right)$ also satisfies the condition (R) of Definition \ref{cat} by Lemma \ref{hom of naive functor}. Hence, $\Xi ( f )$ is a morphism of ${\mathfrak C}
\left( \Sigma \right)$. 

We prove its converse. Suppose $\left( V _1, \{ F ^{\tau} _1 \left( \bullet \right)
\} \right), \left( V _2, \{ F ^{\tau} _2 \left( \bullet \right)
\} \right) \in {\bf Ob} {\mathfrak C} \left(
\Sigma, 1 \right)$ and let $f ^{\prime} : \left( V _1, \{ F ^{\tau} _1 \left( \bullet \right)
\} \right) \rightarrow \left( V _2, \{ F ^{\tau} _2 \left( \bullet \right)
\} \right)$
be a morphism in ${\mathfrak C} \left(
\Sigma, 1 \right)$. Then, we have a morphism $\Phi ( f ^{\prime} )$ of $EV \left(
\Sigma, 1 \right) _{c}$. By the definition of morphisms of ${\mathfrak C} \left(
\Sigma, 1 \right)$, we have
$$\mathrm{Coker} f ^{\prime} := \left( \mathrm{coker} f ^{\prime}, \{ F _2 ^{\tau}
\left( \bullet \right) / f ^{\prime} \left( F _1 ^{\tau}
\left( \bullet \right) \right) \} \right) \in
{\bf Ob} {\mathfrak C} \left(
\Sigma, 1 \right).$$
Our construction of $\Phi$ coincides with that of Klyachko's when restricted to $T \left( \Sigma \right)$. Therefore, the condition (R) of Definition \ref{cat} yields
$$\Phi \left( \left( V _1, \{ F ^{\tau} _1 \left( \bullet \right)
\} \right) \right) \MID _{T \left( \sigma \right)} \cong \Phi \left( \left( f ^{\prime} ( V _1 ), \{ f ^{\prime} ( F _1 ^{\tau}
\left( \bullet \right) ) \} \right) \right) \MID _{T \left( \sigma \right)} \oplus \Phi \left( \mathrm{Coker} f ^{\prime} \right) \MID _{T \left( \sigma \right)}$$
for each $\sigma \in \Sigma$. Thus, both $\mathrm{coker} \Phi ( f ^{\prime} ) \MID _{T \left( \Sigma \right)}$ and $\mathrm{Im} \Phi ( f ^{\prime} ) \MID _{T \left( \Sigma \right)}$ are vector bundles. Hence, $( G \times G ) T \left( \Sigma \right) = X \left( \Sigma \right)$ asserts that both $\mathrm{coker} \Phi ( f ^{\prime} )$ and $\mathrm{Im} \Phi ( f ^{\prime} ) $ are vector bundles. By the similar argument, $\mathrm{ker} \Phi ( f ^{\prime} )$ is also a vector bundle. Therefore, $\Phi \left( f ^{\prime} \right)$ is a morphism of $EV \left( \Sigma \right)$. By the arguments before Proposition
\ref{reverse functor Phi}, we obtain the category equivalence

$$\Xi : EV \left( \Sigma \right) \rightarrow {\mathfrak C} \left(
\Sigma \right)$$
from the equivalence of the ambient categories.

\section{Consequences of the main theorem}\label{sectioncor}

\subsection{Comparison with Klyachko's category}\label{compare}
In this subsection, we compare our category ${\mathfrak C} \left(
\Sigma \right)$ to that of Klyachko's.

We have a natural $\tilde{T} \times \tilde{T}$-equivariant embedding $T \left( \Sigma \right) \hookrightarrow X \left( \Sigma \right)$. We denote the category of $\tilde{T} \times \tilde{T}$-equivariant vector bundles on $T \left( \Sigma \right)$ by $EV \left( \Sigma \right) _T$ (to distinguish it from $EV \left( \Sigma \right)$, the category of $\tilde{G} \times \tilde{G}$-equivariant vector bundles on $X \left( \Sigma \right)$). Then, pullback defines a functor
$$\mathrm{rest} _{\Sigma} : EV \left( \Sigma \right) \rightarrow EV \left( \Sigma \right) _T.$$
Let ${\mathfrak C} \left( \Sigma \right) _T$ be the category
${\mathfrak C} \left( \Sigma \right)$ in \S \ref{NoteCat} which is obtained by
replacing $G$ with $T$. We define a functor restricting the
$\tilde{G}$-action to the $\tilde{T}$-action by
$$\mathrm{rest} _{\Sigma} ^{\prime} : {\mathfrak C} \left( \Sigma \right) \rightarrow {\mathfrak C} \left( \Sigma \right) _T.$$
Then, we have $\Xi \circ \mathrm{rest} _{\Sigma} =  \mathrm{rest}
_{\Sigma} ^{\prime} \circ \Xi$ since the construction of $\Xi$ in \S \ref{CatEV} is in terms of $T$-orbits
of $e$. Define a full-subcategory ${\mathfrak C} \left( \Sigma \right) _T ^-$ of ${\mathfrak C} \left( \Sigma \right) _T$ as follows:
$${\bf Ob} {\mathfrak C} \left( \Sigma \right) _T ^- := \{ \left( V, \{ F ^{\tau} \left( \bullet \right) \} \right) \in {\bf Ob} {\mathfrak C} \left( \Sigma \right) _T ; V \text{ is trivial as a } \tilde{T} _s\text{-module} \}.$$
Here we regard $\tilde{T} \cap \tilde{G} _s$ as a subgroup of $\tilde{T} \times \tilde{Z} ( T )$.

In this case, we can naturally regard $V$ as a $\tilde{G} _s$-module with trivial $\tilde{G} _s$-action. Since the transversality condition of Definition \ref{adm} is a void condition for a direct sum of trivial $\tilde{G} _s$-modules, we have the following section of $\mathrm{rest} ^{\prime} _{\Sigma}$.

$$\mathrm{sect} _{\Sigma} : {\bf Ob} {\mathfrak C} \left( \Sigma \right) _T ^- \ni \left( V, \{ F ^{\tau} \left( \bullet \right) \} \right) \mapsto \left( V, \{ F ^{\tau} \left( \bullet \right) \} \right) \in {\bf Ob} {\mathfrak C} \left( \Sigma \right).$$

$\mathrm{sect} _{\Sigma}$ naturally gives rise to a functor which we denote by the same letter. Its inverse functor is $\mathrm{rest} _{\Sigma} ^{\prime}$.

We review the affine local description of Klyachko. Notice that the left $T$-action and (the inverse of) the right $T$-action on $T \left( \Sigma \right)$ are the same since $T$ is commutative. As a result, a $T \times T$-equivariant vector bundle on $T \left( \Sigma \right)$ is a vector bundle on $T \left( \Sigma \right)$ with two commutative $T$-equivariant structures.

\begin{corollary}[cf. Klyachko \cite{Kl} 6.1.5 2)]\label{split cor of Klyachko}
For each $\sigma \in \Sigma ( r )$, every $T \times T$-equivariant vector bundle on $T \left( \sigma \right)$ splits into a direct sum of $T \times T$-equivariant line bundles. $\Box$
\end{corollary}

Now we state an analogous statement of Corollary \ref{split cor of Klyachko} in the wonderful setting, which the author does not know any non-equivariant counterpart.

\begin{theorem}\label{split via fiber}
Let $G$ be an adjoint semisimple group and let $X = X _0 = X \left( \Sigma _0 \right)$ be its wonderful compactification as defined in \cite{DP}. Then, a $\tilde{G} \times \tilde{G}$-equivariant vector bundle ${\mathcal E}$ on $X$ splits into a direct sum of $\tilde{G} \times \tilde{G}$-equivariant line bundles if, and only if, ${\mathcal E} \otimes _X k ( e )$ is trivial as a $\tilde{G}$-module.
\end{theorem}

\begin{proof}
We have $\left( V, \{ F ^{\tau}
\left( \bullet \right) \} _{ \tau \in \Sigma ( 1 )} \right) :=  \Xi \left(
{\mathcal E} \right) \in  {\bf Ob}{\mathfrak C} \left( \Sigma _0 \right) ^- _T$. By assumption, we have $\tilde{G} _s \cong \tilde{G}$. Hence, $V$ is a direct sum of trivial $\tilde{G} \times 1$-modules for each $\left( V, \{ F ^{\tau} \left( \bullet \right) \} _{ \tau \in \Sigma ( 1 )} \right) \in {\bf Ob}{\mathfrak C} \left( \Sigma _0 \right) ^- _T$. Let $V = \bigoplus _{\chi \in Z ( G ) ^{\vee}} V _{\chi}$ be the isotypical decomposition of $V$ with respect to the $1 \times Z ( G )$-action. By Corollary \ref{isotypical decomposition for C}, we have
$$\left( V, \{ F ^{\tau} \left( \bullet \right) \} _{ \tau \in \Sigma _0 ( 1 )} \right) = \bigoplus _{\chi \in \tilde{Z} ( G ) ^{\vee}} \left( V _{\chi}, \{ F ^{\tau} \left( \bullet \right) _{\chi} \} _{ \tau \in \Sigma _0 ( 1 )} \right).$$
$\Sigma _0$ consists of a unique $r$-dimensional cone $\sigma$ and its faces. For each $\chi \in \tilde{Z} ( G ) ^{\vee}$, we have a basis $B^{\sigma} _{\chi}$ of $V _{\chi}$ such that $F ^{\tau} \left( n \right) _{\chi}$ is spanned by a subset of $B ^{\sigma} _{\chi}$ for every $\tau \in \sigma ( 1 )$ and every $n \in \Z$. Thus, $\left( V, \{ F ^{\tau} \left( \bullet \right) \} _{ \tau \in \Sigma _0 ( 1 )} \right) \in {\bf Ob}{\mathfrak C} \left( \Sigma _0 \right) ^- _T$ splits into a direct sum of one-dimensional objects. Sending by $\Phi \circ \mathrm{sect} _{\Sigma _0}$, we obtain the result.
\end{proof}

We formulate corollaries of the above result by using the following well-known result.

\begin{lemma}[Minimal dimension of modules cf. \cite{Ens}]
Let $G$ be an adjoint simple group. Then, the minimal dimension $d _G$ of a nontrivial representation of a simply connected cover $\tilde{G}$ of $G$ is the following: 1) $d _{A _r} = r + 1$, 2) $d _{B _r} = 2 r + 1$, 3) $d _{C _r} = 2 r$, 4) $d _{D _r} = 2 r$, 5) $d _{E _6} = 27$, 6) $d _{E _7} = 56$, 7) $d _{E _8} = 248$, 8) $d _{F _4} = 26$, 9) $d _{G _2} = 7$.
\end{lemma}

Together with Theorem \ref{split via fiber}, we obtain the following

\begin{corollary}
Let $G$ be an adjoint simple group and let $X = X _0$ be its wonderful compactification.
Then, every $\tilde{G} \times \tilde{G}$-equivariant vector bundle of rank less than $d _G$ splits into a direct sum of line bundles.
\end{corollary}

Since $d _{G}$ is always greater than the rank $r = \mathrm{rk} G$ of a
simple group $G$, we have the following.

\begin{corollary}\label{split via rank-statement}
Let $G$ be an adjoint simple group and let $X = X _0$ be its wonderful compactification.
Then, every $\tilde{G} \times \tilde{G}$-equivariant vector bundle of rank less than or equals to $\mathrm{rk} G$ splits into a direct sum of line bundles.
\end{corollary}

\subsection{Kostant's Problem}\label{canonical extension}
The first part of this subsection is independent of the other parts of this paper. A general reference for the material in this subsection is the book of Borel-Wallach \cite{BW}.

Kostant \cite{Ko} raised the question of existence of a canonical extension
of an equivariant vector bundle on a (complexified) symmetric
space to its wonderful compactification. His question is connected to the
existence of a vector bundle corresponding to the asymptotic behavior of the minimal
$K$-type of a given unitary representation. More precisely, he seeks
an algebraic framework to handle the celebrated Casselman theorem (cf. \cite{BW} Chapter X Theorem 2.4) about the $\mathfrak{n} _0$-coinvariants of an unitary
representation using the boundary behavior of equivariant vector
bundles. In our setting, his general
conjecture is as follows.

\begin{conjecture}[Kostant's Problem]
Let $G _0$ be a real reductive linear Lie group, and $K _0$ its maximal compact subgroup. We fix an Iwasawa decomposition $G _0 = K _0 A _0 N _0$ of $G _0$. Here $A _0$ is an abelian group which normalizes $N _0$. We define $G$, $K$,... to be the complexification of $G _0$, $K _0$,... respectively. Let $e \in G / K$ be the point corresponding to $[K] \in G / K$. Then, for each irreducible
$K$-module $V$ which appears as the minimal $K _0$-type of an unitary
representation of $G _0$, we have a $G$-equivariant vector bundle ${\mathcal E} _V$ on the wonderful compactification $Y$ of $G / K$ with the following properties:
\begin{enumerate}
\item There exists an isomorphism ${\mathcal E} _V \MID _{G / K} \cong G \times _K V$ of $G$-equivariant vector bundles;
\item For every $v \in V \cong {\mathcal E} _V \otimes _{{\mathcal O} _Y} k ( e )$, and every one-parameter subgroup $\tau : \mathbb{G} _m \rightarrow A$, the limit value $\lim_{t \rightarrow \infty} \tau ( t ) v$ exists in the total space $V \left( {\mathcal E} _V \right)$ of ${\mathcal E} _V$;
\item For every $G$-equivariant vector bundle ${\mathcal E}$ on $Y$ with the above two properties, we have a $G$-equivariant embedding ${\mathcal E} \hookrightarrow {\mathcal E} _V$.
\end{enumerate}
\end{conjecture}

\begin{remark}
The wonderful compactification of a symmetric space is an algebraic analogue of the Oshima compactification \cite{Os} in the real-analytic setting. Sato \cite{SS} described the analogous complex-analytic compactification in the case of adjoint semisimple groups, which coincides with the wonderful compactification in the sense of De Concini-Procesi \cite{DP}.
\end{remark}

In the rest of this subsection, we assume that $X ( = X
_0 = X \left( \Sigma _0 \right) )$ is the wonderful compactification of an adjoint
semisimple group over $\C$ and use the notation and terminology introduced in \S \ref{NoteAlg} and \S \ref{NoteVar}. 

We define the canonical extension of a $G \times G$-equivariant vector bundle
corresponding to $G$-module $V$ on $G \cong ( G \times G ) / \triangle ^d
( G )$ as follows. By Lemma \ref{Basic Situation}, $V \otimes
{\mathcal O} _G$ is the only $G \times G$-equivariant vector bundle on
$G$ whose identity fiber is isomorphic to $V$ (as a $G$-module). Let $V = \oplus _{\lambda \in X ^* \left( T \right)} V ^{\lambda}$ be the $T$-isotypical decomposition of $V$.

Then, we define the canonical filtration of $V$ with respect to $\varpi \in \Sigma _0 ( 1 )$ as follows:

$$F ^{\varpi} _{can} \left( n, V \right):= \bigoplus _{\left< \varpi, \lambda \right> \ge 2 n} V ^{\lambda} \subset V \text{ for each } n \in \Z.$$

We regard $V$ as a $\tilde{G} \times \tilde{Z} ( G )$-module via natural surjection. 

\begin{lemma}
For each $G$-module $V$, $\left( V, \{ F ^{\varpi} _{can} \left( \bullet, V \right) \} _{\varpi \in \Sigma _0 ( 1 )} \right)$ is an object of ${\mathfrak C} \left( \Sigma \right)$.
\end{lemma}

\begin{proof}
Since $F ^{\varpi} _{can} \left( n, V \right)$ is a direct sum of $T$-isotypical components for every $\varpi \in \Sigma _0 ( 1 )$, $\{ F ^{\varpi} \left( n \right) \} _{\varpi \in \Sigma _0 ( 1 ), n \in \Z}$ clearly forms a distributive lattice. Moreover, it is a $\mathfrak{p} ^{\varpi}$-stable decreasing filtration since $\mathfrak{p} ^{\varpi} V ^{\lambda} \subset \oplus _{\left< \varpi, \mu \right> \ge 0} V ^{\lambda + \mu}$. Similarly, we have
$$e _{\alpha} : \bigoplus _{\left< \varpi, \lambda \right> \ge 2 n} V ^{\lambda} \rightarrow \bigoplus _{\left< \varpi, \lambda \right> \ge 2 n} V ^{\lambda + \alpha} \subset \bigoplus _{\left< \varpi, \lambda \right> \ge 2 n + \left< \varpi, \alpha \right>} V ^{\lambda} \subset \bigoplus _{\left< \varpi, \lambda \right> \ge 2 n + 2 \left< \varpi, \alpha \right>} V ^{\lambda}$$
for every negative root $\alpha$. Hence, $F ^{\varpi} _{can} \left( \bullet, V \right)$ is a $\varpi$-transversal filtration.
\end{proof}

Now we define the canonical extension ${\mathcal E} _V$ of a $G$-module $V$ by
$${\mathcal E} _V := \Phi \left( ( V, \{ F ^{\varpi} _{can} \left( \bullet, V \right) \} _{\varpi \in \Sigma _0 ( 1 )} ) \right).$$
This is a $G \times G$-equivariant vector bundle on $X$ which is isomorphic to $V \otimes {\mathcal O} _G$ when restricted to $G$. Since $G \times G$ is the complexification of a complex group $G$
regarded as a real group, $G$(=$\left( G \times G \right)/ \triangle ^d \left( G \right)$) is a symmetric space. Here the group is $G \times G$ and its
(complexified) Cartan involution $\vartheta$ swaps its first and
second factor. Iwasawa's abelian subgroup is defined as the
exponential of the maximal abelian
subalgebra in the set of semisimple elements of $-1$-eigenpart of the induced $\vartheta$-action on
$\mathfrak{g} \oplus \mathfrak{g}$. Thus, in our case, it is the
anti-diagonal embedding of the torus $T$ into $G \times G$ defined by $t
\mapsto ( t, t ^{- 1} )$. Let $a : T
\curvearrowright X$ be the action coming from $T \times G \ni ( t, g )
\mapsto t. g. t \in G$. By the above arguments, this is the
(complexified) action of Iwasawa's abelian subgroup.

Now we present an answer to Kostant's problem in this case.

\begin{theorem}\label{Answer to Kostant}
For every $G$-module $V$, ${\mathcal E} _V$ satisfies the following
properties:
\begin{enumerate}
\item ${\mathcal E} _V \otimes _X k ( e ) \cong V \boxtimes k$ as a
$\tilde{G} \times Z \left( G \right)$-module;
\item For every $v \in {\mathcal E} _V \otimes _X k ( e )$, there exist a
limit value $\lim _{t
\rightarrow 0} a ( t ) v$ in $V \left( {\mathcal E} _V \right)$, the total space of ${\mathcal E} _V$;
\item Let ${\mathcal E}$ be another $G \times G$-equivariant vector bundle with the above two properties. Then, we have a $G \times G$-equivariant inclusion ${\mathcal E} \hookrightarrow {\mathcal E} _V$. In other word, ${\mathcal E} _V$ is maximal among the $G \times G$-equivariant vector bundles with the
above two properties.
\end{enumerate}
\end{theorem}

\begin{proof}
1) follows from the definition of ${\mathcal E} _V$. Since the $a$-action preserves $T \left( \Sigma _0 \right)$, we restrict our attention to $T \left( \Sigma _0 \right)$. Let $\varpi \in \Sigma _0 ( 1 )$. $\mathbb{G} _m$ acts on ${\mathcal O} _X \left( D _{\varpi} \right) \otimes _X k ( x _{\varpi} )$ via $a \circ \varpi$-action by degree $- 2$. For each weight $\lambda$ of $V$, $V ^{\lambda} \otimes {\mathcal O} _X \MID _{T \left( \Sigma _0 \right)}$ is a $T \times T$-equivariant vector subbundle of $V \otimes {\mathcal O} _X \MID _{T \left( \Sigma _0 \right)}$. $\mathbb{G} _m$ acts on $V ^{\lambda} \otimes {\mathcal O} _X \MID _{T \left( \Sigma _0 \right)} \otimes _{{\mathcal O} _{T \left( \Sigma _0 \right)}} k ( x _{\varpi} )$ via $a \circ \varpi$-action by degree $\left< \varpi, \lambda \right>$. Hence, the condition of the convergence with respect to the $\varpi$-direction ($s \rightarrow 0$ in $s \in \mathbb{A} ^1 \supset \mathbb{G} _m$) is given as
\begin{eqnarray}
V ^{\lambda} \cap F ^{\varpi} \left( \lfloor \frac{\left< \varpi, \lambda \right>}{2} \rfloor + 1 \right) = \{ 0 \}\label{convergence}
\end{eqnarray}
by Corollary \ref{def of rho} and Theorem \ref{Klyachko's Proposition}. Hence, if the filtration $F ^{\varpi} \left( \bullet
\right)$ of $V$ satisfies
\begin{eqnarray}
V ^{\lambda} \subset F ^{\varpi} \left( \lfloor \frac{\left< \varpi, \lambda \right>}{2} \rfloor \right)\label{maximality}
\end{eqnarray}
and (\ref{convergence}) for every $\lambda \in X ^* ( T )$, then the corresponding vector bundle ${\mathcal E}
\MID _{T ( \varpi )}$ is a maximal $T \times T$-equivariant vector bundle on $T ( \varpi )$ such that:
\begin{itemize}
\item[a.] ${\mathcal E} \MID _{T  ( \varpi )}$ is contained in $V \otimes {\mathcal O} _X \left( D \right) \MID _{T  ( \varpi )}$ for a sufficiently large equivariant divisor $D$ (see Corollary \ref{sandwich});
\item[b.] Limit value $\lim _{t \rightarrow 0} ( \varpi \times (
\varpi ) ^{- 1} ) ( t ) v$ exists in $V \left( {\mathcal E} \MID _{T  ( \varpi )} \right)$ for every $v \in {\mathcal E} \MID _{T  ( \varpi )} \otimes _{{\mathcal O} _{T  ( \varpi )}} k ( e )$.
\end{itemize}
We can check (\ref{convergence}) and (\ref{maximality}) $T$-isotypical componentwise because $F ^{\varpi} ( n )$ is $T$-stable for each $n \in \Z$. $F ^{\xi} _{can} \left( \bullet, V \right)$ satisfies
(\ref{convergence}) and (\ref{maximality}) for every fundamental
coweight $\xi$. Therefore, ${\mathcal E} _V \MID _{T (
\varpi )}$ is maximal as a $T$-equivariant vector subbundle of $V \otimes {\mathcal O} _X \left( D \right) \MID _{T  ( \varpi )}$ which satisfies
the conditions a) and b). Thus, every element of ${\mathcal E} _V \MID
_{T \left( \Sigma _0 \right)} \otimes _{{\mathcal O} _{T \left( \Sigma _0 \right)}} k ( e )$ converges when $t \in
T$ approaches to zero in an arbitrary way. This yields 2). For every $G \times G$-equivariant vector bundle ${\mathcal E}$ such that ${\mathcal E} \MID _{T \left( \varpi \right)}$ satisfies the conditions a) and b), we have a $T \times T$-equivariant embedding ${\mathcal E} \MID _{T \left( \varpi \right)} \subset {\mathcal E} _V \MID _{T \left( \varpi \right)}$. By Theorem \ref{Main Theorem}, we can recover ${\mathcal E}$ from the condition a) and a family of filtrations $\{ F ^{\xi} \left( \bullet, {\mathcal E} \right) \} _{\xi \in \Sigma _0 ( 1 )}$ as vector subspaces of $V$. As a result, we have a $G \times G$-equivariant embedding ${\mathcal E} \hookrightarrow {\mathcal E}
_V$.
\end{proof}

{\bf Acknowledgement}
The original version of this paper was written under the supervision of Professor Hisayosi Matumoto. Professor Michel Brion gave me valuable comments and suggestions on many versions of this paper. It is no doubt that this paper does not exist without them. The author wants to express deep and sincere gratitude to both of them. The author also wants to express his thanks to Professors Valery Alexeev, Toshiyuki Katsura, Bertram Kostant, Friedrich Knop, Takayuki Oda, Toshio Oshima, Jiro
Sekiguchi, Tonny A. Springer, Toshiyuki Tanisaki, T\^ohru Uzawa, Drs. Shuusi Harasita, Satoshi Mochizuki, Takao Nishikawa, Kenei Suzuki, Alexis Tchoudjem, and Mr. Yugi Kato.

\end{document}